\documentclass[10pt,letterpaper]{amsart}
\usepackage{charter}
\usepackage[OT2,OT1]{fontenc}
\usepackage[latin9]{inputenc}
\usepackage{array}
\usepackage{units}
\usepackage{textcomp}
\usepackage{mathrsfs}
\usepackage{bm}
\usepackage{multirow}
\usepackage{amstext}
\usepackage{amsthm}
\usepackage{amssymb}
\usepackage{cancel}
\usepackage{stmaryrd}
\usepackage[all]{xy}

\makeatletter

\pdfpageheight\paperheight
\pdfpagewidth\paperwidth

\providecommand{\tabularnewline}{\\}

\numberwithin{equation}{section}
\numberwithin{figure}{section}
\theoremstyle{plain}
\newtheorem{thm}{\protect\theoremname}[section]
\theoremstyle{plain}
\newtheorem{lem}[thm]{\protect\lemmaname}
\theoremstyle{plain}
\newtheorem{cor}[thm]{\protect\corollaryname}
\theoremstyle{plain}
\newtheorem{prop}[thm]{\protect\propositionname}
\theoremstyle{definition}
\newtheorem{defn}[thm]{\protect\definitionname}
\theoremstyle{remark}
\newtheorem{rem}[thm]{\protect\remarkname}

\usepackage{amsfonts}
\usepackage{mathrsfs}
\usepackage{graphicx}
\usepackage{amscd}
\usepackage{enumerate}

\newcommand{\cyr}{
\renewcommand\rmdefault{wncyr} \renewcommand\sfdefault{wncyss} \renewcommand\encodingdefault{OT2} \normalfont
\selectfont
}
\DeclareTextFontCommand{\textcyr}{\cyr}    

   \def\@settitle
{
\begin{center} \baselineskip14\p@\relax \LARGE\@title \end{center}
}

\usepackage[normalem]{ulem}
\usepackage[mathscr]{eucal}
\numberwithin{equation}{section}

\usepackage{verbatim}
\usepackage{amscd}
\usepackage{bm}\usepackage[mathscr]{eucal}
\usepackage[all]{xy}
\usepackage{enumerate}
\usepackage{color}
\usepackage{colortbl}
\usepackage{url}

\input xyimport.tex

\title{{\bf  Key varieties for prime $\mQ$-Fano threefolds \\ defined by Jordan algebras of cubic forms. \\ Part II}}

\author{Hiromichi Takagi}

\address{Department of Mathematics, Gakushuin University, 
Mejiro, Toshima-ku, Tokyo 171-8588, Japan}
\email{hiromici@math.gakushuin.ac.jp}

\newcommand{\sO}{\mathcal{O}}

\newcommand{\mA}{\mathbb{A}}
\newcommand{\mC}{\mathbb{C}}

\newcommand{\mP}{\mathbb{P}}
\newcommand{\mQ}{\mathbb{Q}}

\numberwithin{equation}{section}

 \newcounter{myparagraph}[subsection]

\makeatother

\providecommand{\corollaryname}{Corollary}
\providecommand{\definitionname}{Definition}
\providecommand{\lemmaname}{Lemma}
\providecommand{\propositionname}{Proposition}
\providecommand{\remarkname}{Remark}
\providecommand{\theoremname}{Theorem}

\begin{document}
\maketitle 
\begin{abstract}
Subsequent to the previous paper \cite{Tak5}, we are concerned with
the classification of complex prime $\mathbb{Q}$-Fano $3$-folds
of anti-canonical codimension 4 which are produced, as weighted complete
intersections of appropriate weighted projectivizations of certain
affine varieties related with $\mathbb{P}^{2}\times\mathbb{P}^{2}$-fibrations.
Such affine varieties or their appropriate weighted projectivizations
(possibly allowing some coordinates have weights $0$) are called
key varieties for prime $\mathbb{Q}$-Fano 3-folds. 

The purpose of this paper is to give new constructions of a $14$-dimensional
affine variety $\Upsilon_{\mathbb{A}}^{14}$ and a $15$-dimensional
affine variety $\Pi_{\mathbb{A}}^{15}$ related with $\mathbb{P}^{2}\times\mathbb{P}^{2}$-fibration,
which were constructed and were shown to be key varieties in the papers
\cite{Tak3,Tak4} and \cite{Tay}. It is well-known that the affine
cone of the Segre embedded $\mathbb{P}^{2}\times\mathbb{P}^{2}$ is
defined as the null loci of the so called $\sharp$-mapping of a 9-dimensional
nondegenerate quadratic Jordan algebra $J$ of a cubic form. Inspired
with this fact, we construct $\Upsilon_{\mathbb{A}}^{14}$ and $\Pi_{\mathbb{A}}^{15}$
in the same way coordinatizing $J$ with 9 and 10 parameters, respectively.
The coordinatization with $9$ parameters is derived by using a fixed
primitive idempotent, and the associated Peirce decomposition. The
coordinatization with $10$ parameters is derived from the construction
of a quadratic Jordan subalgebra generated by $\sharp$-products of
two elements due to Petersson \cite{Pe1}. 
\end{abstract}

\maketitle
\markboth{$\mQ$-Fano 3-fold and Jordan algebra }{Hiromichi Takagi}
{\small{}{\tableofcontents{}}}{\small\par}

2020\textit{ Mathematics subject classification}: 14J45, 14E30, 17C10,
17C50

\textit{Key words and phrases}: $\mQ$-Fano $3$-fold, key variety,
$\mP^{2}\times\mP^{2}$-fibration, Jordan algebra.

\section{\textbf{Introduction\label{sec:Introduction}}}

\subsection{Background for classification of $\mQ$-Fano 3-folds\label{subsec:Background-for-classification}}

A complex projective variety is called a \textit{$\mQ$-Fano variety}
if it is a normal variety with only terminal singularities, and its
anti-canonical divisor is ample. The classification of $\mQ$-Fano
3-folds is one of the central problems in Mori theory for projective
3-folds. Although the classification is far from completion, several
systematic classification results have been obtained so far. In the
online database \cite{GRDB}, a huge table of candidates of $\mQ$-Fano
3-folds are given. One kind of the classification results is that
with respect to the codimension ${\rm {ac}_{X}}$ of a $\mQ$-Fano
3-fold $X$ in the weighted projective space of the minimal dimension
determined by the anti-canonical graded ring of $X$. We call ${\rm {ac}_{X}}$
\textit{the anti-canonical codimension}, which can be considered as
a generalization of the genus of a $\mQ$-Fano 3-fold. 

Subsequent to our previous paper \cite{Tak5}, we are concerned with
the classification of prime $\mQ$-Fano 3-folds with ${\rm {ac}_{X}=4}$,
where we recall that a \textit{prime} $\mQ$-Fano 3-fold is one whose
anticanonical divisor generates the group of numerical equivalence
classes of $\mQ$-Cartier Weil divisors. In \cite{CD1,CD2}, Coughlan
and Ducat construct systematically prime $\mQ$-Fano 3-folds with
${\rm {ac}_{X}=4}$ as weighted complete intersections in weighted
projectivizations of the cluster varieties of types $C_{2}$ and $G_{2}^{(4)}$,
which are related with $\mP^{2}\times\mP^{2}$- and $\mP^{1}\times\mP^{1}\times\mP^{1}$-fibrations,
respectively. This type of construction is modeled after the classification
of smooth prime Fano 3-folds by Gushel\'\ \cite{Gu1,Gu2} and Mukai
\cite{Mu2}. The idea of them is describing $\mQ$-Fano 3-folds $X$
in appropriate projective varieties $\Sigma$ of larger dimensions
and usually with larger symmetries; in the case of \cite{Mu2}, $\Sigma$
are rational homogeneous spaces, and in the case of \cite{CD1,CD2},
$\Sigma$ are several weighted projectivizations of the cluster varieties.
Such $\Sigma$ or the affine cone of $\Sigma$ are informally called
\textit{key variety} as in the title of this paper. 

In our previous paper \cite{Tak5}, we construct a $13$-dimensional
affine variety $\mathscr{H}_{\mA}^{13}$ and show that a weighted
projectivization of $\mathscr{H}_{\mA}^{13}$ produces as weighted
complete intersections $\mQ$-Fano 3-folds of genus 3 with only three
$\nicefrac{1}{2}(1,1,1)$-singularities belonging to the class No.5.4
in \cite{Tak1}. Moreover, we show that the cluster variety of type
$C_{2}$ is a subvariety of $\mathscr{H}_{\mA}^{13}$, hence the prime
$\mQ$-Fano 3-folds constructed by \cite{CD1,CD2} or more general
ones are obtained from weighted projectivizations of $\mathscr{H}_{\mA}^{13}$.

The purpose of this paper is to give new constructions of the two
affine varieties $\Upsilon_{\mA}^{14}$ and $\Pi_{\mA}^{15}$ related
with $\mP^{2}\times\mP^{2}$-fibration, which were constructed and
were shown to be key varieties for prime $\mQ$-Fano 3-folds $X$
with ${\rm {ac}_{X}=4}$ in the papers \cite{Tak3,Tak4} and \cite{Tay}.
Our new perspective for the constructions is the theory of the quadratic
Jordan algebra of a cubic form.

\subsection{Quadratic Jordan algebra of a cubic form and its coordinatizations\label{subsec:Jordan-algebra-of}}

Let $J$ be a 9-dimensional nondegenerate quadratic Jordan algebra
of a cubic form over a field $\mathsf{k}$ with reasonable assumptions
(for details, we refer to \cite[Sect.2, Assumption (a)--(c)]{Tak5},
which originally come from \cite{R}). The quadratic Jordan algebra
$J$ is endowed with a quadratic mapping called $\sharp$-mapping,
and a bi-linear product called $\sharp$-product defined by the $\sharp$-mapping.
In ibid., we have obtained a coordinatization of $J$ with 8 parameters
using three complementary primitive idempotent and the associated
Peirce decomposition. In this paper, we obtain two new coordinatizations
of $J$.

To obtain the first coordinatization, we use one primitive idempotent
$\upsilon$ and the associated Peirce decomposition of $J$: 

\[
J=J_{\nicefrac{1}{2}}\oplus J_{0}\oplus J_{1},
\]
where $J_{1}=\mathsf{k\upsilon}$ and $J_{0}\simeq J_{\nicefrac{1}{2}}\simeq\mathsf{k}^{4}$.
Assuming $\mathsf{k}$ is algebraically closed and taking a suitable
basis of each of the subspaces $J_{i}\,(i=\nicefrac{1}{2},0,1)$,
we show that the $\sharp$-products between elements of these basis
can be described by introducing $9$ parameters indexed as $z,s_{i},t_{i}\,(1\leq i\leq4)$
(Subsection \ref{subsec:Coordination-of-9param}). 

The second coordinatization is based on \cite[Prop.6.6]{Pe1}, in
which Petersson gives a Jordan subalgebra of $J$ generated by $\sharp$-products
of two elements. We show that such a subalgebra coincides with $J$
under a generality condition of the two elements, and derive a coordinatization
of $J$ with 8 parameters (Subsection \ref{subsec:Coordination-of-Perterson}).
Further, assuming the characteristic of $\mathsf{k}$ is not $3$
and changing the bases of $J$ and the parameters, we obtain another
coordinatization of $J$, where the number of parameters becomes $10$
due to increase of the choices of bases (Subsection \ref{subsec:Change-of-base}).
We denote by 
\begin{equation}
t_{1},t_{2},t_{123},t_{124},t_{125},t_{126},t_{135},t_{136},t_{245},t_{246}\label{eq:tenparam}
\end{equation}
the $10$ parameters. 

\vspace{5pt}

Hereafter in Section \ref{sec:Introduction}, we assume for simplicity
that $\mathsf{k}=\mC$, the complex number field.

\subsection{Affine varieties $\Upsilon_{\mA}^{14}$ and $\Pi_{\mA}^{15}$\label{subsec:Affine-varieties-UpPi}}

It is well-known that the affine cone $C(\mP^{2}\times\mP^{2})$ of
(the Segre embedded) $\mP^{2}\times\mP^{2}$ is defined as the null
loci of the $\sharp$-mapping of $J$. Inspired with this fact, we
have constructed the affine variety $\mathscr{H}_{\mA}^{13}$ in \cite{Tak5}
in the same way as $C(\mP^{2}\times\mP^{2})$ coordinatizing $J$
with $8$ parameters. In this paper, we construct the affine varieties
$\Upsilon_{\mA}^{14}$ and $\Pi_{\mA}^{15}$ in the same way using
the two coordinatizations of $J$ mentioned in Subsection \ref{subsec:Jordan-algebra-of}. 

Let $\mA_{\Upsilon}$ be the $18$-dimensional affine space containing
$J\simeq\mA^{9}$ with $z,s_{i},t_{i}$ as additional $9$ coordinates,
and $\mA_{\Pi}$ the $19$-dimensional affine space containing $J\simeq\mA^{9}$
with (\ref{eq:tenparam}) as additional $10$ coordinates. Then we
define $\Upsilon_{\mA}^{14}\subset\mA_{\Upsilon}^{18}$ and $\Pi_{\mA}^{15}\subset\mA_{\Pi}^{19}$
as the null loci of the $\sharp$-mapping. We refer to Corollary \ref{cor:Upsilonshar}
and Definition \ref{def:Upsilon} for the explicit equations of $\Upsilon_{\mA}^{14}$,
and to Proposition \ref{prop:The--mapping-Pi} and Definition \ref{def:Pi15}
for the explicit equations of $\Pi_{\mA}^{15}$. 

Here we summarize properties of $\Upsilon_{\mA}^{14}$ and $\Pi_{\mA}^{15}$
as follows, where we include some properties for $\Upsilon_{\mA}^{14}$
obtained by \cite{Tay} for convenience:
\begin{thm}
\label{thm:UpPiKeyMain} The affine variety $\Upsilon_{\mA}^{14}$
and $\Pi_{\mA}^{15}$ have the following properties:

\begin{enumerate}[$(1)$]

\item There exist a ${\rm GL}_{2}\times{\rm GL}_{2}$-action on $\Upsilon_{\mA}^{14}$
and a ${\rm GL}_{2}$-action on $\Pi_{\mA}^{15}$ (Propositions \ref{prop:For-any-elementGL2GL2Up}
and \ref{prop:GL}).

\item We set $\Sigma:=\Upsilon_{\mA}^{14}$ or $\Pi_{\mA}^{15}$.
Let $\mA_{{\rm B}}$ be the affine space with the coordinates $z,s_{i},t_{i}\,(1\leq i\leq4)$
for $\Sigma=\Upsilon_{\mA}^{14}$, and the affine space with the coordinates
(\ref{eq:tenparam}) for $\Sigma=\Pi_{\mA}^{15}$. Let $\widehat{\Sigma}$
be the variety in $\mP(J)\times\mA_{\mathsf{B}}$ with the same equations
as $\Sigma$. A general fiber of the natural projection $\rho\colon\widehat{\Sigma}\to\mA_{\mathsf{B}}$
is isomorphic to $\mP^{2}\times\mP^{2}$ (Proposition \ref{prop:P2P2Up}
and \ref{prop:P2P2Pi}).

\item The varieties $\Upsilon_{\mA}^{14}$ and $\Pi_{\mA}^{15}$
are respectively $14$- and $15$-dimensional irreducible normal varieties
with only factorial Gorenstein terminal singularities (\cite{Tay}
for Gorensteinness and irreducibility of $\Upsilon_{\mA}^{14}$, Theorem
\ref{thm:PropUp}, Propositions \ref{prop:916}, \ref{prop:UFDPi}
and \ref{prop:MoreSing} for others).

\item Each of the ideals of $\Upsilon_{\mA}^{14}$ and $\Pi_{\mA}^{15}$
is generated by $9$ elements with $16$ relations (\cite{Tay} for
$\Upsilon_{\mA}^{14}$, Proposition \ref{prop:916} for $\Pi_{\mA}^{15}$)
. 

\end{enumerate}
\end{thm}

\subsection{Relation with previous works}

Due to the theory of unprojection, the affine variety $\Upsilon_{\mA}^{14}$
is constructed in \cite{Tay} and \cite{Tak2} independently, and
$\Pi_{\mA}^{15}$ is constructed in \cite{Tak3} inspired by \cite{Tay}.
We give new constructions of them in this paper. In this paper and
the previous one \cite{Tak5}, we provide to $\mathscr{H}_{\mA}^{13}$,
$\Upsilon_{\mA}^{14}$ and $\Pi_{\mA}^{15}$ a unified perspective
of the quadratic Jordan algebra of a cubic form.

In \cite{Tay}, candidates of prime $\mQ$-Fano 3-folds are constructed
from $\Upsilon_{\mA}^{14}$ and $\Pi_{\mA}^{15}$ (they should be
proved to have Picard number $1$). In a future, we will construct
prime $\mQ$-Fano 3-folds from $\Upsilon_{\mA}^{15}$ using its description
in this paper. In \cite{Tak3}, we construct prime $\mQ$-Fano 3-folds
from $\Pi_{\mA}^{15}$. 

In the paper \cite{Tak4}, we construct key varieties for prime $\mQ$-Fano
3-folds related with $\mP^{1}\times\mP^{1}\times\mP^{1}$-fibration,
and study the relation of them with the $G_{2}^{(4)}$-cluster variety.
Our constructions of them are based on the theory of Freudenthal triple
system. 

\vspace{5pt}

\noindent\textbf{Notation:} Throughout this paper, we denote by $J$
the quadratic Jordan algebra of a cubic form over a field $\mathsf{k}$
with the properties (a)--(c) as in \cite[Sect.2]{Tak5}. We also
refer to ibid. the very basic properties of $J$, which are taken
from \cite{Mc,R}.

\subsection{Structure of the paper}

Sections \ref{sec:CoordJ9param} and \ref{sec:Coord10paramJ} are
devoted to studying the quadratic Jordan algebra of a cubic form,
and, in the remaining part, we are concerned with the classification
of prime $\mQ$-Fano $3$-folds. 

\vspace{3pt}

\noindent Section \ref{sec:CoordJ9param}: In Subsection \ref{subsec:Facts-about-a single},
we collect basic facts about a single primitive idempotent of $J$
and the associated Peirce decomposition $J=J_{\nicefrac{1}{2}}\oplus J_{0}\oplus J_{1}$.
In the remaining part of Section \ref{sec:CoordJ9param}, we assume
that $\mathsf{k}$ is algebraically closed. In Subsection \ref{subsec:A-good-2-dimensional},
we find a good $2$-dimensional subspace $V$ of $J_{\nicefrac{1}{2}}$
with respect to the $\sharp$-product (Theorem \ref{thm:TypeII1Jordan}).
This result is the key to obtaining the desired coordinatization of
$J$ with 9 parameters in Subsection \ref{subsec:Coordination-of-9param}
(Propositions \ref{prop:vpi}, \ref{prop:pichipsi}, and \ref{prop:J12sha}).
In Subsection \ref{subsec:-mapping-and-cubicformUp}, we write down
the $\sharp$-mapping and the cubic form of $J$ explicitly using
the coordinatization obtained in Subsection \ref{subsec:Coordination-of-9param}
(Corollary \ref{cor:Upsilonshar}). 

\vspace{3pt}

\noindent Section \ref{sec:Affine-varietyUp}: We assume that $\mathsf{k}=\mC$.
The description of the $\sharp$-mapping in Section \ref{sec:CoordJ9param}
leads us to define the affine variety $\Upsilon_{\mA}^{14}$ (Definition
\ref{def:Upsilon}), which surprisingly coincides with the variety
defined in \cite{Tak2} and \cite{Tay}. After reviewing the properties
of $\Upsilon_{\mA}^{14}$ obtained by \cite{Tay} (Theorem \ref{thm:PropUp}),
we add more properties summarized as in Theorem \ref{thm:UpPiKeyMain},
which will be crucial to construct prime $\mQ$-Fano $3$-folds from
$\Upsilon_{\mA}^{14}$ in a future. 

\vspace{3pt}

\noindent Section \ref{sec:Coord10paramJ}: In Subsection \ref{subsec:Coordination-of-Perterson},
we obtain a coordinatization of $J$ with $8$ parameters due to the
result of Petersson (\cite[Prop.6.6]{Pe1}). In the remaining part
of Section \ref{sec:Coord10paramJ}, we assume that the characteristic
of $\mathsf{k}$ is not $3$. Taking a base change of $J$ and introducing
new parameters, we arrive at the desired coordinatization of $J$
in Subsection \ref{subsec:Change-of-base}. We write down the $\sharp$-product
of $J$ (Proposition \ref{prop:The--mapping-Pi}).

\vspace{3pt}

\noindent Section \ref{sec:Affine-varietyPi}: We assume that $\mathsf{k}=\mC$.
The description of the $\sharp$-mapping in Section \ref{sec:Coord10paramJ}
leads us to define the affine variety $\Pi_{\mA}^{15}$ (Definition
\ref{def:Pi15}), which also surprisingly coincides with the variety
defined in \cite{Tak3}. In the remaining part of Section \ref{sec:Affine-varietyPi},
we show Theorem \ref{thm:UpPiKeyMain} for $\Pi_{\mA}^{15}$.

\vspace{3pt}

\noindent Section \ref{sec:Singularities-ofUpPi}: We study the singularities
of $\Upsilon_{\mA}^{14}$ and $\Pi_{\mA}^{15}$. 

\vspace{3pt}

This paper contains several assertions which can be proved by straightforward
computations; we often omit such computation processes. Some computations
are difficult by hand but are easy within a software package. In our
computations, we use intensively the software systems Mathematica
\cite{W} and \textsc{Singular} \cite{DGPS}. 

\vspace{3pt}

\noindent\textbf{Acknowledgment:} I am very grateful to Professor
Shigeru Mukai; inspired by his article \cite{Mu1}, I was led to Jordan
algebra to describe key varieties. I owe many important calculations
in the paper to Professor Shinobu Hosono. I wish to thank him for
his generous cooperations. I would like to dedicate this paper to
my father, who passed away while I was obtaining part of the results.
This work is supported in part by Grant-in Aid for Scientific Research
(C) 16K05090. 

\vspace{5pt}

Before getting into the main part of the paper, we mention that unusual
notation numberings introduced in Subsections \ref{subsec:Facts-about-a single}--\ref{subsec:Coordination-of-9param}
and Subsection \ref{subsec:Change-of-base} come from the notation
for $\Upsilon_{\mA}^{14}$ in \cite{Tak2} and the notation $\Pi_{\mA}^{15}$
in \cite{Tak3} respectively. Finally, the notation become not so
weird in Corollary \ref{cor:Upsilonshar} and Definition \ref{def:Upsilon},
and Proposition \ref{prop:The--mapping-Pi} and Definition \ref{def:Pi15},
respectively.

\section{\textbf{Coordinatization of $J$ with $9$ parameters \label{sec:CoordJ9param}}}

\subsection{Facts about a single primitive idempotent\label{subsec:Facts-about-a single}}

A nonzero element $e\in J$ is called an \textit{idempotent }if it
satisfies $e\bullet e=e$, where $\empty\bullet\empty$denotes the
Jordan product of $J$, and an idempotent is called\textit{ primitive
}if it also satisfies $e^{\sharp}=0$. We define 
\begin{align*}
U_{x}y & :=T(x,y)x-x^{\sharp}\sharp y\,(x,y\in J),\\
U_{x_{1},x_{2}}(y) & :=U_{\mathfrak{1}}(y)-U_{x_{1}}(y)-U_{x_{2}}(y)\,(x_{1},x_{2},y\in J).
\end{align*}

By \cite{R}, $J$ contains a primitive idempotent $\upsilon$. We
set 
\begin{equation}
\pi_{3}:=\mathfrak{1}-\upsilon,\label{eq:pi3-0}
\end{equation}
which is also an idempotent. By \cite[the line just below (23)]{R},
we have 
\begin{equation}
T(\pi_{3})=T(\pi_{3},\pi_{3})=2.\label{eq:Tpi3}
\end{equation}
The primitive idempotent $\upsilon$ gives the following decomposition
of $J$, which is called the\textit{ Peirce decomposition} (associated
to $\upsilon$): 
\[
J=J_{1/2}\oplus J_{0}\oplus J_{1},
\]
 where 
\[
J_{1/2}:=U_{\upsilon,\pi_{3}}(J),\,J_{0}:=U_{\pi_{3}}(J),\,J_{1}:=U_{\upsilon}(J)
\]
(\cite[(24)]{R}). The subspaces $J_{l}$ are called the \textit{Peirce
subspaces}. For an element $\pi\in J$, we set 
\[
\pi^{\perp}:=\{T(*,\pi)=0\}.
\]
 The following facts are known from \cite{R}:
\begin{thm}
\label{thm:Peirceracine2}The following hold:
\end{thm}

\begin{enumerate}[$($a$)$]

\item $J_{1}=\mathsf{k}\upsilon$, $\dim J_{\nicefrac{1}{2}}=\dim J_{0}$
(\cite[(27), the arg. under (43)]{R}).

\item $\pi_{3}\in J_{0}$ (\cite[(27)]{R}). 

\item The quadratic form $S|_{J_{0}}$ is nondegenerate (\cite[p.97, the 5th line of the last par.]{R}).

\item $T(*,\pi_{3})=S(*,\pi_{3})$ on $J_{0}$ (\cite[the proof of Lem.2]{R}). 

\end{enumerate}
\begin{lem}
\label{Lem:J1212}For $\pi\in J_{0}$ and $\rho\in J_{\nicefrac{1}{2}}$,
we have $\pi\sharp\rho\in J_{\nicefrac{1}{2}}$. Therefore we obtain
the linear transformation 
\[
f_{\pi}\colon J_{\nicefrac{1}{2}}\to J_{\nicefrac{1}{2}}
\]
 defined by taking $\sharp$ product with $\pi\in J_{0}$.
\end{lem}

\begin{proof}
Indeed, it suffices to check $T(\pi\sharp\rho,\upsilon)=0$ and $T(\pi\sharp\rho,\phi)=0$
for $\phi\in J_{0}$ since the Peirce decomposition is orthogonal
with respect to the nondegenerate bilinear trace $T$ by \cite[(25)]{R}.
The first one follows since $T(\pi\sharp\rho,\upsilon)=T(\pi,\rho\sharp\upsilon)=0$,
where the first equality follows from \cite[(12)]{R}, and the second
equality follows from \cite[(28)]{R}. The second one follows since
$T(\pi\sharp\rho,\phi)=T(\rho,\pi\sharp\phi)=0$, where the second
equality follows from \cite[(25)]{R} and the bilinearization of \cite[(31)]{R}.
\end{proof}
\begin{lem}
\label{lem:SJ0-is-nondegenerate.}The quadratic form $S|_{J_{0}\cap\pi_{3}^{\perp}}$
is nondegenerate. 
\end{lem}

\begin{proof}
Indeed, since $T(\pi_{3},\pi_{3})=T(\pi_{3},\pi_{3}+\upsilon)=2\not=0$
by \cite[(23)]{R}, we have $\text{\ensuremath{\dim} }J_{0}\cap\pi_{3}^{\perp}=3$
and $J_{0}=\mathsf{k}\pi_{3}\oplus(\pi_{3}^{\perp}|_{J_{0}})$. Let
$\pi$ be an element of $J_{0}$ and we write $\pi=\pi_{0}+\alpha\pi_{3}$
with some $\pi_{0}\in\pi_{3}^{\perp}$ and $\alpha\in\mathsf{k}$.
We have $S(\pi,\pi)=S(\pi_{0},\pi_{0})+2\alpha S(\pi_{0},\pi_{3})+\alpha^{2}S(\pi_{3},\pi_{3}).$
By Theorem \ref{thm:Peirceracine2} (d), we see that $S(\pi_{0},\pi_{3})=T(\pi_{0},\pi_{3})=0$
and $S(\pi_{3},\pi_{3})=T(\pi_{3},\pi_{3})=2.$ Thus we have $S(\pi,\pi)=S(\pi_{0},\pi_{0})+2\alpha^{2}.$
From this and Theorem \ref{thm:Peirceracine2} (c), we see that $S|_{J_{0}\cap\pi_{3}^{\perp}}$
is nondegenerate.
\end{proof}
\vspace{5pt}

Hereafter in Section \ref{sec:CoordJ9param}, \textit{we assume that
$\dim J=9$ and $\mathsf{k}$ is an algebraically closed field.} 

\subsection{A good $2$-dimensional subspace of $J_{\nicefrac{1}{2}}$\label{subsec:A-good-2-dimensional}}

By Theorem \ref{thm:Peirceracine2} (a), we have

\[
\dim J_{\nicefrac{1}{2}}=\dim J_{0}=4.
\]

In the following theorem, we show that there exists a good $2$-dimensional
subspace in $J_{\nicefrac{1}{2}}$ with respect to the $\sharp$-product
with $J_{0}\cap\pi_{3}^{\perp}$; this is the key result to obtain
a coordinatization of $J$ with $9$ parameters in Subsection \ref{subsec:Coordination-of-9param}. 
\begin{thm}
\label{thm:TypeII1Jordan}

\begin{enumerate}[$(1)$]

\item The subalgebra ${\rm Cl}$ of ${\rm End}(J_{1/2})$ generated
by the identity map and $f_{\pi}\,(\pi\in J_{0}\cap\pi_{3}^{\perp})$
is the Clifford algebra associated to the $3$-dimensional space $J_{0}\cap\pi_{3}^{\perp}$
and the quadric form $-S|_{J_{0}\cap\pi_{3}^{\perp}}$.

\item Let $\pi_{1},\pi_{2},\pi_{4}$ be any basis of $J_{0}\cap\pi_{3}^{\perp}$
. There exists a $2$-dimensional subspace $V\subset J_{1/2}$ such
that 
\[
\pi_{i}\sharp V\subset V\,(i=2,4),\dim(\pi_{1}\sharp V)=2\ \text{and}\ J_{1/2}=V\oplus(\pi_{1}\sharp V).
\]

\end{enumerate}
\end{thm}

\begin{proof}
For $\pi\in J_{0}\cap\pi_{3}^{\perp}$, it holds that 
\begin{equation}
T(\pi)=T(\pi,\pi_{3}+\upsilon)=T(\pi,\pi_{3})=0.\label{eq:TpiJ0}
\end{equation}
For $\chi\in J_{\nicefrac{1}{2}}$, we have $T(\chi)=0$ by \cite[(26)]{R}.
We also have $T(\pi\sharp\chi)=T(\pi)T(\chi)-T(\pi,\chi)=0$ by \cite[(13), (25)]{R}.
Moreover, $\pi^{\sharp}\sharp\chi=S(\pi)\upsilon\sharp\chi=0$ by
\cite[(28), (31)]{R}. Therefore, by \cite[(21)]{R}, we have 
\begin{equation}
\pi\sharp(\pi\sharp\chi)=-T(\pi^{\sharp})\chi=-S(\pi)\chi,\label{eq:pipi}
\end{equation}
where the last equality follows from \cite[(16)]{R}. Thus we have
the assertion (1). 

\vspace{5pt}

The proof of (2) consists of several steps. 

\vspace{3pt}

\noindent \textbf{Step 1.} \textit{We find two candidates of good
coordinatizations of $J_{\nicefrac{1}{2}}$ as in (\ref{eq:M1}) and
(\ref{eq:M2twocase}) based upon representation theory of matrix algebra. }

Since $\dim(J_{0}\cap\pi_{3}^{\perp})=3$, it is well-known that the
Clifford algebra ${\rm Cl}$ is the direct sum of two copies of the
algebra $E$ of linear endomorphisms of a 2-dimensional vector space.
Therefore, by \cite[Thm.3.3.1]{Rep}, we may write $J_{\nicefrac{1}{2}}=W_{1}\oplus W_{2}$,
where $W_{i}\,(i=1,2)$ is an irreducible $2$-dimensional $E$-module.
Let $M_{1}\oplus M_{2}=\left(\begin{array}{cc}
M_{1} & O\\
O & M_{2}
\end{array}\right)$ be the representation matrix of the map $f_{\pi}\colon V_{\nicefrac{1}{2}}\to V_{\nicefrac{1}{2}}$
with respect to a basis $v_{1}$, $w_{1}$ of $W_{1}$ and a basis
$v_{2},w_{2}$ of $W_{2}$. Let $\rho_{1},\,\rho_{2},\,\rho_{3}$
be a basis of $J_{0}\cap\pi_{3}^{\perp}$, and $q_{1},q_{2},q_{3}$
the associated coordinates of $J_{0}\cap\pi_{3}^{\perp}$. By (\ref{eq:pipi}),
$M_{1}^{2}=M_{2}^{2}=-(S|_{J_{0}\cap p_{3}^{\perp}})E$. Since $S|_{J_{0}\cap p_{3}^{\perp}}$
is of rank $3$ by Lemma \ref{lem:SJ0-is-nondegenerate.}, neither
$M_{1}$ nor $M_{2}$ are multiples of the identity matrix by (\ref{eq:pipi}).
Therefore, by (\ref{eq:pipi}), we may write $M_{1}=\left(\begin{array}{cc}
f_{1} & f_{2}\\
f_{3} & -f_{1}
\end{array}\right)$ and $M_{2}=\left(\begin{array}{cc}
g_{1} & g_{2}\\
g_{3} & -g_{1}
\end{array}\right)$, where $f_{i},\,g_{i}\,(i=1,2,3)$ are linear forms with respect
to $q_{1},q_{2},q_{3}$ and it holds that $\det M_{1}=\det M_{2}=S|_{J_{0}\cap p_{3}^{\perp}}$.
Since $S|_{J_{0}\cap p_{3}^{\perp}}$ is of rank $3$ by Lemma \ref{lem:SJ0-is-nondegenerate.},
$f_{1},f_{2},f_{3}$ are linearly independent. Therefore, by a base
change of $J_{0}\cap\pi_{3}^{\perp}$, we may assume that 
\begin{equation}
M_{1}=\left(\begin{array}{cc}
q_{1} & q_{2}\\
q_{3} & -q_{1}
\end{array}\right).\label{eq:M1}
\end{equation}
Let $\iota\colon J_{0}\cap\pi_{3}^{\perp}\to J_{0}\cap\pi_{3}^{\perp}$
be the linear transformation defined by $q_{i}\mapsto g_{i}\,(i=1,2,3)$.
Since $\det M_{2}=\det M_{1}=-q_{1}^{2}-q_{2}q_{3}$, we have $g_{1}^{2}+g_{2}g_{3}=q_{1}^{2}+q_{2}q_{3}$.
This implies that $g_{1}$, $g_{2}$, $g_{3}$ are linearly independent.
Therefore $\iota$ is an isomorphism. Let $C:=\{q{}_{1}^{2}+q_{2}q_{3}=0\}$
in $\mP(J_{0}\cap\pi_{3}^{\perp})$. Since the conic $C$ is stable
by the projectivization of $\iota$, there exist $h\in{\rm GL}_{2}$
and $\alpha\in\mathsf{k^{\times}}$ such that the transformation $\left(\begin{array}{cc}
q_{1} & q_{2}\\
q_{3} & -q_{1}
\end{array}\right)\mapsto\left(\begin{array}{cc}
g_{1} & g_{2}\\
g_{3} & -g_{1}
\end{array}\right)$ coincides with $\left(\begin{array}{cc}
q_{1} & q_{2}\\
q_{3} & -q_{1}
\end{array}\right)\mapsto\alpha h^{-1}\left(\begin{array}{cc}
q_{1} & q_{2}\\
q_{3} & -q_{1}
\end{array}\right)h$ (this can be seen by using the Veronese embedding $\mP^{1}\to C$).
Changing the coordinates of $W_{2}$ by this $h$, we may assume that
$M_{2}=\alpha\left(\begin{array}{cc}
q_{1} & q_{2}\\
q_{3} & -q_{1}
\end{array}\right)$. Since $\det M_{2}=-q_{1}^{2}-q_{2}q_{3}$, we have $\alpha=\pm1.$
Therefore, we have 
\begin{equation}
M_{2}=\left(\begin{array}{cc}
q_{1} & q_{2}\\
q_{3} & -q_{1}
\end{array}\right)or\left(\begin{array}{cc}
-q_{1} & -q_{2}\\
-q_{3} & q_{1}
\end{array}\right).\label{eq:M2twocase}
\end{equation}

\vspace{1cm}

To step further, we need some auxiliary construction. By the form
of $M_{1}$ and $M_{2}$, we have $S(\rho_{1})=-1$. Then we may verify
that 
\begin{equation}
e_{1}:=\upsilon,e_{2}:=\nicefrac{1}{2}(\rho_{1}+\pi_{3}),e_{3}:=\nicefrac{1}{2}(-\rho_{1}+\pi_{3})\label{eq:e1e2e3}
\end{equation}
 are primitive idempotents satisfying the assumptions of $e_{1},e_{2},e_{3}$
in \cite{R}. We take the Peirce decomposition of $J$ as in ibid.
with respect to $e_{1},\,e_{2},\,e_{3}$, whose Peirce subspaces satisfy
\begin{align*}
 & J_{1}=J_{11},J_{0}=J_{22}\oplus J_{33}\oplus J_{23},J_{\nicefrac{1}{2}}=J_{12}\oplus J_{13},\\
 & J_{ii}=\mathsf{k}e_{i}\,(i=1,2,3).
\end{align*}
Note that, for $x\in J_{\nicefrac{1}{2}}$, we have 
\begin{equation}
\pi_{3}\sharp x=1\sharp x-\upsilon\sharp x=-x\label{eq:pi3x}
\end{equation}
 by the condition ($\sharp3$) of the $\sharp$-product as in \cite[Sect.2]{Tak5}
since $\upsilon\sharp x=0$ and $T(x)=0$ by \cite[(28)]{R}. Note
also that, for $x\in J_{12}$ (resp. $x\in J_{13}$), $e_{2}\sharp x=0$
(resp. $e_{3}\sharp x=0)$ by \cite[(33)]{R}. Therefore, by (\ref{eq:e1e2e3}),
$J_{12}$ and $J_{13}$ are contained in the $1$-eigenspace and the
$-1$-eigenspace of the map $f_{\rho_{1}}$, respectively. 

\vspace{10pt}

\noindent \textbf{Step 2.} \textit{We show that the latter case of
(\ref{eq:M2twocase}) occurs.} 

Assume, for a contradiction, that the former occurs. Then, by the
shape of $M_{1}\oplus M_{2}$ for $f_{\rho_{1}}$, the elements $v_{1},v_{2}$
form a basis of the $1$-eigenspace of the map $f_{\rho_{1}}$, and
the elements $w_{1},w_{2}$ form a basis of the $-1$-eigenspace of
the map $f_{\rho_{1}}$. Hence $J_{12}$ and $J_{13}$ coincides with
the $1$-eigenspace and the $-1$-eigenspace of the map $f_{\rho_{1}}$,
respectively since $\dim J_{12}=\dim J_{13}=2$ (note that, by the
proof of \cite[(43)]{R}, $J_{12}$, $J_{13}$, $J_{23}$ are mutually
isomorphic). We show that $\rho_{2},\rho_{3}$ form a basis of $J_{23}$.
Let $x$ be an element of $J_{23}$. Since $J_{23}\subset J_{0},$
we may write $x=\alpha\rho_{1}+\beta\rho_{2}+\gamma\rho_{3}+\delta\pi_{3}$
with $\alpha,\beta,\gamma,\delta\in\mathsf{k.}$ Note that, by the
shape of the matrix $M_{1}\oplus M_{2},$ we have 

\begin{equation}
\begin{cases}
\rho_{1}\sharp v_{1}=v_{1},\,\rho_{1}\sharp v_{2}=v_{2},\,\rho_{1}\sharp w_{1}=-w_{1},\,\rho_{1}\sharp w_{2}=-w_{2},\\
\rho_{2}\sharp v_{1}=0,\,\rho_{2}\sharp v_{2}=0,\,\rho_{2}\sharp w_{1}=v_{1},\,\rho_{2}\sharp w_{2}=v_{2},\\
\rho_{3}\sharp v_{1}=w_{1},\,\rho_{3}\sharp v_{2}=w_{2},\,\rho_{3}\sharp w_{1}=0,\,\rho_{3}\sharp w_{2}=0.
\end{cases}\label{eq:biprod}
\end{equation}
By (\ref{eq:pi3x}) and (\ref{eq:biprod}), we have $x\sharp v_{1}=(\alpha-\delta)v_{1}+\gamma w_{1}.$
Since $J_{23}\sharp J_{12}\subset J_{13}$ by \cite[(34)]{R}, we
have $\alpha=\delta$. Computing $x\sharp w_{1}$ similarly, we obtain
$\alpha=-\delta.$ Therefore, we have $\alpha=\delta=0$, and this
implies that $J_{23}\subset\langle\rho_{2},\rho_{3}\rangle$. Since
$\dim J_{23}=2$, we have $J_{23}=\langle\rho_{2},\rho_{3}\rangle$
as desired. Now we compute $v_{i}\sharp w_{j}\,(i=1,2,j=1,2)$. Since
$J_{12}\sharp J_{13}\subset J_{23}$ by \cite[(34)]{R}, we may write
$v_{i}\sharp w_{j}=\alpha_{ij}\rho_{2}+\beta_{ij}\rho_{3}$ with some
$\alpha_{ij},\beta_{ij}\in\mathsf{k}$. Then, by (\ref{eq:biprod}),
we have 
\begin{align*}
 & v_{i}\sharp(v_{i}\sharp w_{j})=v_{i}\sharp(\alpha_{ij}\rho_{2}+\beta_{ij}\rho_{3})=\beta_{ij}w_{i}.
\end{align*}
On the other hand, we have $v_{i}\sharp(v_{i}\sharp w_{j})=-S(v_{i})w_{j}$
by \cite[(36)]{R}. Therefore we have $\beta_{ij}w_{i}=-S(v_{i})w_{j.}$
From this for $i\not=j$, we obtain $\beta_{ij}=0$ and $S(v_{i})=0$,
and then we obtain $\beta_{ij}=0$ also for $i=j$ since $S(v_{i})=0$.
In a similar way, we also obtain $\alpha_{ij}=0$ for any $i=1,2$
and $j=1,2$. Therefore we have $v_{i}\sharp w_{j}=0$ for any $i=1,2$
and $j=1,2$. This means that $J_{12}\sharp J_{13}=0$, and hence
we have $S|_{J_{12}}=0$ or $S|_{J_{13}}=0$ by \cite[(35)]{R}. Thus
we have $J_{12}^{\sharp}=\{0\}$ or $J_{13}^{\sharp}=\{0\}$ by \cite[(33), 5th line]{R}.
This contradicts \cite[p.98, the 4th line from the bottom]{R}. Therefore
we have $M_{2}=\left(\begin{array}{cc}
-q_{1} & -q_{2}\\
-q_{3} & q_{1}
\end{array}\right)$ in (\ref{eq:M2twocase}). 

\vspace{5pt}

\noindent \textbf{Step 3.} \textit{Using geometry of the Grassmannian
${\rm G}(2,J_{\nicefrac{1}{2}})$, we express the condition as in
(\ref{eq:Gammap}) that a $2$-dimensional subspace $V\subset J_{\nicefrac{1}{2}}$
to be $f_{\pi}$-stable for an element $\pi\in J_{0}\cap\pi_{3}^{\perp}$.} 

Let $\pi=q_{1}\rho_{1}+q_{2}\rho_{2}+q_{3}\rho_{3}$ be an element
of $J_{0}\cap\pi_{3}^{\perp}$ such that $S(\pi)\not=0$ (such an
element exists by Lemma \ref{lem:SJ0-is-nondegenerate.}). By (\ref{eq:pipi}),
the map $f_{\pi}$ has two distinct eigenvalues $\pm\sqrt{-S(\pi)}.$
Let $E_{+}(\pi)$ and $E_{-}(\pi)$ be the $\sqrt{-S(\pi)}$- and
the $-\sqrt{-S(\pi)}$-eigenspaces for $f_{\pi}$, respectively. 

Assume that $q_{3}\not=0$. We express the eigenvectors of $f_{\pi}$
by the coordinates associated to the basis $v_{1},w_{1},v_{2},w_{2}$
of $J_{\nicefrac{1}{2}}$. We set $R:=\sqrt{q_{1}^{2}+q_{2}q_{3}}$
for simplicity of notation. By explicit computations using the shape
of $M_{1}\oplus M_{2}$, we see that 
\[
\phi_{1}:=\,\empty^{t}\!\left(\begin{array}{cccc}
q_{1}-R & q_{3} & 0 & 0\end{array}\right),\phi_{2}:=\,\empty^{t}\!\left(\begin{array}{cccc}
0 & 0 & q_{1}+R & q_{3}\end{array}\right)
\]
form a basis of $E_{+}(\pi)$, and 
\[
\psi_{1}:=\,\empty^{t}\!\left(\begin{array}{cccc}
q_{1}+R & q_{3} & 0 & 0\end{array}\right),\psi_{2}:=\,\empty^{t}\!\left(\begin{array}{cccc}
0 & 0 & q_{1}-R & q_{3}\end{array}\right)
\]
form a basis of $E_{-}(\pi)$. Let $L$ be a $2$-dimensional subspace
of $J_{\nicefrac{1}{2}}$. When $L\not=E_{+}(\pi),E_{-}(\pi)$, $L$
is $f_{\pi}$-stable if and only if $L$ intersects both $E_{+}(\pi)$
and $E_{-}(\pi)$ nontrivially. Such an $L$ is spanned by some linearly
independent two vectors $\alpha\phi_{1}+\beta\phi_{2}$ and $\gamma\psi_{1}+\delta\psi_{2}\,(\alpha,\beta,\gamma,\delta\in\mathsf{k})$.
Let $p_{ij}\,(1\leq i<j\leq4)$ be the Pl\"ucker coordinates for
${\rm G}(2,J_{\nicefrac{1}{2}})$ with respect to the basis $v_{1},w_{1},v_{2},w_{2}$
of $J_{\nicefrac{1}{2}}$. Since the basis $\alpha\phi_{1}+\beta\phi_{2}$,
$\gamma\psi_{1}+\delta\psi_{2}$ of $L$ forms the $2\times4$ matrix
\[
\left(\begin{array}{cccc}
\alpha(q_{1}-R) & \alpha q_{3} & \beta(q_{1}+R) & \beta q_{3}\\
\gamma(q_{1}+R) & \gamma q_{3} & \delta(q_{1}-R) & \delta q_{3}
\end{array}\right),
\]
the point $[L]$ of ${\rm G}(2,J_{\nicefrac{1}{2}})$ corresponding
to $L$ is contained in the locus 
\begin{equation}
\Gamma_{\pi}:={\rm G}(2,J_{\nicefrac{1}{2}})\cap\{p_{14}=p_{23},\,2q_{1}p_{23}+q_{2}p_{24}-q_{3}p_{13}=0\}.\label{eq:Gammap}
\end{equation}
From this calculation, we also see that, if $q_{3}\not=0$, then 
\begin{equation}
2\text{-dimensional}\ f_{p}\text{-stable subspaces \ensuremath{\not=E_{+}(\pi),E_{-}(\pi)} of \ensuremath{J_{\nicefrac{1}{2}}} are parameterized by}\ \Gamma_{\pi}.\label{eq:fpstable}
\end{equation}
 By similar calculations, we see that this fact (\ref{eq:fpstable})
holds even in the case that $q_{3}=0$. Note that the 2-dimensional
subspaces $W_{1}=\langle v_{1},w_{1}\rangle$ and $W_{2}=\langle v_{2},w_{2}\rangle$
correspond to two points of $\Gamma_{\pi}$.

\vspace{10pt}

To show the assertion (2), we may assume that $S(\pi_{2})\not=0$
and $S(\pi_{4})\not=0$. Indeed, there exist $\alpha,\beta,\gamma,\delta\in\mathsf{k}$
such that $\alpha\pi_{2}+\beta\pi_{4}$ and $\gamma\pi_{2}+\delta\pi_{4}$
are linearly independent (equivalently, $\alpha\delta-\beta\gamma\not=0)$
and $S(\alpha\pi_{2}+\beta\pi_{4})\not=0$ and $S(\gamma\pi_{2}+\delta\pi_{4})\not=0$
since $S|_{J_{0}\cap\pi_{3}^{\perp}}$ is of rank 3 by Lemma \ref{lem:SJ0-is-nondegenerate.}.
If a 2-dimensional subspace $L$ of $J_{\nicefrac{1}{2}}$ is stable
for $f_{\alpha\pi_{2}+\beta\pi_{4}}$ and $f_{\gamma\pi_{2}+\delta\pi_{4}}$,
then $L$ is also stable for $f_{\pi_{2}}$ and $f_{\pi_{4}}.$ Thus,
if such an $L$ satisfies the condition of (2) for $\alpha\pi_{2}+\beta\pi_{4},\,\gamma\pi_{2}+\delta\pi_{4},\,\pi_{1}$
instead of $\pi_{2},\,\pi_{4},\,\pi_{1}$, then $L$ satisfies the
condition of (2) for $\pi_{2},\,\pi_{4},\,\pi_{1}$. 

\vspace{10pt}

\noindent \textbf{Step 4. }\textit{Finally, we finish the proof of
(2).} 

By (\ref{eq:fpstable}), $\Gamma_{\pi_{2}}\cap\Gamma_{\pi_{4}}$ parameterizes
$2$-dimensional $f_{\pi_{2}}$- and $f_{\pi_{4}}$-stable subspaces.
Note that $\Gamma_{\pi_{2}}\cap\Gamma_{\pi_{4}}$ is a conic by (\ref{eq:Gammap})
since $\pi_{2}$ and $\pi_{4}$ are linearly independent. Therefore,
there exists a 2-dimensional subspace $V\not=W_{1},\,W_{2}$ of $J_{\nicefrac{1}{2}}$
which corresponds to a point of $\Gamma_{\pi_{2}}\cap\Gamma_{\pi_{4}}$.
It remains to show that $\dim f_{\pi_{1}}(V)=2$ and $V\cap f_{\pi_{1}}(V)=\{0\}.$ 

\vspace{3pt}

\noindent \textbf{Case $S(\pi_{1})\not=0$}: In this case, it holds
that $\dim f_{\pi_{1}}(V)=2$ since $f_{\pi_{1}}$ is an isomorphism
by (\ref{eq:pipi}). We show that $V$ is not $f_{\pi_{1}}$-stable.
Assume for a contradiction that $V$ is also $f_{\pi_{1}}$-stable.
Then $V$ is $f_{\pi}$-stable for any $\pi\in J_{\nicefrac{1}{2}}\cap\pi_{3}^{\perp}$
since $\pi_{1},\pi_{2},\pi_{4}$ span $J_{\nicefrac{1}{2}}\cap\pi_{3}^{\perp}$.
Let $\pi_{a}:=\rho_{2}-\rho_{3}$, and $\pi_{b}:=\rho_{2}+\rho_{3}$.
Note that $S(\pi_{a})\not=0$ and $S(\pi_{b})\not\not=0$. Then the
point $[V]\in{\rm G}(2,J_{\nicefrac{1}{2}})$ corresponding to $V$
is contained in $\Gamma_{\rho_{1}}\cap\Gamma_{\pi_{a}}\cap\Gamma_{\pi_{b}}$.
By an explicit calculation, we see that $\Gamma_{\rho_{1}}\cap\Gamma_{\pi_{a}}\cap\Gamma_{\pi_{b}}$
consists of the two points corresponding to $W_{1}$ and $W_{2}$.
This contradicts the choice of $V$. Therefore $V$ is not $f_{\pi_{1}}$-stable.
If $\ell:=V\cap f_{\pi_{1}}(V)$ is 1-dimensional, then $\ell$ is
contained in an $f_{\pi_{1}}$-eigenspace, say E, since $f_{\pi_{1}}^{2}$
is a scalar map by (\ref{eq:pipi}). Therefore $V$ intersects $E$
non-trivially. For a fixed 2-dimensional subspace $F$ of $J_{\nicefrac{1}{2}}$,
it is well-known that 2-dimensional subspaces of $J_{\nicefrac{1}{2}}$
intersecting $F$ nontrivially are parameterized by a hyperplane section
$H_{F}$ of ${\rm G}(2,J_{\nicefrac{1}{2}})$, and $H_{F}$ is singular
at the point $[F]$ of ${\rm G}(2,J_{\nicefrac{1}{2}})$. Since $\Gamma_{\pi_{1}}$
parameterizes 2-dimensional subspaces of $J_{\nicefrac{1}{2}}$ intersecting
the two eigenspaces $E$, $E'$ of $f_{\pi_{1}}$ nontrivially, it
holds that $\Gamma_{\pi_{1}}=H_{E}\cap H_{E'}$. Note that $H_{E}\not={\rm G}(2,J_{\nicefrac{1}{2}})\cap\{p_{14}=p_{23}\}$
since $H_{E}$ is singular while ${\rm G}(2,J_{\nicefrac{1}{2}})\cap\{p_{14}=p_{23}\}$
is nonsingular. Since $\Gamma_{\pi_{1}}$ is also contained in $\{p_{14}=p_{23}\}$,
it holds that $\Gamma_{\pi_{1}}=H_{E}\cap\{p_{14}=p_{23}\}$. Therefore
$[V]$ is contained in $\Gamma_{\pi_{1}}$, which means that $V$
is $f_{\pi_{1}}$-stable, a contradiction. 

\noindent \textbf{Case $S(\pi_{1})=0$}: We can choose some $\epsilon\in\mathsf{k}$
such that $S(\pi_{1}+\epsilon\pi_{2})\not=0$ since $S|_{J_{0}\cap\pi_{3}^{\perp}}$
is of rank 3 by Lemma \ref{lem:SJ0-is-nondegenerate.}. Therefore,
by the discussion above, we have $\dim f_{\pi_{1}+\epsilon\pi_{2}}(V)=2$
and $V\cap f_{\pi_{1}+\epsilon\pi_{2}}(V)=\{0\}.$ This is equivalent
to that, for a basis $\tau_{1},\tau_{2}$ of $V$, 4 vectors $\tau_{1},\,\tau_{2},\,f_{\pi_{1}+\epsilon\pi_{2}}(\tau_{1}),\,f_{\pi_{1}+\epsilon\pi_{2}}(\tau_{2})$
are linearly independent. Since $f_{\pi_{1}+\epsilon\pi_{2}}(\tau_{i})=f_{\pi_{1}}(\tau_{i})+\epsilon f_{\pi_{2}}(\tau_{i})$,
and $f_{\pi_{2}}(\tau_{i})\in V\,(i=1,2)$, we see that $\tau_{1},\,\tau_{2},\,f_{\pi_{1}}(\tau_{1}),\,f_{\pi_{1}}(\tau_{2})$
are also linearly independent. Therefore we have $\dim f_{\pi_{1}}(V)=2$
and $V\cap f_{\pi_{1}}(V)=\{0\}$ as desired. 
\end{proof}

\subsection{Coordinatization of $J$ with $9$ parameters\label{subsec:Coordination-of-9param}}

In this subsection, we coordinatize $J$ introducing 9 parameters
and determining the $\sharp$-products between elements of suitable
bases of the Peirce spaces $J_{\nicefrac{1}{2}},J_{0},J_{1}$. 
\begin{lem}
\label{Lem:J0J0}For an element $\sigma_{0}\in J_{0}$, it holds that
\begin{equation}
\upsilon\sharp\sigma_{0}=T(\sigma_{0},\pi_{3})\pi_{3}-\sigma_{0}.\label{eq:usigma}
\end{equation}
 In particular, $\upsilon\sharp\sigma_{0}\in J_{0}$.
\end{lem}

\begin{proof}
By the condition ($\sharp3$) of the $\sharp$-mapping as in \cite[Sect.2]{Tak5},
we have $\mathfrak{1}\sharp\sigma_{0}=T(\sigma_{0})\mathfrak{1}-\sigma_{0}$.
Since $T(\sigma_{0},\upsilon)=0$ by \cite[(25)]{R} and the fact
that $\mathfrak{1}=\upsilon+\pi_{3}$, we have $\upsilon\sharp\sigma_{0}=T(\sigma_{0},\pi_{3})\pi_{3}-\sigma_{0}+T(\sigma_{0},\pi_{3})\upsilon-\pi_{3}\sharp\sigma_{0}.$
Bilinearizing \cite[(16), (31)]{R}, we have $\pi_{3}\sharp\sigma_{0}=S(\pi_{3},\sigma_{0})\upsilon$.
By the fact that $T(\pi_{3})=2$ and $T(\sigma_{0},\upsilon)=0$,
we have $S(\pi_{3},\sigma_{0})=T(\pi_{3})T(\sigma_{0})-T(\pi_{3},\sigma_{0})=T(\pi_{3})T(\sigma_{0},\pi_{3})-T(\pi_{3},\sigma_{0})=T(\pi_{3},\sigma_{0}).$
Therefore we obtain the desired result.
\end{proof}
By \cite[(28), (30), (31)]{R}, and Lemmas \ref{Lem:J1212} and \ref{Lem:J0J0},
we obtain the following:
\begin{cor}
\label{cor:sharpformula}Let $\sigma\in J$ be any element and we
write 
\[
\sigma=u\upsilon+\sigma_{0}+\sigma_{\nicefrac{1}{2}}
\]
 with $u\in\mathsf{k}$, $\sigma_{0}\in J_{0}$ and $\sigma_{\nicefrac{1}{2}}\in J_{\nicefrac{1}{2}}$.
It holds that 
\begin{align*}
\sigma^{\sharp} & =S(\sigma_{0})\upsilon+\left((\sigma_{\nicefrac{1}{2}})^{\sharp}+u(\upsilon\sharp\sigma_{0})\right)+\sigma_{0}\sharp\sigma_{\nicefrac{1}{2}},
\end{align*}
where $(\sigma_{\nicefrac{1}{2}})^{\sharp}+u(\upsilon\sharp\sigma_{0})\in J_{0}$
and $\sigma_{0}\sharp\sigma_{\nicefrac{1}{2}}\in J_{\nicefrac{1}{2}}$. 
\end{cor}

Let $\upsilon,\pi_{1},\pi_{2},\pi_{3},\pi_{4}$ be as above. Let $V\subset J_{\nicefrac{1}{2}}$
be as in Theorem \ref{thm:TypeII1Jordan} (2), and let $\psi_{1},\psi_{2}$
form a basis of $V$. We set 
\begin{equation}
\chi_{1}:=-\pi_{1}\sharp\psi_{1}\ \text{and}\ \chi_{2}:=-\pi_{1}\sharp\psi_{2},\label{eq:chaidef}
\end{equation}
which are also elements of $J_{\nicefrac{1}{2}}$. By ibid., we see
that 
\begin{equation}
\upsilon,\pi_{1},\pi_{2},\pi_{3},\pi_{4},\psi_{1},\psi_{2},\chi_{1},\chi_{2}\label{eq:Jbasis}
\end{equation}
 form a basis of $J$, where $\pi_{1},\dots,\pi_{4}$ form a basis
of $J_{0}$, and $\psi_{1},\psi_{2},\chi_{1},\chi_{2}$ form a basis
of $J_{\nicefrac{1}{2}}$. 

We will express $\sigma^{\sharp}$ for any $\sigma\in J$ using the
coordinates with respect to this basis of $J$. For this, we determine
the $\sharp$-products among the elements of this basis in Propositions
\ref{prop:vpi}, \ref{prop:pichipsi}, and \ref{prop:J12sha}. 
\begin{prop}
\label{prop:vpi}It holds that 
\begin{equation}
\begin{cases}
\upsilon^{\sharp} & =0,\\
\upsilon\sharp\psi_{i}=\upsilon\sharp\chi_{j} & =0\,(i,j=1,2),\\
v\sharp\pi_{3} & =\pi_{3},\\
v\sharp\pi_{i} & =-\pi_{i}\,(i=1,2,4).
\end{cases}\label{eq:vpi}
\end{equation}
 
\end{prop}

\begin{proof}
Since $\upsilon$ is a primitive idempotent, we have $\upsilon^{\sharp}=0$.
Since $\psi_{i},\chi_{j}\in J_{\nicefrac{1}{2}}\,(i,j=1,2)$, we have
$\upsilon\sharp\psi_{i}=\upsilon\sharp\chi_{j}=0\,(i,j=1,2)$ by \cite[(28)]{R}.
By Lemma \ref{Lem:J0J0} and (\ref{eq:Tpi3}), we have $v\sharp\pi_{3}=\pi_{3}$.
By Lemma \ref{Lem:J0J0} and the assumption $\pi_{i}\in J_{0}\cap\pi_{3}^{\perp}\,(i=1,2,4)$,
we have $v\sharp\pi_{i}=-\pi_{i}$. 
\end{proof}
\begin{prop}
\label{prop:pichipsi}It holds that 
\begin{equation}
\pi_{3}\sharp\chi_{i}=-\chi_{i}\ \text{and}\ \pi_{3}\sharp\psi_{i}=-\psi_{i}\,(i=1,2).\label{eq:pi3}
\end{equation}
Let 
\begin{equation}
s_{1}:=S(\pi_{1},\pi_{2}),\,t_{1}:=-S(\pi_{1},\pi_{4}),\label{eq:s1}
\end{equation}
 and 
\begin{equation}
z:=-S(\pi_{1}).\label{eq:zpi1}
\end{equation}
For some $s_{2},s_{3},s_{4},t_{2},t_{3},t_{4}\in\mathsf{k}$, it holds
that

\begin{equation}
\begin{cases}
\pi_{1}\sharp\pi_{2} & =s_{1}\upsilon,\\
\pi_{1}\sharp\pi_{4} & =t_{1}\upsilon,\\
\pi_{1}^{\sharp} & =-z\upsilon.
\end{cases}\label{eq:pi1}
\end{equation}

\begin{equation}
\begin{cases}
\pi_{2}\sharp\psi_{1}=-s_{3}\psi_{1}+s_{2}\psi_{2},\\
\pi_{2}\sharp\psi_{2}=-s_{4}\psi_{1}+s_{3}\psi_{2}.
\end{cases}\label{eq:pi2psi}
\end{equation}

\begin{equation}
\begin{cases}
\pi_{4}\sharp\psi_{1}=t_{3}\psi_{1}-t_{2}\psi_{2},\\
\pi_{4}\sharp\psi_{2}=t_{4}\psi_{1}-t_{3}\psi_{2}.
\end{cases}\label{eq:pi4psi}
\end{equation}

\begin{equation}
\begin{cases}
\pi_{2}^{\sharp} & =(-s_{3}^{2}+s_{2}s_{4})\upsilon,\\
\pi_{4}^{\sharp} & =(-t_{3}^{2}+t_{2}t_{4})\upsilon,\\
\pi_{2}\sharp\pi_{4} & =(2s_{3}t_{3}-s_{2}t_{4}-t_{2}s_{4})\upsilon.
\end{cases}\label{eq:pisha}
\end{equation}
\end{prop}

\begin{equation}
\begin{cases}
\pi_{2}\sharp\chi_{1}=s_{3}\chi_{1}-s_{2}\chi_{2}+s_{1}\psi_{1},\\
\pi_{2}\sharp\chi_{2}=s_{4}\chi_{1}-s_{3}\chi_{2}+s_{1}\psi_{2},\\
\pi_{4}\sharp\chi_{1}=-t_{3}\chi_{1}+t_{2}\chi_{2}-t_{1}\psi_{1},\\
\pi_{4}\sharp\chi_{2}=-t_{4}\chi_{1}+t_{3}\chi_{2}-t_{1}\psi_{2}.
\end{cases}\label{eq:pi2pi4chi}
\end{equation}

\begin{equation}
\pi_{1}\sharp\chi_{1}=-z\psi_{1},\,\pi_{1}\sharp\chi_{2}=-z\psi_{2}.\label{eq:z}
\end{equation}

\begin{proof}
By (\ref{eq:s1}) and (\ref{eq:zpi1}), we have (\ref{eq:pi1}) by
\cite[(31)]{R} and its bilinearization. 

By Theorem \ref{thm:TypeII1Jordan} (2) and (\ref{eq:pipi}), we have
(\ref{eq:pi2psi}) with some $s_{2},s_{3},s_{4}\in\mathsf{k}$, and
(\ref{eq:pi4psi}) using some $t_{2},t_{3},t_{4}\in\mathsf{k}$. 

By (\ref{eq:pipi}), (\ref{eq:pi2psi}) and (\ref{eq:pi4psi}), we
have 
\begin{equation}
S(\pi_{2})=-s_{3}^{2}+s_{2}s_{4},\,S(\pi_{4})=-t_{3}^{2}+t_{2}t_{4}.\label{eq:Spi2pi4}
\end{equation}
 Bilinearizing (\ref{eq:pipi}), we have 
\begin{equation}
\pi\sharp(\rho\sharp\chi)+\rho\sharp(\pi\sharp\chi)=-S(\pi,\rho)\chi\label{eq:pirho}
\end{equation}
for $\pi,\rho\in\pi_{3}^{\perp}$ and $\chi\in J_{\nicefrac{1}{2}}$.
By (\ref{eq:pi2psi}), (\ref{eq:pi4psi}), and (\ref{eq:pirho}),
we have 
\begin{equation}
S(\pi_{2},\pi_{4})=2s_{3}t_{3}-s_{2}t_{4}-t_{2}s_{4}.\label{eq:Sp12pi4}
\end{equation}
Therefore, we have (\ref{eq:pisha}) by \cite[(31)]{R} and its bilinearization. 

By (\ref{eq:pirho}), we have 
\[
\pi_{j}\sharp\chi_{i}=-\pi_{j}\sharp(\pi_{1}\sharp\psi_{i})=\pi_{1}\sharp(\pi_{j}\sharp\psi_{i})+S(\pi_{1},\pi_{j})\psi_{i}\,(i=1,2,j=2,4).
\]
Therefore we obtain (\ref{eq:pi2pi4chi}) from (\ref{eq:s1}), (\ref{eq:pi2psi}),
and (\ref{eq:pi4psi}). 

By (\ref{eq:pipi}) and (\ref{eq:zpi1}), we have (\ref{eq:z}). 
\end{proof}
Hereafter, we change $\pi_{2},\pi_{3},\psi_{1},\psi_{2}$ if necessary
such that the following hold:

\begin{align}
S(\pi_{2}) & \not=0,\,S(\pi_{4})\not=0.\label{eq:Spi2pi4not0}
\end{align}

\begin{equation}
\text{All of\ }s_{2},s_{3},s_{4},t_{2},t_{3},t_{4}\ \text{are nonzero.}\label{eq:nonzerocoeff}
\end{equation}

Indeed, by Lemma \ref{lem:SJ0-is-nondegenerate.}, we may take $\pi_{2}$
and $\pi_{4}$ satisfying (\ref{eq:Spi2pi4not0}). Then general $\psi_{1}$,
$\psi_{2}$ satisfy (\ref{eq:nonzerocoeff}).

\vspace{5pt}

To determine the $\sharp$-product on $J_{\nicefrac{1}{2}}$ is more
involved.
\begin{lem}
\label{Lem:psisharpnotzero}It holds that $\psi_{1}^{\sharp}\not=0$
and $\psi_{2}^{\sharp}\not=0$.
\end{lem}

\begin{proof}
We only show $\psi_{1}^{\sharp}\not=0$ since we can prove $\psi_{2}^{\sharp}\not=0$
similarly. Assume for a contradiction that $\psi_{1}^{\sharp}=0$.
Then, by \cite[(21), (25), (28)]{R}, we have 
\begin{equation}
\psi_{1}\sharp(\psi_{1}\sharp\pi_{i})=0\,(i=1,2).\label{eq:phiphipi}
\end{equation}
 Thus, from (\ref{eq:pi2psi}), we have $s_{2}(\psi_{1}\sharp\psi_{2})=0$.
Since $s_{2}\not=0$, we obtain $\psi_{1}\sharp\psi_{2}=0.$ By (\ref{eq:phiphipi}),
we have 
\begin{equation}
\psi_{1}\sharp\chi_{1}=-\psi_{1}\sharp(\psi_{1}\sharp\pi_{1})=0.\label{eq:phichi}
\end{equation}
By \cite[(11), (25)]{R}, we have $(\psi_{1}\sharp\pi_{i})\sharp\chi_{1}=-(\psi_{1}\sharp\pi_{2})\sharp(\psi_{1}\sharp\pi_{1})=0\,(i=2,4)$.
Thus, from (\ref{eq:pi2psi}), and (\ref{eq:phichi}), we have $s_{2}(\chi_{1}\sharp\psi_{2})=0$.
Since $s_{2}\not=0$, we obtain $\chi_{1}\sharp\psi_{2}=0$. Linearizing
\cite[(21)]{R} completely, the formula \cite[(7) in the pf. of Prop.6.6]{Pe1}
is derived. Applying this by setting $u=\pi_{1},\,v=\psi_{1},\,w=\psi_{2}$,
we obtain $\chi_{1}\sharp\psi_{1}+\chi_{2}\sharp\psi_{1}=0$ (for
this calculation, we also use $\psi_{1}\sharp\psi_{2}=0$). Thus,
by (\ref{eq:phichi}), we have 
\[
\chi_{2}\sharp\psi_{1}=0.
\]
 Therefore, the assumption that $\psi_{1}^{\sharp}=0$ implies that
$\psi_{1}\sharp\sigma=0$ for any $\sigma\in J_{\nicefrac{1}{2}}$.
Then, by \cite[(13), (28)]{R}, we have $T(\psi_{1},\sigma)=T(\psi_{1})T(\sigma)-T(\psi_{1}\sharp\sigma)=0$
for any $\sigma\in J_{\nicefrac{1}{2}}$. Therefore, by \cite[(25)]{R}
and nondegeneracy of $J$, we have $\psi_{1}=0$, a contradiction.
\end{proof}
\begin{prop}
\label{prop:J12sha}By replacing $\psi_{1}$ and $\psi_{2}$ with
their nonzero multiples if necessary, it holds that 
\begin{equation}
\begin{cases}
\psi_{1}^{\sharp} & =t_{2}\pi_{2}+(s_{3}t_{2}-s_{2}t_{3})\pi_{3}+s_{2}\pi_{4},\\
\psi_{1}\sharp\psi_{2} & =2t_{3}\pi_{2}+(t_{2}s_{4}-s_{2}t_{4})\pi_{3}+2s_{3}\pi_{4},\\
\psi_{2}^{\sharp} & =t_{4}\pi_{2}+(s_{4}t_{3}-s_{3}t_{4})\pi_{3}+s_{4}\pi_{4},\\
\chi_{1}^{\sharp} & =(s_{2}t_{1}-s_{1}t_{2})\pi_{1}-z\left(t_{2}\pi_{2}-(s_{3}t_{2}-s_{2}t_{3})\pi_{3}+s_{2}\pi_{4}\right),\\
\chi_{1}\sharp\chi_{2} & =2(s_{3}t_{1}-s_{1}t_{3})\pi_{1}-z\left(2t_{3}\pi_{2}-(s_{4}t_{2}-s_{2}t_{4})\pi_{3}+2s_{3}\pi_{4}\right),\\
\chi_{2}^{\sharp} & =(s_{4}t_{1}-s_{1}t_{4})\pi_{1}-z\left(t_{4}\pi_{2}-(s_{4}t_{3}-s_{3}t_{4})\pi_{3}+s_{4}\pi_{4}\right),\\
\chi_{1}\sharp\psi_{1} & =2(s_{3}t_{2}-s_{2}t_{3})\pi_{1}+(s_{2}t_{1}-s_{1}t_{2})\pi_{3},\\
\chi_{2}\sharp\psi_{2} & =2(s_{4}t_{3}-s_{3}t_{4})\pi_{1}+(s_{4}t_{1}-s_{1}t_{4})\pi_{3},\\
\chi_{1}\sharp\psi_{2} & =(s_{4}t_{2}-s_{2}t_{4})\pi_{1}-t_{1}\pi_{2}+(s_{3}t_{1}-s_{1}t_{3})\pi_{3}-s_{1}\pi_{4},\\
\chi_{2}\sharp\psi_{1} & =(s_{4}t_{2}-s_{2}t_{4})\pi_{1}+t_{1}\pi_{2}+(s_{3}t_{1}-s_{1}t_{3})\pi_{3}+s_{1}\pi_{4}.
\end{cases}\label{eq:xy}
\end{equation}
\end{prop}

\begin{proof}
~

\noindent$\underline{\psi_{1}^{\sharp},\,\psi_{2}^{\sharp},\psi_{1}\sharp\psi_{2}.}$
By \cite[(14), (28)]{R}, we have 
\begin{equation}
\psi_{1}\sharp\psi_{1}^{\sharp}=-N(\psi_{1})\mathfrak{1}-T(\psi_{1}^{\sharp})\psi_{1}.\label{eq:phi1phi1ad}
\end{equation}
 Since $\psi_{1}^{\sharp}\in J_{0}$, we may write 
\begin{equation}
\psi_{1}^{\sharp}=a_{1}\pi_{1}+a_{2}\pi_{2}+a_{3}\pi_{3}+a_{4}\pi_{4}\label{eq:phi1ad}
\end{equation}
 with some $a_{1},a_{2},a_{3},a_{4}\in\mathsf{k}$. Since $T(\pi_{i})=0\,(i=1,2,4)$
by (\ref{eq:TpiJ0}) and $T(\pi_{3})=2$ by (\ref{eq:Tpi3}), we have
$T(\psi_{1}^{\sharp})=2a_{3}.$ Inserting (\ref{eq:phi1ad}) in (\ref{eq:phi1phi1ad})
and using (\ref{eq:chaidef}), (\ref{eq:pi3}), (\ref{eq:pi2psi})
and (\ref{eq:pi4psi}), we have 

\begin{align*}
 & -a_{1}\chi_{1}+a_{2}\left(-s_{3}\psi_{1}+s_{2}\psi_{2}\right)-a_{3}\psi_{1}+a_{4}\left(t_{3}\psi_{1}-t_{2}\psi_{2}\right)\\
 & =-N(\psi_{1})\mathfrak{1}-2a_{3}\psi_{1}.
\end{align*}
Therefore, since $\mathfrak{1}\in J_{1}\oplus J_{0},$ we obtain $N(\psi_{1})=0$,
$a_{1}=0,\,a_{2}s_{2}-a_{4}t_{2}=0,\,a_{3}=a_{2}s_{3}-a_{4}t_{3}.$
By Lemma \ref{Lem:psisharpnotzero} and the assumption (\ref{eq:nonzerocoeff}),
we have 
\begin{equation}
\psi_{1}^{\sharp}=\alpha(t_{2}\pi_{2}+(s_{3}t_{2}-s_{2}t_{3})\pi_{3}+s_{2}\pi_{4})\label{eq:a2a3a4}
\end{equation}
for some $\alpha\in\mathsf{k}^{\times}$. Similarly, computing $\psi_{2}\sharp\psi_{2}^{\sharp}$,
we obtain 

\begin{equation}
\psi_{2}^{\sharp}=\beta\left(t_{4}\pi_{2}+(s_{4}t_{3}-s_{3}t_{4})\pi_{3}+s_{4}\pi_{4}\right)\label{eq:a2a3a4V2}
\end{equation}
for some $\beta\in\mathsf{k}^{\times}$.

By the formula \cite[(21)]{R}, we obtain
\begin{equation}
\psi_{1}\sharp(\psi_{1}\sharp\pi_{2})=S(\psi_{1}^{\sharp},\pi_{2})\pi_{3}-T(\psi_{1}^{\sharp})\pi_{2}.\label{eq:psi1psi1pi2}
\end{equation}
Indeed, it holds that $T(\psi_{1})=0$ by \cite[(28)]{R}, and $T(\pi_{2})=0$
by (\ref{eq:TpiJ0}). It holds that $T(\psi_{1}\sharp\pi_{2})=0$
by \cite[(13), (25), (28)]{R}. Since $\psi_{1}^{\sharp},\pi_{2}\in J_{0}$,
we have $\psi_{1}^{\sharp}\sharp\pi_{2}=S(\psi_{1}^{\sharp},\pi_{2})\upsilon$
by the bilinearization of \cite[(31)]{R}. By \cite[(13), the bilin. of (16)]{R},
we have $T(\psi_{1}^{\sharp}\sharp\pi_{2})=-S(\psi_{1}^{\sharp},\pi_{2}).$
Therefore, $T(\psi_{1}^{\sharp}\sharp\pi_{2})\mathfrak{1}-\psi_{1}^{\sharp}\sharp\pi_{2}=S(\psi_{1}^{\sharp},\pi_{2})\pi_{3}$.
Now (\ref{eq:psi1psi1pi2}) follows. Therefore, since $s_{2}\not=0$,
we obtain 

\begin{equation}
\psi_{1}\sharp\psi_{2}=\alpha(2t_{3}\pi_{2}+(t_{2}s_{4}-s_{2}t_{4})\pi_{3}+2s_{3}\pi_{4}).\label{eq:psi1psi2}
\end{equation}

Similarly, computing $\psi_{1}\sharp(\psi_{1}\sharp\pi_{4})$, we
obtain 

\[
\psi_{1}\sharp\psi_{2}=\beta(2t_{3}\pi_{2}+(t_{2}s_{4}-s_{2}t_{4})\pi_{3}+2s_{3}\pi_{4}).
\]

Therefore we have $\alpha=\beta$ by the assumption (\ref{eq:nonzerocoeff}).
Now we will multiply $\psi_{1}$ and $\psi_{2}$ with some nonzero
constants to obtain the presentations of $\psi_{1}^{\sharp}$, $\psi_{2}^{\sharp}$,
$\psi_{1}\sharp\psi_{2}$ as in (\ref{eq:xy}). If we replace $\psi_{1}$
and $\psi_{2}$ with $\gamma\psi_{1}$and $\delta\psi_{2}$ respectively
($\gamma$, $\delta\not=0$), then, by (\ref{eq:pi2psi}) and (\ref{eq:pi4psi}),
$s_{3}$ and $t_{3}$ do not change and $s_{2}$, $t_{2}$, $s_{4}$,
$t_{4}$ change to $s_{2}':=\frac{\gamma}{\delta}s_{2}$, $t_{2}':=\frac{\gamma}{\delta}t_{2}$,
$s_{4}':=\frac{\delta}{\gamma}s_{4}$, $t_{4}':=\frac{\delta}{\gamma}t_{4}$
respectively. By (\ref{eq:a2a3a4}), (\ref{eq:a2a3a4V2}) and (\ref{eq:psi1psi2}),
we have 
\begin{align*}
(\gamma\psi_{1})^{\sharp} & =\alpha\gamma\delta\left(t_{2}'\pi_{2}+(s_{3}t_{2}'-s_{2}'t_{3})\pi_{3}+s_{2}'\pi_{4}\right)\\
(\delta\psi_{2})^{\sharp} & =\alpha\gamma\delta\left(t_{4}'\pi_{2}+(s_{4}'t_{3}-s_{3}t_{4}')\pi_{3}+s_{4}'\pi_{4}\right)\\
(\gamma\psi_{1})\sharp(\delta\psi_{2}) & =\alpha\gamma\delta\left(2t_{3}\pi_{2}+(t_{2}'s_{4}'-s_{2}'t_{4}')\pi_{3}+2s_{3}\pi_{4}\right).
\end{align*}

Therefore choosing $\gamma$, $\delta$ such that $\alpha\gamma\delta=1$,
we obtain the presentations of $\psi_{1}^{\sharp}$, $\psi_{2}^{\sharp}$,
$\psi_{1}\sharp\psi_{2}$ as in (\ref{eq:xy}).

\vspace{3pt}

\noindent$\underline{\chi_{1}\sharp\psi_{1}\ \text{and}\ \chi_{2}\sharp\psi_{2}.}$
By the formula \cite[(21)]{R}, we obtain 
\[
\chi_{1}\sharp\psi_{1}=-(\pi_{1}\sharp\psi_{1})\sharp\psi_{1}=T(\psi_{1}^{\sharp})\pi_{1}-S(\pi_{1},\psi_{1}^{\sharp})\pi_{3}
\]
in a similar way to get (\ref{eq:psi1psi1pi2}). Then we obtain the
presentation of $\chi_{1}\sharp\psi_{1}$ as in (\ref{eq:xy}) using
the presentation of $\psi_{1}^{\sharp}$. We obtain the presentation
of $\chi_{2}\sharp\psi_{2}$ similarly.

\vspace{3pt}

\noindent$\underline{\chi_{1}^{\sharp},\chi_{2}^{\sharp}\ \text{and}\ \chi_{1}\sharp\chi_{2}.}$
Using the formula \cite[(10)]{R} and \cite[(25)]{R}, we have
\begin{align*}
\chi_{1}^{\sharp} & =(\pi_{1}\sharp\psi_{1})^{\sharp}=-\pi_{1}^{\sharp}\sharp\psi_{1}^{\sharp}+T(\psi_{1}^{\sharp},\pi_{1})\pi_{1}.
\end{align*}
From this and the presentation of $\psi_{1}^{\sharp}$, we obtain
the presentation of $\chi_{1}^{\sharp}$ as in (\ref{eq:xy}). We
obtain the presentation of $\chi_{2}^{\sharp}$ similarly. 

Using the formula \cite[(11)]{R}, we have 

\[
\chi_{1}\sharp\chi_{2}=(\pi_{1}\sharp\psi_{1})\sharp(\pi_{1}\sharp\psi_{2})=-\pi_{1}^{\sharp}\sharp(\psi_{1}\sharp\psi_{2})+T(\psi_{1}\sharp\psi_{2},\pi_{1})\pi_{1},
\]
and obtain the presentation of $\chi_{1}\sharp\chi_{2}$ as in (\ref{eq:xy}). 

\vspace{3pt}

\noindent$\underline{\chi_{1}\sharp\psi_{2}\ \text{and}\ \chi_{2}\sharp\psi_{1}.}$
Since $s_{2}\not=0$, we have 
\begin{align*}
\chi_{1}\sharp(\psi_{1}\sharp\pi_{2}) & =-(\psi_{1}\sharp\pi_{1})\sharp(\psi_{1}\sharp\pi_{2})\\
 & =\psi_{1}^{\sharp}\sharp(\pi_{1}\sharp\pi_{2})-T(\psi_{1}^{\sharp},\pi_{1})\pi_{2}-T(\psi_{1}^{\sharp},\pi_{2})\pi_{1}
\end{align*}
by the formula \cite[(11)]{R} and \cite[(25)]{R}. Then we obtain
the presentation of $\chi_{1}\sharp\psi_{2}$ as in (\ref{eq:xy})
using (\ref{eq:pi2psi}). 
\end{proof}

\subsection{$\sharp$-mapping and cubic form\label{subsec:-mapping-and-cubicformUp}}

Now we can write down $\sigma^{\sharp}$ for any $\sigma\in J$ explicitly.
To state the result, we prepare some notation as follows:

\begin{align*}
\bm{s} & :=\,\empty^{t}\!\left(\begin{array}{cccc}
s_{1} & s_{2} & s_{3} & s_{4}\end{array}\right),~\bm{t}:=\,\empty^{t}\!\left(\begin{array}{cccc}
t_{1} & t_{2} & t_{3} & t_{4}\end{array}\right).\\
\Delta_{ij} & :\text{the}\ (i,j)\text{-minor of the matrix}\ \left(\begin{array}{cccc}
y_{1} & y_{2} & zx_{1} & zx_{2}\\
x_{1} & x_{2} & y_{1} & y_{2}
\end{array}\right).\\
\bm{\Delta} & :=\,\empty^{t}\!\left(\begin{array}{cccc}
\Delta_{12} & \Delta_{13} & 2\Delta_{23} & \Delta_{24}\end{array}\right).\\
M_{1} & :=\left(\begin{array}{cccc}
0 & -x_{1}^{2} & -2x_{1}x_{2} & -x_{2}^{2}\\
 & 0 & -2x_{1}y_{1} & -x_{1}y_{2}-x_{2}y_{1}\\
 &  & 0 & -2x_{2}y_{2}\\
 &  &  & 0
\end{array}\right),\\
M_{3} & :=\left(\begin{array}{cccc}
0 & x_{1}y_{1} & x_{1}y_{2}+x_{2}y_{1} & x_{2}y_{2}\\
 & 0 & y_{1}^{2}+zx_{1}^{2} & y_{1}y_{2}+zx_{1}x_{2}\\
 &  & 0 & y_{2}^{2}+zx_{2}^{2}\\
 &  &  & 0
\end{array}\right),\\
 & \text{which are 4\ensuremath{\times}4 skew-symmetric matrices. }\\
M & :=\text{\text{the \ensuremath{5\times5}} skew symmetric matrix \ensuremath{(m_{ij})} with}\\
 & m_{12}=x_{1},\,m_{13}=x_{2},\,m_{14}=y_{1},\,m_{15}=y_{2},\\
 & m_{23}=-p_{1},\,m_{24}=p_{4}t_{4}-p_{2}s_{4},\,m_{25}=p_{3}+(p_{2}s_{3}-p_{4}t_{3}),\\
 & m_{34}=-p_{3}+(p_{2}s_{3}-p_{4}t_{3}),\,m_{35}=p_{4}t_{2}-p_{2}s_{2},\\
 & m_{45}=-zp_{1}+p_{2}s_{1}-p_{4}t_{1}.\\
{\rm Pf}_{k} & :\text{the Pfaffian of the \ensuremath{4\times4}submatrix of \ensuremath{M} obtained by deleleting \ensuremath{k}th row and culumn.}
\end{align*}

The following result follows from Propositions \ref{prop:vpi}, \ref{prop:pichipsi},
and \ref{prop:J12sha}:
\begin{cor}
\label{cor:Upsilonshar}In the situation of Corollary \ref{cor:sharpformula},
we further write
\begin{align*}
\sigma_{0} & =p_{1}\pi_{1}+p_{2}\pi_{2}+p_{3}\pi_{3}+p_{4}\pi_{4},\\
\sigma_{\nicefrac{1}{2}} & =x_{1}\chi_{1}+x_{2}\chi_{2}+y_{1}\psi_{1}+y_{2}\psi_{2},
\end{align*}
and 
\begin{align*}
(\sigma_{\nicefrac{1}{2}})^{\sharp}+u(\upsilon\sharp\sigma_{0}) & =p_{1}^{\sharp}\pi_{1}+p_{2}^{\sharp}\pi_{2}+p_{3}^{\sharp}\pi_{3}+p_{4}^{\sharp}\pi_{4},\\
\sigma_{0}\sharp\sigma_{\nicefrac{1}{2}} & =x_{1}^{\sharp}\chi_{1}+x_{2}^{\sharp}\chi_{2}+y_{1}^{\sharp}\psi_{1}+y_{2}^{\sharp}\psi_{2}.
\end{align*}
We also set 
\[
u^{\sharp}:=S(\sigma_{0}).
\]
 It holds that {\small{}
\begin{equation}
\begin{cases}
p_{1}^{\sharp}=-up_{1}-\,\empty^{t}\bm{t}M_{1}\bm{s},\\
p_{2}^{\sharp}=-up_{2}+\,\empty^{t}\bm{t}\bm{\Delta},\\
p_{3}^{\sharp}=up_{3}+\,\empty^{t}\bm{t}M_{3}\bm{s}\\
p_{4}^{\sharp}=-up_{4}+\,\empty^{t}\bm{s}\bm{\Delta},\\
x_{1}^{\sharp}={\rm Pf}_{5},\,x_{2}^{\sharp}={\rm Pf}_{4},\\
y_{1}^{\sharp}={\rm Pf}_{3},\,y_{2}^{\sharp}={\rm Pf}_{2}.\\
u^{\sharp}=-{\rm Pf}_{1},
\end{cases}\label{eq:sigmasharp}
\end{equation}
}The cubic form defining $J$ is 
\begin{align}
N_{\Upsilon} & :=\nicefrac{2}{3}[\left(2z\mathsf{S}\!\mathsf{T}_{\!23}\,x_{1}+z\mathsf{S}\!\mathsf{T}_{\!24}\,x_{2}+\mathsf{S}\!\mathsf{T}_{\!12}\,y_{1}+\mathsf{S}\!\mathsf{T}_{\!13}\,y_{2}\right)(x_{1}^{\sharp})\label{eq:NUp}\\
 & +\nicefrac{2}{3}\left(z\mathsf{S}\!\mathsf{T}_{\!24}\,x_{1}+2z\mathsf{S}\!\mathsf{T}_{\!34}\,x_{2}+\mathsf{S}\!\mathsf{T}_{\!13}\,y_{1}+\mathsf{S}\!\mathsf{T}_{\!14}\,y_{2}\right)(x_{2}^{\sharp})\nonumber \\
 & +\nicefrac{2}{3}\left(\mathsf{S}\!\mathsf{T}_{\!12}\,x_{1}+\mathsf{S}\!\mathsf{T}_{\!13}\,x_{2}+2\mathsf{S}\!\mathsf{T}_{\!23}\,y_{1}+\mathsf{S}\!\mathsf{T}_{\!24}\,y_{2}\right)(y_{1}^{\sharp})\nonumber \\
 & +\nicefrac{2}{3}\left(\mathsf{S}\!\mathsf{T}_{\!13}\,x_{1}+\mathsf{S}\!\mathsf{T}_{\!14}\,x_{2}+\mathsf{S}\!\mathsf{T}_{\!24}\,y_{1}+2\mathsf{S}\!\mathsf{T}_{\!34}\,y_{2}\right)(y_{2}^{\sharp})\nonumber \\
 & +\nicefrac{1}{3}\,u(u^{\sharp})+\nicefrac{1}{3}(2zp_{1}-s_{1}p_{2}+t_{1}p_{4})(p_{1}^{\sharp})+\nicefrac{2}{3}\,p_{3}(p_{3}^{\sharp})\nonumber \\
 & +\nicefrac{1}{3}\left(-s_{1}p_{1}+2(s_{3}^{2}-s_{2}s_{4})p_{2}+(s_{4}t_{2}-2s_{3}t_{3}+s_{2}t_{4})p_{4}\right)(p_{2}^{\sharp})\nonumber \\
 & +\nicefrac{1}{3}\left(t_{1}p_{1}+(s_{4}t_{2}-2s_{3}t_{3}+s_{2}t_{4})p_{2}+2(t_{3}^{2}-t_{2}t_{4})p_{4}\right)(p_{4}^{\sharp}),\nonumber 
\end{align}
where we set $\mathsf{S}\!\mathsf{T}_{\!ij}:=s_{i}t_{j}-s_{j}t_{i}$.
We note that the denominator $3$ in (\ref{eq:NUp}) only appears
in the expression and is canceled such that r.h.s. of (\ref{eq:NUp})
has only integer coefficients (thus, it works also when the characteristic
of $\mathsf{k}$ is $3$).
\end{cor}

Similarly to \cite[Prop.3.5]{Tak5}, we also have the following converse
statement, for which we do not need to assume that $\mathsf{k}$ is
algebraically closed:
\begin{prop}
\label{Prop:reveseJ2}Let $J$ be a $9$-dimensional $\mathsf{k}$-vector
space with coordinates $u$, $x_{i},y_{i}$ $(i=1,2)$, $p_{j}\,(1\leq j\leq4)$.
Let $N_{\Upsilon}$ be the cubic form (\ref{eq:NUp}) on $J$ as in
the statement of Corollary \ref{cor:Upsilonshar} with some constants
$z,s_{k},t_{k}\,(1\leq k\leq4).$ We also define $\sigma^{\sharp}$
for $\sigma\in J$ as in the statement of Corollary \ref{cor:Upsilonshar}.
The vector space $J$ has the structure of the quadratic Jordan algebra
of the cubic form $N_{\Upsilon}$.
\end{prop}

\section{\textbf{Affine variety $\Upsilon_{\mA}^{14}$}\label{sec:Affine-varietyUp}}

In this section, we assume for simplicity that $\mathsf{k}=\mC$.

\subsection{Definition and some basic properties of{\small{} $\Upsilon_{\mA}^{14}$\label{subsec:Definition-of-theUpsilon}}}
\begin{defn}
\label{def:Upsilon} Let $\mA_{\Upsilon}$ be the affine $18$-space
with coordinates 
\begin{equation}
u,x_{i},y_{i}(i=1,2),p_{j}\,(1\leq j\leq4).\label{eq:UpCoord}
\end{equation}
 We define a morphism $\mA_{\Upsilon}\to\mA_{\Upsilon}$ by $\mA_{\Upsilon}\ni\sigma\mapsto\sigma^{\sharp}\in\mA_{\mathscr{H}}$
as in Corollary \ref{cor:sharpformula}. In $\mA_{\Upsilon}$, we
define $\Upsilon_{\mA}^{14}$ to be the affine scheme defined by the
condition 
\[
\sigma^{\sharp}=0,
\]
 namely, the vanishing of the r.h.s. of (\ref{eq:sigmasharp}).
\end{defn}

Surprisingly, the affine scheme $\Upsilon_{\mA}^{14}$ coincides with
the one constructed by R.Taylor in her PhD thesis \cite[Subsec. 3.1.2]{Tay}
(as for the notation, we follow our previous paper \cite[Sect. 7]{Tak2}).
Her construction is based on the theory of unprojection. She studies
the properties of $\Upsilon_{\mA}^{14}$ much in details. Here we
summarize her results:
\begin{thm}
\label{thm:PropUp}Let $S_{\Pi}$ be the polynomial ring over $\mC$
whose variables are (\ref{eq:UpCoord}). Let $I_{\Pi}$ be the ideal
of the polynomial ring $S_{\Pi}$ generated by the $9$ equations
of $\Upsilon_{\mathbb{A}}^{14}$. The following assertions hold:

\begin{enumerate}[$(1)$]

\item We give nonnegative weights for coordinates of $S_{\Upsilon}$
such that all the equations of $\Upsilon_{\mA}^{14}$ are weighted
homogeneous, and we denote by $w(*)$ the weight of the monomial $*$.
The ideal of $I_{\Upsilon}$ of $\Upsilon_{\mathbb{A}}^{14}$ has
the following graded minimal $S_{\Pi}$-free resolution: 
\[
0\leftarrow R_{\Pi}\leftarrow P_{0}\leftarrow P_{1}\leftarrow P_{2}\leftarrow P_{3}\leftarrow P_{4}\leftarrow0.
\]

\item $I_{\Upsilon}$ is a Gorenstein ideal of codimension $4$. 

\item $\Upsilon_{\mA}^{14}$ is irreducible and reduced, thus $I_{\Upsilon}$
is a prime ideal.

\end{enumerate}
\end{thm}

She also constructs several quasi-smooth $\mQ$-Fano $3$-folds from
$\Upsilon_{\mA}^{14}$, which are expected to be prime. In the sequel,
we obtain further properties of $\Upsilon_{\mA}^{14}$, based on which
we will show that Taylor's example of $\mQ$-Fano 3-folds are actually
prime in a future.

\subsection{Weights for variables and equations\label{subsec:Weights-for-variables-Up}}

We assign weights for variables of the polynomial ring $S_{\Upsilon}$
such that all the $9$ equations of $\Upsilon_{\mA}^{14}$ are homogeneous.
Moreover, we assume that all the variables are not zero allowing some
of them are constants. Then it is easy to derive the following relations
between the weights of variables of $S_{\Upsilon}$:

\begin{align*}
w(y_{1}) & =w(x_{1})+\nicefrac{w(z)}{2},\,w(y_{2})=w(x_{2})+\nicefrac{w(z)}{2},\\
w(p_{3}) & =w(p_{1})+\nicefrac{w(z)}{2},\\
w(s_{1}) & =w(p_{1})-w(p_{2})+w(z),\\
w(s_{2}) & =w(p_{1})-w(p_{2})+w(x_{2})-w(x_{1})+\nicefrac{w(z)}{2},\\
w(s_{3}) & =w(p_{1})-w(p_{2})+\nicefrac{w(z)}{2},\\
w(s_{4}) & =w(p_{1})-w(p_{2})+w(x_{1})-w(x_{2})+\nicefrac{w(z)}{2,}\\
w(t_{1}) & =w(p_{1})-w(p_{4})+w(z),\\
w(t_{2}) & =w(p_{1})-w(p_{4})+w(x_{2})-w(x_{1})+\nicefrac{w(z)}{2},\\
w(t_{3}) & =w(p_{1})-w(p_{4})+\nicefrac{w(z)}{2},\\
w(t_{4}) & =w(p_{1})-w(p_{4})+w(x_{1})-w(x_{2})+\nicefrac{w(z)}{2,}\\
w(u) & =w(p_{1})-w(p_{2})-w(p_{4})+w(x_{1})+w(x_{2})+\nicefrac{3w(z)}{2}.
\end{align*}

We need this result in the proof of Lemma \ref{lem:prime-Up}.

\subsection{${\rm GL}_{2}\times{\rm GL}_{2}$-action on $\Upsilon_{\mA}^{14}$}
\begin{prop}
\label{prop:For-any-elementGL2GL2Up}For any element $g\in\mathrm{GL}_{2}$,
we set {\small{}
\[
\hat{g}:=(1/\det g)\left(\begin{array}{ccc}
d^{2} & -cd & c^{2}\\
-2bd & bc+ad & -2ac\\
b^{2} & -ab & a^{2}
\end{array}\right),
\]
and denote by }$g^{\dagger}$ the adjoint matrix of $g$. The group
$\mathrm{GL_{2}\times{\rm GL}_{2}}$ acts on the affine variety $\Upsilon_{\mA}^{14}$
by the following rule:\begin{itemize}

\item For any $g\in\mathrm{GL}_{2}$, we set 

\begin{align*}
\left(\begin{array}{ccc}
t_{2} & t_{3} & t_{4}\\
s_{2} & s_{3} & s_{4}
\end{array}\right) & \mapsto\left(\begin{array}{ccc}
t_{2} & t_{3} & t_{4}\\
s_{2} & s_{3} & s_{4}
\end{array}\right)\widehat{g},\\
\left(\begin{array}{cc}
x_{1} & x_{2}\\
y_{1} & y_{2}
\end{array}\right) & \mapsto\left(\begin{array}{cc}
x_{1} & x_{2}\\
y_{1} & y_{2}
\end{array}\right)g,\\
\left(\begin{array}{cccc}
p_{1} & p_{2} & p_{3} & p_{4}\end{array}\right) & \mapsto\det g\left(\begin{array}{cccc}
p_{1} & p_{2} & p_{3} & p_{4}\end{array}\right).
\end{align*}

\item For any $h\in\mathrm{GL}_{2}$, we set 

\begin{align*}
\left(\begin{array}{cccc}
t_{1} & t_{2} & t_{3} & t_{4}\\
s_{1} & s_{2} & s_{3} & s_{4}
\end{array}\right) & \mapsto h\left(\begin{array}{cccc}
t_{1} & t_{2} & t_{3} & t_{4}\\
s_{1} & s_{2} & s_{3} & s_{4}
\end{array}\right),\,\left(\begin{array}{c}
p_{2}\\
p_{4}
\end{array}\right)\mapsto h\left(\begin{array}{c}
p_{2}\\
p_{4}
\end{array}\right),\\
p_{1} & \mapsto\det h\,p_{1},p_{3}\mapsto\det h\,p_{2}.
\end{align*}

\end{itemize}
\end{prop}

\begin{proof}
We can check directly that the above rule define two actions of ${\rm GL}_{2}$
on $\Upsilon_{\mA}^{14}$ and the two actions mutually commute. Therefore
we have the ${\rm GL}_{2}\times{\rm GL_{2}}$-action on $\Upsilon_{\mA}^{14}$.
\end{proof}

\subsection{Charts of $\Upsilon_{\mA}^{14}$\label{subsec:Charts-ofUp}}

Let $\mathbb{A}_{Z\!S\!T}$ be the affine $9$-space with the coordinates
$z,s_{i},t_{i}\,(1\leq i\leq4)$. In this section, we consider $\mathbb{A}_{Z\!S\!T}$
is contained in the ambient affine space $\mA_{\Upsilon}$ of $\Upsilon_{\mathbb{A}}^{14}$.
Then we see that $\mathbb{A}_{Z\!S\!T}\subset\Upsilon_{\mathbb{A}}^{14}$. 

For a coordinate $*$, we call the open subset of $\Upsilon_{\mA}^{14}$
with $*\not=0$ \textit{the $*$-chart}. In this subsection, we describe
the $*$-chart with $*=u$, $p_{i}\,(1\leq i\leq4)$, $x_{1},x_{2}$,
$y_{1},y_{2}$.

\vspace{3pt}

\noindent $\underline{\text{\ensuremath{u}-chart}}$: By the equation
$p_{1}^{\sharp}=p_{2}^{\sharp}=p_{3}^{\sharp}=p_{4}^{\sharp}=0$,
we may erase the coordinates $p_{1},p_{2},p_{3},p_{4}$. Therefore
the $u$-chart is smooth.

\noindent $\underline{\text{\ensuremath{p_{1}}-chart}}$: By the
equation $x_{1}^{\sharp}=x_{2}^{\sharp}=p_{1}^{\sharp}=u^{\sharp}=0$,
we may erase the coordinates $u,y_{1},y_{2},z$. Therefore the $p_{1}$-chart
is smooth.

\vspace{3pt}

Each of the equations of the $p_{2}$-, $p_{4}$-, $x_{1}$-, $x_{2}$-charts
is presented as five Pfaffians $4\times4$ principal submatrix of
the $5\times5$ skew-symmetric matrix $(m_{ij})$ in the following
form: 

\vspace{3pt}

\noindent $\underline{\text{\ensuremath{p_{2}}-, \ensuremath{p_{4}}-chart}}$:
The $5\times5$ skew-symmetric matrix is $M$ introduced before Corollary
\ref{cor:Upsilonshar} with $p_{2}=1$ and $p_{4}=1$ respectively.

\vspace{3pt}

By these descriptions, each of the $p_{2}$-, $p_{4}$-, $x_{1}$-,
$x_{2}$-charts has $c({\rm G}(2,5))$-singularities along $\{m_{ij}=0\,(1\leq i<j\leq5)\}$. 

\vspace{3pt}

\noindent $\underline{\text{\ensuremath{x_{1}}-chart}}$:

\begin{align*}
m_{12} & =u,\,m_{13}=s_{2}+s_{3}x_{2}+s_{4}x{}_{2}^{2},\,m_{14}=-t_{2}-t_{3}x_{2}-t_{4}x{}_{2}^{2},\\
m_{15} & =y_{2}-x_{2}y_{1},\,m_{23}=-s_{1}+s_{3}y_{1}+s_{4}x_{2}y_{1}+s_{4}y_{2},\\
m_{24} & =t_{1}-t_{3}y_{1}-t_{4}y_{2}-t_{4}x_{2}y_{1},\,m_{25}=z-y{}_{1}^{2},\\
m_{34} & =p_{1},\,m_{35}=-p_{4},\,m_{45}=p_{2}.
\end{align*}

\vspace{3pt}

\noindent $\underline{\text{\ensuremath{x_{2}}-chart}}$: 
\begin{align*}
m_{12} & =u,\,m_{13}=s_{4}+s_{3}x_{1}+s_{2}x{}_{1}^{2},\,m_{14}=-t_{4}-t_{3}x_{1}-t_{2}x{}_{1}^{2},\\
m_{15} & =y_{1}-x_{1}y_{2},\,m_{23}=s_{1}+s_{2}y_{1}+s_{3}y_{2}+s_{2}x_{1}y_{2},\\
m_{24} & =-t_{1}-t_{2}y_{1}-t_{3}y_{2}-t_{2}x_{1}y_{2},\,m_{25}=z-y{}_{2}^{2},\\
m_{34} & =-p_{1},\,m_{35}=-p_{4},\,m_{45}=p_{2}.
\end{align*}

\vspace{3pt}

\noindent $\underline{\text{\ensuremath{y_{1}}-chart}}$: Since we
already describe the $x_{1}$-chart, we have only to investigate the
$y_{1}$-chart near $x_{1}=0$. Then, by the equation $p_{2}^{\sharp}=p_{4}^{\sharp}=x_{1}^{\sharp}=y_{1}^{\sharp}=0$,
we may erase the coordinates $p_{1},p_{3},s_{2},t_{2}$ allowing $zx_{1}^{2}-y_{1}^{2}$
as denominators. Therefore the $y_{1}$-chart is smooth near $x_{1}=0$. 

\vspace{3pt}

\noindent $\underline{\text{\ensuremath{y_{2}}-chart}}$: Similarly
to the case of the $y_{1}$-chart, we can check that the $y_{2}$-chart
is smooth near $x_{2}=0$. 

\vspace{3pt}

\noindent $\underline{\text{\ensuremath{p_{3}}-chart}}$: By the
equation $u^{\sharp}=0$, we see that the $p_{3}$-chart is contained
in one of the charts considered above.

\vspace{5pt}

Let $\Delta$ be the union of the singular locus of the $*$-chart
with $*=u,p_{i}\,(1\leq i\leq4),x_{1},x_{2},y_{1},y_{2}.$ By the
above descriptions of the charts of $\Upsilon_{\mA}^{14}$, we have
the following:
\begin{prop}
\label{prop:The-singular-locusUp} 

The variety $\Upsilon_{\mA}^{14}$ has $c({\rm G}(2,5))$-singularities
along $\Delta$, where we call a singularity isomorphic to the vertex
of the cone over ${\rm G}(2,5)$ a \textup{$c({\rm G}(2,5))$-singularity. }
\end{prop}

\begin{cor}
\label{cor:The-variety-normal UP}The variety $\Upsilon_{\mA}^{14}$
is normal.
\end{cor}

\begin{proof}
This follows from Theorem \ref{thm:PropUp} (2) and (3), and Proposition
\ref{prop:The-singular-locusUp}.
\end{proof}
We will describe the singularities of $\Upsilon_{\mA}^{14}$ along
$\mA_{ST}$ in Section \ref{sec:Singularities-ofUpPi}.

\subsection{Factoriality of $\Upsilon_{\mathbb{A}}^{14}$\label{sec:Factoriality-ofUp}}

In this section, we show that the affine coordinate ring $R_{\Upsilon}$
of $\Upsilon_{\mathbb{A}}^{14}$ is a UFD. 
\begin{lem}
\label{lem:prime-Up} We denote by $\bar{p_{1}}$ the image of $p_{1}$
in $R_{\Upsilon}$. The element of $\bar{p_{1}}\in R_{\Upsilon}$
is a prime element, equivalently, $\Upsilon_{\mathbb{A}}^{14}\cap\{p_{1}=0\}$
is irreducible and reduced. 
\end{lem}

\begin{proof}
First we show that the assertion follows from the normality of $\Upsilon_{\mathbb{A}}^{14}\cap\{p_{1}=0\}.$
Note that the polynomial ring $S_{\Upsilon}$ is positively graded
with some weights of variables such that $I_{\Upsilon}$ is homogeneous
and $p_{1}$ is semi-invariant (cf.~Subsection \ref{subsec:Weights-for-variables-Up}).
Then the ring $R_{\Upsilon}$ is also graded with these weights. In
this situation, $p_{1}$ defines an ample divisor on the corresponding
weighted projectivization $\Upsilon_{\mathbb{P}}^{13}$ of $\Upsilon_{\mathbb{A}}^{14}$.
Note that $\Upsilon_{\mathbb{P}}^{13}$ is irreducible and normal
since so is $\Upsilon_{\mathbb{A}}^{14}$ by Theorem \ref{thm:PropUp}
and Corollary \ref{cor:The-variety-normal UP}, and $\Upsilon_{\mathbb{P}}^{13}$
is a geometric $\mathbb{C^{*}}$-quotient of $\Upsilon_{\mathbb{A}}^{14}$.
By a general property of an ample divisor on an irreducible projective
normal variety, the Cartier divisor $\Upsilon_{\mathbb{P}}^{13}\cap\{p_{1}=0\}$
is connected and hence so is $\Upsilon_{\mathbb{A}}^{14}\cap\{p_{1}=0\}$.
Therefore the normality of $\Upsilon_{\mathbb{A}}^{15}\cap\{p_{1}=0\}$
implies that $\Upsilon_{\mathbb{A}}^{14}\cap\{p_{1}=0\}$ is irreducible
and reduced.

Now we show that $\Upsilon_{\mathbb{A}}^{14}\cap\{p_{1}=0\}$ is normal.
Since it is Gorenstein by Theorem \ref{thm:PropUp} (2) and $\bar{p_{1}}$
is not a zero divisor by ibid. (3), we have only to show that $\Upsilon_{\mathbb{A}}^{14}\cap\{p_{1}=0\}$
is regular in codimension 1. By the Pfaffians presentations of the
$x_{1}$-, $x_{2}$-, $p_{2}$- and $p_{4}$-charts as in Subsection
\ref{subsec:Charts-ofUp}, we see that $\Upsilon_{\mathbb{A}}^{14}\cap\{p_{1}=0\}$
is regular in codimension 1 in these charts. Moreover, we can check
the codimension of $\Upsilon_{\mathbb{A}}^{14}\cap\{p_{1}=x_{1}=x_{2}=p_{2}=p_{4}=0\}$
in $\Upsilon_{\mathbb{A}}^{14}\cap\{p_{1}=0\}$ is greater than 1.
Thus $\Upsilon_{\mathbb{A}}^{14}\cap\{p_{1}=0\}$ is regular in codimension
1.
\end{proof}
In the proof of Lemma \ref{lem:prime-Up}, we have proved the following: 
\begin{cor}
\label{cor:boundarynormal.-Up}$\Upsilon_{\mathbb{A}}^{14}\cap\{p_{1}=0\}$
is normal. 
\end{cor}

\begin{prop}
\label{prop:UFD-Up}The affine coordinate ring $R_{\Upsilon}$ of
$\Upsilon_{\mathbb{A}}^{14}$ is a UFD. 
\end{prop}

\begin{proof}
Using the description of the $p_{1}$-chart as in Subsection \ref{subsec:Charts-ofUp}
and Lemma \ref{lem:prime-Up}, the proof of \cite[Prop.4.9]{Tak5}
work verbatim.
\end{proof}
\begin{prop}
\label{prop:Pic1-1}Let $\Upsilon_{\mathbb{P}}^{13}$ be the weighted
projectivization of $\Upsilon_{\mathbb{\mathbb{A}}}^{14}$ with some
positive weights of coordinates. It holds that

\vspace{3pt}

\noindent $(1)$ any prime Weil divisors on $\Upsilon_{\mathbb{P}}^{13}$
are the intersections between $\Upsilon_{\mathbb{P}}^{13}$ and weighted
hypersurfaces. In particular, $\Upsilon_{\mathbb{P}}^{13}$ is $\mathbb{Q}$-factorial
and have Picard number one, and

\noindent $(2)$ let $X$ be a quasi-smooth threefold such that $X$
is a codimension $10$ weighted complete intersection in $\Upsilon_{\mathbb{P}}^{13}$.
Assume moreover that $X\cap\{p_{1}=0\}$ is a prime divisor. Then
any prime Weil divisor on $X$ is the intersection between $X$ and
a weighted hypersurface. In particular, $X$ is $\mathbb{Q}$-factorial
and has Picard number one. 
\end{prop}

\begin{proof}
The proofs of \cite[Cor.4.10]{Tak5} work verbatim using Propositions
\ref{prop:UFD-Up}. 
\end{proof}

\subsection{Relation between $\Upsilon_{\mathbb{A}}^{14}$ and $\mathscr{H}_{\mathbb{A}}^{13}$. }

We establish the following relation between $\Upsilon_{\mathbb{A}}^{14}$
and $\mathscr{H}_{\mathbb{A}}^{13}$. We define some notation. Let
$\mathbb{A}_{Z\!S\!T}$ be the affine $9$-space as in the section
\ref{subsec:Charts-ofUp} (but we do not consider here that $\mathbb{A}_{Z\!S\!T}\subset\Upsilon_{\mathbb{A}}^{14}$).
Let $\mathbb{A}_{S\!T}$ be the affine space with the coordinates
$s_{i},t_{i}\,(1\leq i\leq4)$ and $\mathbb{A}_{A}$ the affine line
with the coordinate $A$. We set $\mathbb{\widetilde{A}}_{Z\!S\!T}=\mathbb{A}_{A}\times\mathbb{A}_{S\!T}$.
Let $b\colon\mathbb{\widetilde{A}}_{Z\!S\!T}\to\mathbb{A}_{Z\!S\!T}$
be the morphism defined with $z=A^{2}$ and the remaining coordinates
being unchanged. The morphism $b$ is a finite morphism of degree
two. Let $\widetilde{\Upsilon}_{\mathbb{A}}^{14}:=\Upsilon_{\mathbb{A}}^{14}\times_{\mathbb{A}_{Z\!S\!T}}\widetilde{\mathbb{A}}_{Z\!S\!T}$,
and $b_{\Upsilon}\colon\widetilde{\Upsilon}_{\mathbb{A}}^{15}\to\Upsilon_{\mathbb{A}}^{15}$
and $p_{\Upsilon}\colon\widetilde{\Upsilon}_{\mathbb{A}}^{14}\to\mathbb{\widetilde{A}}_{Z\!S\!T}$
the naturally induced morphisms. Let $p_{\mathscr{H}}\colon\mathscr{H}_{\mathbb{A}}^{13}\to\mathbb{A_{\mathsf{P}}}$
be the natural projection. Let $U_{A}$ be the open subset $\{A\not=0\}$
of $\mathbb{A}_{A}$. Note that the morphism $b$ is unramified on
$U_{A}\times\mathbb{A}_{S\!T}\subset\widetilde{\mathbb{A}}_{Z\!S\!T}$. 
\begin{prop}
\label{prop:isomHUp}There exists an isomorphism from $U_{A}\times\mathbb{A}_{S\!T}$
to $U_{A}\times\mathbb{A}_{\mathsf{P}}$, and an isomorphism from
$\widetilde{\Upsilon}_{\mathbb{A}}^{14}\cap\{A\not=0\}$ to $U_{A}\times\mathscr{H}_{\mathbb{A}}^{13}$
which fit into the following commutative diagram:{\small{}
\[
\xymatrix{\widetilde{\Pi}_{\mathbb{A}}^{15}\cap\{A\not=0\}\ar[r]\ar[d]_{p_{\Upsilon}} & U_{A}\times\mathscr{H}_{\mathbb{A}}^{13}\ar[d]^{p_{\mathscr{H}}\times{\rm id}}\\
U_{A}\times\mathbb{A}_{S\!T}\ar[r] & U_{A}\times\mathbb{A}_{\mathsf{P}}
}
\]
}{\small\par}
\end{prop}

\begin{proof}
We can directly check that a morphism from $U_{A}\times\mathbb{A}_{S\!T}$
to $U_{A}\times\mathbb{A}_{\mathsf{P}}$, and a morphism from $\widetilde{\Upsilon}_{\mathbb{A}}^{14}\cap\{A\not=0\}$
to $U_{A}\times\mathscr{H}_{\mathbb{A}}^{13}$ can be defined by the
following equalities and they are actually isomorphisms with the desired
properties.
\begin{align*}
p_{111} & =t_{4},p_{112}=s_{4},p_{121}=-\nicefrac{(t_{1}+At_{3})}{2A},p_{211}=-\nicefrac{(t_{1}-At_{3})}{2A},\\
p_{122} & =-\nicefrac{(s_{1}+As_{3})}{2A},p_{212}=-\nicefrac{(s_{1}-As_{3})}{2A},p_{221}=-t_{2},p_{222}=-s_{2},\\
x_{11} & =Ax_{1}+y_{1},x_{21}=-(Ax_{2}+y_{2}),\\
x_{12} & =-(Ax_{1}-y_{1}),x_{22}=-(Ax_{2}-y_{2}),\\
x_{13} & =p_{2},x_{23}=p_{4},\\
u_{1} & =Ap_{1}+p_{3}+\nicefrac{(-p_{2}s_{1}+p_{4}t_{1})}{2A},\\
u_{2} & =Ap_{1}-p_{3}+\nicefrac{(-p_{2}s_{1}+p_{4}t_{1})}{2A},\\
u_{3} & =-u.
\end{align*}
\end{proof}

\subsection{$\mathbb{P}^{2}\times\mathbb{P}^{2}$-fibration associated to $\Upsilon_{\mathbb{A}}^{14}$}

Note that any equation of $\text{\ensuremath{\Upsilon_{\mathbb{A}}^{14}}}$
is of degree two if we regard the coordinates of $\mathbb{A}_{Z\!S\!T}$
as constants. Therefore, considering the coordinates $x_{1},x_{2},y_{1},y_{2},p_{1},p_{2},p_{3},p_{4},u$
as projective coordinates, we obtain a quasi-projective variety with
the same equations as $\text{\ensuremath{\Upsilon_{\mathbb{A}}^{14}}}$,
which we denote by $\widehat{\Upsilon}$. We also denote by $\rho_{\Upsilon}\colon\widehat{\Upsilon}\to\mathbb{A}_{Z\!S\!T}$
the natural projection.

We note that the image of $U_{A}\times\mathbb{A}_{S\!T}$ by the morphism
$b$ is the open subset $\{z\not=0\}\times\mathbb{A}_{S\!T}$ of $\mathbb{A}_{Z\!S\!T}$.
As in \cite[Subsec.3.6]{Tak5}, we denote by $D_{\mathscr{H}}$ the
Cayley's hyperdeterminant. Let $\{D'_{\mathscr{H}}=0\}\subset\mA_{Z\!S\!T}$
be the image by $b$ of the pull-back on $U_{A}\times\mathbb{A}_{S\!T}$
of $U_{A}\times\{D_{\mathscr{H}}=0\}$. 
\begin{prop}
\label{prop:P2P2Up}Let $\rho_{\Upsilon}\colon\widehat{\Upsilon}\to\mathbb{A}_{Z\!S\!T}$
be the natural projection. The $\rho_{\Upsilon}$-fibers over points
of $\{D'_{\mathscr{H}}\not=0\}\cap(\{z\not=0\}\times\mathbb{A}_{S\!T})$
are isomorphic to $\mathbb{P}^{2}\times\mathbb{P}^{2}$, and the $\rho_{\Upsilon}$-fibers
over points of a nonempty open subset of $\{D'_{\mathscr{H}}=0\}$
are isomorphic to $\mathbb{P}^{2,2}$, where $\mP^{2,2}$ is defined
in \cite{Fu,Mu1} (see also \cite[Def.4.4]{Tak5}).
\end{prop}

\begin{proof}
Let $\widehat{\Upsilon}'$ be the base change of $\widehat{\Upsilon}$
over $\widetilde{\mathbb{A}}_{Z\!S\!T}$. We have only to describe
the corresponding fibers of $\widehat{\Upsilon}'\to\widetilde{\mathbb{A}}_{Z\!S\!T}$.
Since the isomorphism from $\widetilde{\Upsilon}_{\mathbb{A}}^{14}\cap\{A\not=0\}$
to $\mathscr{H}_{\mathbb{A}}^{13}\times U_{A}$ as in Proposition
\ref{prop:isomHUp} is linear with respect to the coordinates which
are vertical to the fibrations $\rho_{\Upsilon}$ and $\rho_{\mathscr{H}}$,
it descends to an isomorphism from $\widehat{\Upsilon}'\cap\{A\not=0\}$
to the open subset $\widehat{\mathscr{H}}\times U_{A}$. Thus we have
the assertion from \cite[Prop. 4.5]{Tak5}. 
\end{proof}
It is possible to describe all the $\rho_{\Upsilon}$-fibers using
the ${\rm GL}_{2}\times{\rm GL}_{2}$-action on $\Upsilon_{\mA}^{14}$
but it needs more works and pages; Proposition \ref{prop:P2P2Up}
is sufficient for an application in this paper (cf.~Proposition \ref{prop:MoreSing}). 

\section{\textbf{Coordinatization of $J$ with $10$ parameters \label{sec:Coord10paramJ}}}

In this section, we keep considering the quadratic Jordan algebra
$J$ of a cubic form as in the subsection \ref{subsec:Facts-about-a single},\textit{
assuming that $\dim J=9$.}

\subsection{Coordination of $J$ with 8 parameters after a construction due to
Petersson\label{subsec:Coordination-of-Perterson}}

In this subsection, we obtain a coordinatization of $J$ inspired
by a construction due to Petersson \cite{Pe1}, which is different
in nature from those as in \cite[Subsec.3.2]{Tak5} and Subsections
\ref{subsec:Coordination-of-9param}.

We start the construction from two elements $\sigma,\upsilon$ of
$J$ with 

\begin{equation}
\text{T(\ensuremath{\sigma)=T(\upsilon)=0}.}\label{eq:Trace=00003D0}
\end{equation}
We set 
\begin{equation}
\begin{cases}
S_{1}\ :=S(\sigma),\,S_{2}:=S(\upsilon),\,N_{1}:=N(\sigma),\,N_{2}:=N(\upsilon),\\
T_{11}:=T(\sigma,\upsilon),\,T_{12}:=T(\sigma,\upsilon^{\sharp}),\,T_{21}:=T(\sigma^{\sharp},\upsilon),\,T_{22}:=T(\sigma^{\sharp},\upsilon^{\sharp}).
\end{cases}\label{eq:SNT}
\end{equation}
Let
\begin{equation}
V:=\langle\mathfrak{1},\sigma,\sigma^{\sharp},\upsilon,\upsilon^{\sharp},\sigma\sharp\upsilon,\sigma\sharp\upsilon^{\sharp},\sigma^{\sharp}\sharp\upsilon,\sigma^{\sharp}\sharp\upsilon^{\sharp}\rangle.\label{eq:genV}
\end{equation}

The following result due to Pertersson is the starting point of our
consideration:
\begin{prop}[{\cite[Prop.6.6]{Pe1}}]
The vector spacer $V$ is a Jordan subalgebra of $J$. 
\end{prop}

We write an element $x$ of $V$ as
\[
x=x_{1}\mathfrak{1}+x_{2}\sigma+x_{3}\sigma^{\sharp}+x_{4}\upsilon+x_{5}\upsilon^{\sharp}+x_{6}\sigma\sharp\upsilon+x_{7}\sigma\sharp\upsilon^{\sharp}+x_{8}\sigma^{\sharp}\sharp\upsilon+x_{9}\sigma^{\sharp}\sharp\upsilon^{\sharp}.
\]
Then, by (\ref{eq:Trace=00003D0}) and \cite[(13), (16)]{R}, we have
\[
T(x)=3x_{1}+S_{1}x_{3}+S_{2}x_{5}-T_{11}x_{6}-T_{12}x_{7}-T_{21}x_{8}+(S_{1}S_{2}-T_{22})x_{9}.
\]
From this, we may compute directly the bilinear trace $T(x,y)$ for
$x,y\in V$, and then the matrix whose entries are the values of the
bilinear traces between the generators of $V$ as in (\ref{eq:genV})
(we do not write down these since they are lengthy). We set

\begin{align*}
\mathsf{D:=} & 27N_{1}^{2}N{}_{2}^{2}+4N_{1}^{2}S{}_{2}^{3}+4N{}_{2}^{2}S{}_{1}^{3}+6N_{1}N_{2}S_{1}S_{2}T_{11}-S{}_{1}^{2}S{}_{2}^{2}T_{11}^{2}+4N_{1}N_{2}T_{11}^{3}\\
+ & 4N_{1}S_{1}S{}_{2}^{2}T_{12}+4N_{2}S{}_{1}^{2}S_{2}T_{21}+4N_{1}S{}_{2}^{2}T_{11}T_{21}+4N_{2}S{}_{1}^{2}T_{11}T_{12}\\
+ & 4N_{1}S_{2}T_{11}^{2}T_{12}+4N_{2}S_{1}T_{11}^{2}T_{21}-18N_{1}N_{2}T_{12}T_{21}-18N_{1}N_{2}T_{11}T_{22}\\
- & T_{12}^{2}T_{21}^{2}+4N_{1}T_{12}^{3}+4N_{2}T_{21}^{3}+4S{}_{1}^{2}S{}_{2}^{2}T_{22}+2S_{1}S_{2}T_{11}T_{12}T_{21}+2S_{1}S_{2}T_{11}^{2}T_{22}\\
- & 12N_{1}S_{2}T_{12}T_{22}-12N_{2}S_{1}T_{21}T_{22}+2T_{11}T_{12}T_{21}T_{22}\\
+ & 4S_{1}T_{12}^{2}T_{22}+4S_{2}T_{21}^{2}T_{22}-8S_{1}S_{2}T_{22}^{2}-T_{11}^{2}T_{22}^{2}+4T_{22}^{3}.
\end{align*}
We see that the matrix is invertible if and only if $\mathsf{D}\not=0$.
Then, by linear algebra, we obtain the following proposition:
\begin{prop}
\label{prop:J=00003DV}If $\mathsf{D}\not=0$, then $V$ is 9-dimentional,
and hence $J=V.$
\end{prop}

Actually, we show that such $\sigma,\upsilon\in J$ exist:
\begin{prop}
\label{prop:sigmaupsilonExists}There exist $\sigma,\upsilon\in J$
such that $T(\sigma)=T(\upsilon)=0$ and $\mathsf{D}\not=0$.
\end{prop}

\begin{proof}
To show this, we take three complementary idempotents $\upsilon_{1},\upsilon_{2},\upsilon_{3}$
and the associated Peirce decomposition 
\[
J=J_{11}\oplus J_{22}\oplus J_{33}\oplus J_{12}\oplus J_{13}\oplus J_{23}
\]
(see \cite[Subsec.3.1]{Tak5} for a quick review). We choose $\sigma\in J_{12}$
and $\upsilon\in J_{12}\oplus J_{23}$. We write $\upsilon=\upsilon_{12}+\upsilon_{23}$
with $\upsilon_{12}\in J_{12}$ and $\upsilon_{23}\in J_{23}.$ Then,
by \cite[(33)]{R}, we have $\sigma^{\sharp}=S(\sigma)\upsilon_{3}=S_{1}\upsilon_{3}$,
and
\[
\upsilon^{\sharp}=\upsilon_{12}^{\sharp}+\upsilon_{13}^{\sharp}+\upsilon_{12}\sharp\upsilon_{23}=S(\upsilon_{12})\upsilon_{3}+S(\upsilon_{13})\upsilon_{2}+\upsilon_{12}\sharp\upsilon_{23},
\]
and, by \cite[(34)]{R}, $\upsilon_{12}\sharp\upsilon_{23}\in J_{13}$.
We have
\[
T(\sigma)=T(\upsilon)=0
\]
as desired by \cite[(28)]{R}. 

\[
S_{2}=T(\upsilon^{\sharp})=S(\upsilon_{12})+S(\upsilon_{23})
\]
since it holds that $T(\upsilon_{12}\sharp\upsilon_{23})=0$ by \cite[(28)]{R}.
It holds that

\[
N_{1}=N_{2}=0
\]
 by \cite[(32)]{R}, and 

\[
T_{11}=T(\sigma,\upsilon_{12}+\upsilon_{23})=T(\sigma,\upsilon_{12}).
\]
Note that 
\begin{align*}
S(\sigma,\upsilon_{12}): & =S(\sigma+\upsilon_{12})-S(\sigma)-S(\upsilon_{12})\\
 & =T((\sigma+\upsilon_{12})^{\sharp}-\sigma^{\sharp}-\upsilon_{12}^{\sharp})\\
 & =T(\sigma\sharp\upsilon_{12}).
\end{align*}
Therefore, by \cite[(28)]{R}, we have 
\[
T(\sigma,\upsilon_{12})=T(\sigma)T(\upsilon_{12})-T(\sigma\sharp\upsilon_{12})=-S(\sigma,\upsilon_{12})
\]
 By \cite[(33), 2nd line]{R}, we have 
\[
T_{12}=T_{21}=0,
\]
 and 
\[
T_{22}=S_{1}S(\upsilon_{12}).
\]

Therefore we have 
\[
\mathsf{D}=S_{1}^{2}S(\upsilon_{23})^{2}(4S_{1}S(\upsilon_{12})-S(\sigma,\upsilon_{12})^{2}).
\]
By \cite[(33), p.98, 5th line from the bottom]{R}, we may assume
that $S_{1}\not=0$ and $S(\upsilon_{23})\not=0$. We may also choose
$\sigma,\upsilon_{12}$ as a basis of $J_{12}$. We set $S(x_{1}\sigma+x_{2}\upsilon_{12})|_{J_{12}}=ax_{1}^{2}+bx_{1}x_{2}+cx_{2}^{2}$
with $a,b,c\in\mathsf{k}$. Then we have $S(x_{1}\sigma+x_{2}\upsilon_{12},y_{1}\sigma+y_{2}\upsilon_{12})|_{J_{12}}=2ax_{1}y_{1}+b(x_{1}y_{2}+x_{2}y_{1})+2cx_{2}y_{2}$.
Thus we can express
\[
\mathsf{D}=a^{2}S(\upsilon_{23})(4ac-b^{2}).
\]
Now we show that $4ac-b^{2}\not=0$, which finish the proof. By \cite[Cor.4.1.58]{Pe2},
\[
\mathsf{N}:=\{\tau\in J\mid N(\tau)=T(\tau,J)=T(\tau^{\sharp},J)=0\}
\]
is $\text{\{0\}}$ since $J$ is nondegenerate. By \cite[(32), (33), 2nd and 5th lines]{R},
we have 
\[
\mathsf{N}\cap J_{12}=\{\tau\in J_{12}\mid S(\tau,J_{12})=S(\tau)=0\}
\]
since $T(\tau,J_{12})=S(\tau,J_{12})$. Therefore, $\mathsf{N}\cap J_{12}=\{0\}$
implies that $4ac-b^{2}\not=0$.
\end{proof}
Hereafter we choose $\sigma,\upsilon\in J$ as in Proposition \ref{prop:sigmaupsilonExists}
(then the the generators of $V$ as in (\ref{eq:genV}) form a basis
of $V$). By direct computations, we can write the $\sharp$-products
between the generators of $V$ as in (\ref{eq:genV}) as linear combinations
of the generators of $V$ with polynomials of $S_{1},\dots,T_{22}$
as coefficients. We have the following table of the $\sharp$-products,
where $*_{i}$ mean the results are presented below the table :
\begin{flushleft}
{\tiny{}}%
\begin{tabular}{|c|c|c|c|c|c|c|c|c|c|}
\hline 
 & {\tiny{}$\mathfrak{1}$} & {\tiny{}$\sigma$} & {\tiny{}$\sigma^{\sharp}$} & {\tiny{}$\upsilon$} & {\tiny{}$\upsilon^{\sharp}$} & {\tiny{}$\sigma\sharp\upsilon$} & {\tiny{}$\sigma\sharp\upsilon^{\sharp}$} & {\tiny{}$\sigma^{\sharp}\sharp\upsilon$} & {\tiny{}$\sigma^{\sharp}\sharp\upsilon^{\sharp}$}\tabularnewline
\hline 
\multirow{2}{*}{{\tiny{}$\mathfrak{1}$}} & \multirow{2}{*}{{\tiny{}$2\cdot\mathfrak{1}$}} & \multirow{2}{*}{{\tiny{}$-\sigma$}} & {\tiny{}$S_{1}\cdot\mathfrak{1}$} & \multirow{2}{*}{{\tiny{}$-\upsilon$}} & {\tiny{}$S_{2}\cdot\mathfrak{1}$} & {\tiny{}$-T_{11}\cdot\mathfrak{1}$} & {\tiny{}$-T_{12}\cdot\mathfrak{1}$} & {\tiny{}$-T_{21}\cdot\mathfrak{1}$} & {\tiny{}$(S_{1}S_{2}-T_{22})\mathfrak{1}$}\tabularnewline
 &  &  & {\tiny{}$-\sigma^{\sharp}$} &  & {\tiny{}$-\upsilon^{\sharp}$} & {\tiny{}$-\sigma\sharp\upsilon$} & {\tiny{}$-\sigma\sharp\upsilon^{\sharp}$} & {\tiny{}$-\sigma^{\sharp}\sharp\upsilon$} & {\tiny{}$-\sigma^{\sharp}\sharp\upsilon^{\sharp}$}\tabularnewline
\hline 
\multirow{3}{*}{{\tiny{}$\sigma$}} & \multirow{23}{*}{} & \multirow{3}{*}{{\tiny{}$2\sigma^{\sharp}$}} & {\tiny{}$-N_{1}\cdot\mathfrak{1}$} & \multirow{3}{*}{{\tiny{}$\sigma\sharp\upsilon$}} & \multirow{3}{*}{{\tiny{}$\sigma\sharp\upsilon^{\sharp}$}} & {\tiny{}$-T_{21}\cdot\mathfrak{1}+T_{11}\sigma$} & {\tiny{}$(S_{1}S_{2}-T_{22})\cdot\mathfrak{1}$} & {\tiny{}$N_{1}\upsilon$} & {\tiny{}$N_{1}\upsilon^{\sharp}$}\tabularnewline
 &  &  & {\tiny{}$-S_{1}\sigma$} &  &  & {\tiny{}$-S_{1}\upsilon-\sigma^{\sharp}\sharp\upsilon$} & {\tiny{}$+T_{12}\sigma-S_{2}\sigma^{\sharp}$} & {\tiny{}$+T_{11}\sigma^{\sharp}$} & {\tiny{}$+T_{12}\sigma^{\sharp}$}\tabularnewline
 &  &  &  &  &  &  & {\tiny{}$-S_{1}\upsilon^{\sharp}-\sigma^{\sharp}\sharp\upsilon^{\sharp}$} &  & \tabularnewline
\cline{1-1} \cline{3-10} \cline{4-10} \cline{5-10} \cline{6-10} \cline{7-10} \cline{8-10} \cline{9-10} \cline{10-10} 
\multirow{3}{*}{{\tiny{}$\sigma^{\sharp}$}} &  & \multirow{20}{*}{} & \multirow{3}{*}{{\tiny{}$2N_{1}\sigma$}} & \multirow{3}{*}{{\tiny{}$\sigma^{\sharp}\sharp\upsilon$}} & \multirow{3}{*}{{\tiny{}$\sigma^{\sharp}\sharp\upsilon^{\sharp}$}} & {\tiny{}$T_{21}\sigma$} & {\tiny{}$T_{22}\sigma$} & {\tiny{}$(-S_{1}T_{21}-N_{1}T_{11})\mathfrak{1}$} & \multirow{3}{*}{{\tiny{}$*_{1}$}}\tabularnewline
 &  &  &  &  &  & {\tiny{}$+N_{1}\upsilon$} & {\tiny{}$+N_{1}\upsilon^{\sharp}$} & {\tiny{}$+T_{21}\sigma^{\sharp}-N_{1}\sigma\sharp\upsilon$} & \tabularnewline
 &  &  &  &  &  &  &  & {\tiny{}$-S_{1}\sigma^{\sharp}\sharp\upsilon$} & \tabularnewline
\cline{1-1} \cline{4-10} \cline{5-10} \cline{6-10} \cline{7-10} \cline{8-10} \cline{9-10} \cline{10-10} 
\multirow{3}{*}{{\tiny{}$\upsilon$}} &  &  & \multirow{17}{*}{} & \multirow{3}{*}{{\tiny{}$2\upsilon^{\sharp}$}} & {\tiny{}$-N_{2}\cdot\mathfrak{1}$} & {\tiny{}$-T_{12}\cdot\mathfrak{1}-S_{2}\sigma$} & {\tiny{}$N_{2}\sigma$} & {\tiny{}$(S_{1}S_{2}-T_{22})\cdot\mathfrak{1}$} & {\tiny{}$N_{2}\sigma^{\sharp}$}\tabularnewline
 &  &  &  &  & {\tiny{}$-S_{2}\upsilon$} & {\tiny{}$+T_{11}\upsilon-\sigma\sharp\upsilon^{\sharp}$} & {\tiny{}$+T_{11}\upsilon^{\sharp}$} & {\tiny{}$-S_{2}\sigma^{\sharp}+T_{21}\upsilon$} & {\tiny{}$+T_{21}\upsilon^{\sharp}$}\tabularnewline
 &  &  &  &  &  &  &  & {\tiny{}$-S_{1}\upsilon^{\sharp}-\sigma^{\sharp}\sharp\upsilon^{\sharp}$} & \tabularnewline
\cline{1-1} \cline{5-10} \cline{6-10} \cline{7-10} \cline{8-10} \cline{9-10} \cline{10-10} 
\multirow{3}{*}{{\tiny{}$\upsilon^{\sharp}$}} &  &  &  & \multirow{14}{*}{} & \multirow{3}{*}{{\tiny{}$2N_{2}\upsilon$}} & {\tiny{}$N_{2}\sigma$} & {\tiny{}$(-S_{2}T_{12}-N_{2}T_{11})\mathfrak{1}$} & {\tiny{}$N_{2}\sigma^{\sharp}$} & \multirow{2}{*}{{\tiny{}$*_{2}$}}\tabularnewline
 &  &  &  &  &  & {\tiny{}$+T_{12}\upsilon$} & {\tiny{}$+T_{12}\upsilon^{\sharp}-N_{2}\sigma\sharp\upsilon$ } & {\tiny{}$+T_{22}\upsilon$} & \tabularnewline
 &  &  &  &  &  &  & {\tiny{}$-S_{2}\sigma\sharp\upsilon^{\sharp}$} &  & \tabularnewline
\cline{1-1} \cline{6-10} \cline{7-10} \cline{8-10} \cline{9-10} \cline{10-10} 
\multirow{3}{*}{{\tiny{}$\sigma\sharp\upsilon$}} &  &  &  &  & \multirow{11}{*}{} & {\tiny{}$2(T_{12}\sigma+T_{21}\upsilon$} & {\tiny{}$S_{1}N_{2}\mathfrak{1}-T_{11}S_{2}\sigma$} & {\tiny{}$S_{2}N_{1}\mathfrak{1}+T_{22}\sigma$} & \multirow{3}{*}{{\tiny{}$*_{3}$}}\tabularnewline
 &  &  &  &  &  & {\tiny{}$-\sigma^{\sharp}\sharp\upsilon^{\sharp})$} & {\tiny{}$-N_{2}\sigma^{\sharp}+T_{22}\upsilon$} & {\tiny{}$+T_{12}\sigma^{\sharp}-T_{11}S_{1}\upsilon$} & \tabularnewline
 &  &  &  &  &  &  & {\tiny{}$+T_{21}\upsilon^{\sharp}+S_{2}\sigma^{\sharp}\sharp\upsilon$} & {\tiny{}$-N_{1}\upsilon^{\sharp}+S_{1}\sigma\sharp\upsilon^{\sharp}$} & \tabularnewline
\cline{1-1} \cline{7-10} \cline{8-10} \cline{9-10} \cline{10-10} 
\multirow{3}{*}{{\tiny{}$\sigma\sharp\upsilon^{\sharp}$}} &  &  &  &  &  & \multirow{8}{*}{} & {\tiny{}$2(N_{2}T_{11}\sigma+T_{22}\upsilon^{\sharp}$} & \multirow{3}{*}{{\tiny{}$*_{4}$}} & {\tiny{}$T_{21}N_{2}\sigma-N_{1}N_{2}\upsilon$}\tabularnewline
 &  &  &  &  &  &  & {\tiny{}$-N_{2}\sigma^{\sharp}\sharp\upsilon)$} &  & {\tiny{}$-(N_{1}S_{2}+S_{1}T_{12})\upsilon^{\sharp}$}\tabularnewline
 &  &  &  &  &  &  &  &  & {\tiny{}$+T_{11}N_{2}\sigma^{\sharp}+N_{2}S_{1}\sigma\sharp\upsilon$}\tabularnewline
\cline{1-1} \cline{8-10} \cline{9-10} \cline{10-10} 
\multirow{3}{*}{{\tiny{}$\sigma^{\sharp}\sharp\upsilon$}} &  &  &  &  &  &  & \multirow{5}{*}{} & {\tiny{}$2(T_{22}\sigma^{\sharp}+N_{1}T_{11}\upsilon$} & {\tiny{}$-N_{1}N_{2}\sigma+T_{12}N_{1}\upsilon$}\tabularnewline
 &  &  &  &  &  &  &  & {\tiny{}$-N_{1}\sigma\sharp\upsilon^{\sharp})$} & {\tiny{}$-(S_{1}N_{2}+S_{2}T_{21})\sigma^{\sharp}$}\tabularnewline
 &  &  &  &  &  &  &  &  & {\tiny{}$+T_{11}N_{1}\upsilon^{\sharp}+N_{1}S_{2}\sigma\sharp\upsilon$}\tabularnewline
\cline{1-1} \cline{9-10} \cline{10-10} 
\multirow{2}{*}{{\tiny{}$\sigma^{\sharp}\sharp\upsilon^{\sharp}$}} &  &  &  &  &  &  &  & \multirow{2}{*}{} & {\tiny{}$2(N_{2}T_{21}\sigma^{\sharp}+N_{1}T_{12}\upsilon^{\sharp}$}\tabularnewline
 &  &  &  &  &  &  &  &  & {\tiny{}$-N_{1}N_{2}\sigma\sharp\upsilon)$}\tabularnewline
\hline 
\end{tabular}{\tiny\par}
\par\end{flushleft}

\begin{align*}
\ensuremath{*_{1}=\sigma^{\sharp}\sharp(\sigma^{\sharp}\sharp\upsilon^{\sharp})=} & (-N_{1}T_{12}+S_{1}^{2}S_{2}-S_{1}T_{22})\mathfrak{1}-S_{2}N_{1}\sigma\\
- & (S_{1}S_{2}-T_{22})\sigma^{\sharp}-N_{1}\sigma\sharp\upsilon^{\sharp}-S_{1}\sigma^{\sharp}\sharp\upsilon^{\sharp}\\
*_{2}=\upsilon^{\sharp}\sharp(\sigma^{\sharp}\sharp\upsilon^{\sharp})= & (-N_{2}T_{21}+S_{1}S_{2}^{2}-S_{2}T_{22})\mathfrak{1}-S_{1}N_{2}\upsilon\\
- & (S_{1}S_{2}-T_{22})\upsilon^{\sharp}-N_{2}\sigma^{\sharp}\sharp\upsilon-S_{2}\sigma^{\sharp}\sharp\upsilon^{\sharp}\\
*_{3}=(\sigma\sharp\upsilon)\sharp(\sigma^{\sharp}\sharp\upsilon^{\sharp})= & \left(T_{11}(T_{22}-S_{1}S_{2})-T_{21}T_{12}-N_{1}N_{2}\right)\mathfrak{1}\\
- & T_{21}S_{2}\sigma-S_{1}T_{12}\upsilon+(T_{22}-S_{1}S_{2})\sigma\sharp\upsilon\\
- & T_{21}\sigma\sharp\upsilon^{\sharp}-T_{12}\sigma^{\sharp}\sharp\upsilon+T_{11}\sigma^{\sharp}\sharp\upsilon^{\sharp}\\
*_{4}=(\sigma\sharp\upsilon^{\sharp})\sharp(\sigma^{\sharp}\sharp\upsilon)= & \left(T_{11}(S_{1}S_{2}-T_{22})+T_{12}T_{21}-N_{1}N_{2}\right)\mathfrak{1}\\
- & S_{2}T_{11}\sigma^{\sharp}-S_{1}T_{11}u^{\sharp}-T_{22}\sigma\sharp\upsilon\\
+ & T_{21}\sigma\sharp\upsilon^{\sharp}+T_{12}\sigma^{\sharp}\sharp\upsilon-T_{11}\sigma^{\sharp}\sharp\upsilon^{\sharp}
\end{align*}
We may compute the $\sharp$-product on $J$, the $\sharp$-mapping
from $J$, and then the cubic form; we write down the $\sharp$-mapping
in Section \ref{subsec:-mappingPi} using new coordinates of $J$
defined in Section \ref{subsec:Change-of-base} while we omit the
other things since they are lengthy. 

\vspace{1cm}

Hereafter in Section \ref{sec:Coord10paramJ}, \textit{we assume that
the characteristic of $\mathsf{k}$ is not $3$.}

\subsection{Change of base with new parameters\label{subsec:Change-of-base}}

In this subsection, we obtain a coordinatization of $J$ from the
one obtained in the previous subsection changing the base as in there
and introducing new parameters; the number of parameters is 10 here.
We set

\begin{equation}
\sigma_{2}:=\sigma,\,\sigma_{1}:=\sigma_{2}^{\sharp},\,\sigma_{3}:={\rm \mathfrak{1}},\,\upsilon_{1}:=\upsilon.\label{eq:sigma}
\end{equation}
and
\begin{equation}
t_{1}:=S_{1},\,t_{2}:=N_{1}.\label{eq:t1t2}
\end{equation}

\begin{lem}
\label{lem:v2}There exists an element $\upsilon_{2}$ of $\langle\upsilon_{1},\upsilon_{1}^{\sharp},\sigma_{3}\rangle$
such that
\[
\langle\upsilon_{1},\upsilon_{1}^{\sharp},\sigma_{3}\rangle=\langle\upsilon_{1},\upsilon_{2},\sigma_{3}\rangle\ \text{and}\ T(\upsilon_{2})=0.
\]
\end{lem}

\begin{proof}
We have only to find $\upsilon_{2}=\alpha\upsilon_{1}+\beta\upsilon_{1}^{\sharp}+\gamma\sigma_{3}$
such that $\beta\not=0$ and $T(\upsilon_{2})=\beta S_{2}+3\gamma=0$,
which clearly exists.
\end{proof}
Note that 
\begin{equation}
\upsilon_{1}\sharp\upsilon_{1}^{\sharp}=-N(\upsilon_{1})\mathfrak{1}-S(\upsilon_{1})\upsilon_{1}\label{eq:v1v1sha}
\end{equation}
by \cite[(14)]{R} and $T(\upsilon_{1})=0$. Since $\upsilon_{2}$
is a linear combination of $\upsilon_{1},\upsilon_{1}^{\sharp},\sigma_{3}$,
so is $\upsilon_{2}^{\sharp}$ by (\ref{eq:v1v1sha}) and \cite[(5)]{R}.
We set
\begin{equation}
\begin{cases}
t_{123} & :=-T(\sigma_{1},\upsilon_{1})=-T_{21}\\
t_{124} & :=-T(\sigma_{1},\upsilon_{2}),\\
t_{125} & :=T(\sigma_{2},\upsilon_{1})=T_{11},\\
t_{126} & :=T(\sigma_{2},\upsilon_{2}).
\end{cases}\label{eq:t123t126}
\end{equation}
 We may express as 
\begin{align}
 & \upsilon_{1}^{\sharp}=t{}_{136}\upsilon_{1}-t_{135}\upsilon_{2}+\nicefrac{1}{3}S(\upsilon_{1})\mathfrak{1},\label{eq:v1}\\
 & \upsilon_{2}^{\sharp}=t_{246}\upsilon_{1}-t{}_{245}\upsilon_{2}+\nicefrac{1}{3}S(\upsilon_{2})\mathfrak{1}\label{eq:v2}
\end{align}
with some $t_{135},t_{136},\,t_{245},\,t_{246}\in\mathsf{k}$, where
the coefficient of $\mathfrak{1}$ in $\upsilon_{1}^{\sharp}$ is
$\nicefrac{1}{3}S(\upsilon_{1})$ since $T(\upsilon_{1})=T(\upsilon_{2})=0$
by the choice of $\upsilon_{1}$ and $\upsilon_{2}$, and $T(\upsilon_{1}^{\sharp})=S(\upsilon_{1})$
and $T(\mathfrak{1})=3$ by \cite{R}, and similarly for the coefficient
of $\mathfrak{1}$ in $\upsilon_{2}^{\sharp}$ (\textit{here we need
the assumption that the characteristic of $\mathsf{k}$ is $3$}).
Note that $t_{135}\not=0$ by Claim \ref{lem:v2}.
\begin{prop}
\label{prop:S2N2}The following equalities hold:
\begin{align}
S_{2}= & S(\upsilon_{1})=3(t_{135}t_{245}-t_{136}^{2}),\label{eq:Sv1}\\
 & S(\upsilon_{2})=3(t_{246}t_{136}-t_{245}^{2}),\label{eq:Sv2}\\
N_{2}= & N(\upsilon_{1})=2t_{136}^{3}-3t_{135}t_{136}t_{245}+t_{135}^{2}t_{246},\label{eq:N1}\\
\upsilon_{1}\sharp\upsilon_{2}= & 2t_{245}\upsilon_{1}-2t_{136}\upsilon_{2}+(t_{135}t_{246}-t_{136}t_{245})\mathfrak{1}.\label{eq:v1v2}
\end{align}
\end{prop}

\begin{proof}
By (\ref{eq:v1}), we obtain 
\begin{align}
-N(\upsilon_{1})\mathfrak{1}-S(\upsilon_{1})\upsilon_{1} & =\upsilon_{1}\sharp\upsilon_{1}^{\sharp}\label{eq:v1v1}\\
 & =\upsilon_{1}\sharp(t{}_{136}\upsilon_{1}-t_{135}\upsilon_{2}+\nicefrac{1}{3}\,S(\upsilon_{1})\mathfrak{1})\nonumber \\
 & =2t_{136}\upsilon_{1}^{\sharp}-t_{135}\upsilon_{1}\sharp\upsilon_{2}-\nicefrac{1}{3}\,S(\upsilon_{1})\upsilon_{1}\nonumber \\
 & =2t_{136}(t{}_{136}\upsilon_{1}-t_{135}\upsilon_{2}+\nicefrac{1}{3}\,S(\upsilon_{1})\mathfrak{1})-t_{135}\upsilon_{1}\sharp\upsilon_{2}-\nicefrac{1}{3}\,S(\upsilon_{1})\upsilon_{1}.\nonumber 
\end{align}
By (\ref{eq:v1}) and (\ref{eq:v2}), we have 
\begin{align}
N(\upsilon_{1})\upsilon_{1} & =(\upsilon_{1}^{\sharp})^{\sharp}=(t{}_{136}\upsilon_{1}-t_{135}\upsilon_{2}+\nicefrac{1}{3}\,S(\upsilon_{1})\mathfrak{1})^{\sharp}\label{eq:Nv1}\\
 & =(t_{136}\upsilon_{1}-t_{135}\upsilon_{2})\sharp\left(\nicefrac{1}{3}\,S(\upsilon_{1})\mathfrak{1}\right)+(t{}_{136}\upsilon_{1}-t_{135}\upsilon_{2})^{\sharp}+\left(\nicefrac{1}{3}\,S(\upsilon_{1})\mathfrak{1}\right){}^{\sharp}\nonumber \\
 & =(t_{136}\upsilon_{1}-t_{135}\upsilon_{2})\sharp\left(\nicefrac{1}{3}\,S(\upsilon_{1})\mathfrak{1}\right)-t_{136}t_{135}\upsilon_{1}\sharp\upsilon_{2}+t_{136}^{2}\upsilon_{1}^{\sharp}+t_{135}^{2}\upsilon_{2}{}^{\sharp}+(\nicefrac{1}{3}\,S(\upsilon_{1})\mathfrak{1})^{\sharp}\nonumber \\
 & =-\nicefrac{1}{3}S(\upsilon_{1})(t_{136}\upsilon_{1}-t_{135}\upsilon_{2})-t_{136}t_{135}\upsilon_{1}\sharp\upsilon_{2}+t_{136}^{2}\left(t_{136}\upsilon_{1}-t_{135}\upsilon_{2}+\nicefrac{1}{3}\,S(\upsilon_{1})\mathfrak{1}\right)\nonumber \\
 & +t_{135}^{2}\left(t_{246}\upsilon_{1}-t{}_{245}\upsilon_{2}+\nicefrac{1}{3}S(\upsilon_{2})\mathfrak{1}\right)+\nicefrac{1}{9}\,S(\upsilon_{1})^{2}\mathfrak{1}.\nonumber 
\end{align}
Eliminating $t_{135}\upsilon_{1}\sharp\upsilon_{2}$ from (\ref{eq:Nv1})
by (\ref{eq:v1v1}), we obtain an equality. As we mention above, it
holds that $t_{135}\not=0$. Then, considering the coefficient of
$\upsilon_{2}$ in the equality, we have (\ref{eq:Sv1}). Inserting
this to the equality and considering the coefficient of $\upsilon_{1},$we
have (\ref{eq:N1}). Inserting this and (\ref{eq:Sv1}) to the equality
and considering the coefficient of $\mathfrak{1}$, we have (\ref{eq:Sv2}).
Finally, inserting (\ref{eq:Sv1}) and (\ref{eq:N1}) to (\ref{eq:v1v1}),
we have (\ref{eq:v1v2}). 
\end{proof}
We set 
\begin{equation}
\begin{cases}
\pi_{1} & :=\sigma_{2}\sharp\upsilon_{1},\\
\pi_{2} & :=(\sigma_{2}^{\sharp})\sharp\upsilon_{1}+t_{1}\upsilon_{1},\\
\pi_{3} & :=\sigma_{2}\sharp\upsilon_{2},\\
\pi_{4} & :=(\sigma_{2}^{\sharp})\sharp\upsilon_{2}+t_{1}\upsilon_{2}.
\end{cases}\label{eq:pi1pi4}
\end{equation}
Then 
\begin{equation}
\pi_{1},\pi_{2},\pi_{3},\pi_{4},\upsilon_{1},\upsilon_{2},\sigma_{1},\sigma_{2},\sigma_{3}\label{eq:new basis}
\end{equation}
form another basis of $V$. 

The elements of the new basis (\ref{eq:new basis}) of $J$ are written
by those of the original basis (\ref{eq:genV}) of $J$ as follows: 
\begin{lem}
\label{lem:basechange} The following base change rule holds:
\begin{align*}
 & \sigma_{1}=\sigma^{\sharp},\,\sigma_{2}=\sigma,\,\sigma_{3}=\mathfrak{1},\upsilon_{1}=\upsilon,\\
 & \upsilon_{2}=\left(t_{245}-\frac{t_{136}^{2}}{t_{135}}\right)\mathfrak{1}+\frac{t_{136}}{t_{135}}\upsilon-\frac{1}{t_{135}}\upsilon^{\sharp},\\
 & \pi_{1}=\sigma\sharp\upsilon,\,\pi_{2}=\sigma^{\sharp}\sharp\upsilon+t_{1}\upsilon,\\
 & \pi_{3}=-\left(t_{245}-\frac{t_{136}^{2}}{t_{135}}\right)\sigma+\frac{t_{136}}{t_{135}}\sigma\sharp\upsilon-\frac{1}{t_{135}}\sigma\sharp\upsilon^{\sharp},\\
 & \pi_{4}=\frac{2t_{1}(t_{135}t_{245}-t_{136}^{2})}{t_{135}}\mathfrak{1}-\frac{t_{135}t_{245}-t_{136}^{2}}{t_{135}}\sigma^{\sharp}+\frac{t_{1}t_{136}}{t_{135}}\upsilon-\frac{3t_{1}}{t_{135}}\upsilon^{\sharp}\\
 & +\frac{t_{136}}{t_{135}}\sigma^{\sharp}\sharp\upsilon-\frac{1}{t_{135}}\sigma^{\sharp}\sharp\upsilon^{\sharp}.
\end{align*}
\end{lem}

\begin{proof}
The 1st and 3rd lines follows just from the definitions (\ref{eq:sigma})
and (\ref{eq:pi1pi4}). The 2nd line follows from (\ref{eq:v1}) and
the fact that $t_{135}\not=0$, which is mentioned above Proposition
\ref{prop:S2N2}. The 4th and 5th lines follow from the 2nd line and
(\ref{eq:pi1pi4}). 
\end{proof}
\begin{lem}
\label{lem:parameterchange} The following transformation rule of
the parameters holds: 
\begin{align*}
 & S_{1}=t_{1},\\
 & S_{2}=3(t_{135}t_{245}-t_{136}^{2}),\,N_{1}=t_{2},\\
 & N_{2}=2t_{136}^{3}-3t_{135}t_{136}t_{245}+t_{135}^{2}t_{246},\\
 & T_{11}=t_{125},\,T_{21}=-t_{123},\\
 & T_{12}=t_{136}t_{125}-t_{135}t_{126},\\
 & T_{22}=(t_{135}t_{124}-t_{136}t_{123})+t_{1}(t_{135}t_{245}-t_{136}^{2}).
\end{align*}
\end{lem}

\begin{proof}
The 1st and 4th lines follow just from the definitions (\ref{eq:t1t2})
and (\ref{eq:t123t126}). The 2nd and 3rd lines are obtained in Proposition
\ref{prop:S2N2}. Since $T(\sigma_{2},\upsilon_{1}^{\sharp})=T(\sigma_{2},t{}_{136}\upsilon_{1}-t_{135}\upsilon_{2}+\nicefrac{1}{3}S(\upsilon_{1}))$
by (\ref{eq:Sv1}), the 5th line follows from the assumption that
$T(\sigma_{2})=0$ and (\ref{eq:t123t126}). Since $T(\sigma_{2}^{\sharp},\upsilon_{1}^{\sharp})=T(\sigma_{2}^{\sharp},t_{136}\upsilon_{1}-t_{135}\upsilon_{2}+\nicefrac{1}{3}S(\upsilon_{1})\mathfrak{1})$
by (\ref{eq:Sv1}), the 6th line follows from (\ref{eq:sigma}), (\ref{eq:t123t126})
and (\ref{eq:t1t2}).
\end{proof}
By Lemmas \ref{lem:basechange} and \ref{lem:parameterchange}, we
obtain the following from a direct computation:
\begin{lem}
\label{lem:two coord}Let $x$ be an element of $J$ and write 
\begin{align*}
x & =x_{1}\mathfrak{1}+x_{2}\sigma+x_{3}\sigma^{\sharp}+x_{4}\upsilon+x_{5}\upsilon^{\sharp}+x_{6}\sigma\sharp\upsilon+x_{7}\sigma\sharp\upsilon^{\sharp}+x_{8}\sigma^{\sharp}\sharp\upsilon+x_{9}\sigma^{\sharp}\sharp\upsilon^{\sharp}\\
 & =p_{1}\pi_{1}+p_{2}\pi_{2}+p_{3}\pi_{3}+p_{4}\pi_{4}+u_{1}\upsilon_{1}+u_{2}\upsilon_{2}+s_{1}\sigma_{1}+s_{2}\sigma_{2}+s_{3}\sigma_{3}
\end{align*}
in two ways by using the old and new bases. The following transformation
rule of the coordinates holds:
\begin{align*}
x_{1} & =(s_{3}+2p_{4}t_{1}t_{245}+t_{245}u_{2})-(2p_{4}t_{1}+u_{2})\frac{t_{136}^{2}}{t_{135}}\\
x_{2} & =(s_{2}-p_{3}t_{245})+\frac{p_{3}t_{136}^{2}}{t_{135}},\\
x_{3} & =(s_{1}-p_{4}t_{245})+\frac{p_{4}t_{136}^{2}}{t_{135}},\\
x_{4} & =(u_{1}+p_{2}t_{1})+(p_{4}t_{1}+u_{2})\frac{t_{136}}{t_{135}},\\
x_{5} & =-\frac{p_{4}t_{1}+u_{2}}{t_{135}},\,x_{6}=p_{1}+\frac{p_{3}t_{136}}{t_{135}},\\
x_{7} & =-\frac{p_{3}}{t_{135}},\,x_{8}=p_{2}+\frac{p_{4}t_{136}}{t_{135}},\,x_{9}=-\frac{p_{4}}{t_{135}}.
\end{align*}
\end{lem}

\subsection{$\sharp$-mapping\label{subsec:-mappingPi}}

As we have said in the end of Subsection \ref{subsec:Coordination-of-Perterson},
we write down the $\sharp$-mapping of $J$ using the new coordinates
of $J$ defined in Subsection \ref{subsec:Change-of-base}.
\begin{prop}
\label{prop:The--mapping-Pi}The $\sharp$-mapping $J\to J$ is defined
by 
\begin{align*}
(p_{1},p_{2},p_{3},p_{4},u_{1},u_{2},s_{1},s_{2},s_{3}) & \mapsto\\
 & ({\rm ad}p_{1},{\rm ad}p_{2},{\rm ad}p_{3},{\rm ad}p_{4},{\rm ad}u_{1},{\rm ad}u_{2},{\rm ad}s_{1},{\rm ad}s_{2},{\rm ad}s_{3})
\end{align*}
satisfying the following:

\begin{align*}
F_{1}:=3{\rm ad}p_{2} & =u_{1}s_{1}-p_{1}s_{2}-p_{2}s_{3}+t_{126}(p_{2}p_{3}-p_{1}p_{4})+t_{136}(-p_{1}^{2}-t_{1}p_{2}^{2}-p_{2}u_{1})\\
 & +t_{245}(-2p_{1}p_{3}-2t_{1}p_{2}p_{4}-p_{4}u_{1}-p_{2}u_{2})+t_{246}(-p_{3}^{2}-t_{1}p_{4}^{2}-p_{4}u_{2}),
\end{align*}

\begin{align*}
F_{2}:=3{\rm ad}p_{4} & =u_{2}s_{1}-p_{3}s_{2}-p_{4}s_{3}+t_{125}(-p_{2}p_{3}+p_{1}p_{4})+t_{135}(p_{1}^{2}+t_{1}p_{2}^{2}+p_{2}u_{1})\\
 & +t_{136}(2p_{1}p_{3}+2t_{1}p_{2}p_{4}+p_{4}u_{1}+p_{2}u_{2})+t_{245}(p_{3}^{2}+t_{1}p_{4}^{2}+p_{4}u_{2}),
\end{align*}
\begin{align*}
F_{3}:=3{\rm ad}p_{1} & =-t_{2}p_{2}s_{1}+(t_{1}p_{2}+u_{1})s_{2}-p_{1}s_{3}+t_{124}(p_{2}p_{3}-p_{1}p_{4})\\
 & +t_{136}(-p_{2}^{2}t_{2}-p_{1}u_{1})+t_{245}(-2t_{2}p_{2}p_{4}-p_{3}u_{1}-p_{1}u_{2})+t_{246}(-t_{2}p_{4}^{2}-p_{3}u_{2}),
\end{align*}
\begin{align*}
F_{4}:=3{\rm ad}p_{3} & =-t_{2}p_{4}s_{1}+(t_{1}p_{4}+u_{2})s_{2}-p_{3}s_{3}+t_{123}(-p_{2}p_{3}+p_{1}p_{4})+t_{135}(t_{2}p_{2}^{2}+p_{1}u_{1})\\
 & +t_{136}(2t_{2}p_{2}p_{4}+p_{3}u_{1}+p_{1}u_{2})+t_{245}(t_{2}p_{4}^{2}+p_{3}u_{2}),
\end{align*}
\begin{align*}
F_{5}:=-3({\rm ad}u_{1}+t_{1}{\rm ad}p_{2}) & =-t_{2}p_{1}s_{1}+(t_{1}p_{1}-t_{2}p_{2})s_{2}+(t_{1}p_{2}+u_{1})s_{3}\\
 & +t_{123}(t_{1}p_{2}^{2}+p_{1}^{2}+p_{2}u_{1})+t_{124}(t_{1}p_{2}p_{4}+p_{1}p_{3}+p_{2}u_{2})\\
 & -t_{125}(t_{2}p_{2}^{2}+p_{1}u_{1})-t_{126}(t_{2}p_{2}p_{4}+p_{1}u_{2})-t_{136}(t_{1}p_{2}u_{1}-t_{2}p_{1}p_{2}+u_{1}^{2})\\
 & -t_{245}(t_{1}p_{4}u_{1}+t_{1}p_{2}u_{2}-t_{2}p_{2}p_{3}-t_{2}p_{1}p_{4}+2u_{1}u_{2})\\
 & -t_{246}(t_{1}p_{4}u_{2}-t_{2}p_{3}p_{4}+u_{2}^{2}),
\end{align*}
\begin{align*}
F_{6}:=-3({\rm ad}u_{2}+t_{1}{\rm ad}p_{4}) & =-t_{2}p_{3}s_{1}+(t_{1}p_{3}-t_{2}p_{4})s_{2}+(t_{1}p_{4}+u_{2})s_{3}\\
 & +t_{123}(p_{1}p_{3}+t_{1}p_{2}p_{4}+p_{4}u_{1})+t_{124}(p_{3}^{2}+t_{1}p_{4}^{2}+p_{4}u_{2})\\
 & -t_{125}(t_{2}p_{2}p_{4}+p_{3}u_{1})-t_{126}(t_{2}p_{4}^{2}+p_{3}u_{2})+t_{135}(t_{1}p_{2}u_{1}-t_{2}p_{1}p_{2}+u_{1}^{2})\\
 & +t_{136}(t_{1}p_{4}u_{1}+t_{1}p_{2}u_{2}-t_{2}p_{2}p_{3}-t_{2}p_{1}p_{4}+2u_{1}u_{2})\\
 & +t_{245}(t_{1}p_{4}u_{2}-t_{2}p_{3}p_{4}+u_{2}^{2}),
\end{align*}
\begin{align*}
F_{7}:=-3{\rm ad}s_{1} & =-s_{2}^{2}+s_{1}s_{3}+(t_{123}p_{2}+t_{124}p_{4})s_{1}-(t_{125}p_{2}+t_{126}p_{4})s_{2}\\
 & +(t_{136}^{2}-t_{135}t_{245})(p_{1}^{2}-2p_{2}u_{1}-t_{1}p_{2}^{2})+(t_{245}^{2}-t_{136}t_{246})(p_{3}^{2}-2p_{4}u_{2}-t_{1}p_{4}^{2})\\
 & +(t_{135}t_{246}-t_{136}t_{245})(p_{2}u_{2}+p_{4}u_{1}-p_{1}p_{3}+t_{1}p_{2}p_{4})\\
 & +(t_{123}t_{136}-t_{124}t_{135})p_{2}^{2}+2(t_{123}t_{245}-t_{124}t_{136})p_{2}p_{4}+(t_{123}t_{246}-t_{124}t_{245})p_{4}^{2}\\
 & +(t_{126}t_{135}-t_{125}t_{136})p_{1}p_{2}+(t_{126}t_{136}-t_{125}t_{245})(p_{1}p_{4}+p_{2}p_{3})\\
 & +(t_{126}t_{245}-t_{125}t_{246})p_{3}p_{4},
\end{align*}
\begin{align*}
F_{8}:=-3{\rm ad}s_{2} & =(s_{2}s_{3}+t_{1}s_{1}s_{2}-t_{2}s_{1}^{2})+(t_{123}p_{1}+t_{124}p_{3})s_{1}-(t_{125}p_{1}+t_{126}p_{3})s_{2}\\
 & +(t_{136}^{2}-t_{135}t_{245})(-2p_{1}u_{1}-2t_{1}p_{1}p_{2}+t_{2}p_{2}^{2})+(t_{245}^{2}-t_{136}t_{246})(-2p_{3}u_{2}-2t_{1}p_{3}p_{4}+t_{2}p_{4}^{2})\\
 & +(t_{135}t_{246}-t_{136}t_{245})\left(p_{1}u_{2}+p_{3}u_{1}+t_{1}(p_{1}p_{4}+p_{2}p_{3})-t_{2}p_{2}p_{4}\right)\\
 & +(t_{126}t_{135}-t_{125}t_{136})p_{1}^{2}+2(t_{126}t_{136}-t_{125}t_{245})p_{1}p_{3}+(t_{126}t_{245}-t_{125}t_{246})p_{3}^{2}\\
 & +(t_{123}t_{136}-t_{124}t_{135})p_{1}p_{2}+(t_{123}t_{245}-t_{124}t_{136})(p_{1}p_{4}+p_{2}p_{3})\\
 & +(t_{123}t_{246}-t_{124}t_{245})p_{3}p_{4},
\end{align*}
\begin{align*}
F_{9} & :=3({\rm ad}s_{3}+t_{1}{\rm ad}s_{1}-t_{125}{\rm ad}p_{1}+t_{123}{\rm ad}p_{2}-t_{126}{\rm ad}p_{3}+t_{124}{\rm ad}p_{4})\\
 & =(s_{3}^{2}+t_{1}s_{2}^{2}-t_{2}s_{1}s_{2})+(t_{123}u_{1}+t_{124}u_{2})s_{1}-(t_{125}u_{1}+t_{126}u_{2})s_{2}\\
 & +(-t_{136}^{2}+t_{135}t_{245})(2t_{2}p_{1}p_{2}+u_{1}^{2})+(-t_{136}t_{245}+t_{135}t_{246})\left(t_{2}(p_{2}p_{3}+p_{1}p_{4})+u_{1}u_{2}\right)\\
 & +(-t_{245}^{2}+t_{136}t_{246})(2t_{2}p_{3}p_{4}+u_{2}^{2})\\
 & -t_{1}\left((t_{123}t_{136}-t_{124}t_{135})p_{2}^{2}+2(t_{123}t_{245}-t_{124}t_{136})p_{2}p_{4}+(t_{123}t_{246}-t_{124}t_{245})p_{4}^{2}\right)\\
 & +t_{2}\left((t_{125}t_{136}-t_{126}t_{135})p_{2}^{2}+2(t_{125}t_{245}-t_{126}t_{136})p_{2}p_{4}+(t_{125}t_{246}-t_{126}t_{245})p_{4}^{2}\right)\\
 & +(t_{124}t_{135}-t_{123}t_{136})p_{1}^{2}-2(-t_{124}t_{136}+t_{123}t_{245})p_{1}p_{3}+(t_{124}t_{245}-t_{123}t_{246})p_{3}^{2}\\
 & +(t_{124}t_{125}-t_{123}t_{126})(-p_{2}p_{3}+p_{1}p_{4}).
\end{align*}
\end{prop}

\begin{rem}
\begin{enumerate}[(1)]

\item By Proposition \ref{prop:The--mapping-Pi}, the coefficient
$\nicefrac{1}{3}$ is inevitable for the $\sharp$-mapping.
\item It is possible to write down the cubic form but we omit it
since it is lengthy.

\end{enumerate}
\end{rem}

\section{\textbf{Affine variety $\Pi_{\mA}^{15}$}\label{sec:Affine-varietyPi}\protect 
}In this section, we assume for simplicity that $\mathsf{k}=\mC$.

\subsection{Definition of $\Pi_{\mA}^{15}${\small{}\label{subsec:Definition-ofPi}}}
\begin{defn}
\label{def:Pi15} Let $\mA_{\Pi}$ be the affine $19$-space with
coordinates 
\begin{equation}
p_{1},\dots,p_{4},u_{1},u_{2},s_{1},s_{2},s_{3},t_{1},t_{2},t_{123},t_{124},t_{125},t_{126},t_{135},t_{136},t_{245},t_{246}.\label{eq:coordpi}
\end{equation}
 In $\mA_{\Pi}$, we define $\Pi_{\mA}^{15}$ to be the affine scheme
defined by 
\[
F_{1}=\cdots=F_{9}=0.
\]
\end{defn}

Finally we will show that $\Pi_{\mA}^{15}$ is a variety, i.e., it
is irreducible. Surprisingly, again, the affine variety $\Pi_{\mA}^{15}$
coincides with the one constructed in \cite{Tak3} inspired by \cite[Subsec. 5.2]{Tay}
Our construction is based on the theory of unprojection. We investigate
several properties of $\Pi_{\mA}^{15}$ in the following subsections.

\subsection{Weights for variables and equations\label{subsec:Weights-for-variablesPi}}

We assign weights for variables of the polynomial ring $S_{\Pi}$
such that all the equations $F_{1}$--$F_{9}$ are homogeneous. Moreover,
we assume that all the variables are not zero allowing some of them
are constants. Then it is easy to derive the following relations between
the weights of variables of $S_{\Pi}$:

\begin{align*}
 & {\rm \mathbf{Weights\,for\,variables}}\\
 & w(p_{1})=w(p_{3})+w(u_{1})-w(u_{2}),\,w(p_{2})=2w(p_{3})+w(u_{1})-2w(u_{2}),\,w(p_{4})=2w(p_{3})-w(u_{2}),\\
 & w(s_{1})=w(G)-(w(u_{1})+3w(u_{2})-2w(p_{3})),\,w(s_{2})=w(G)-(w(u_{1})+2w(u_{2})-w(p_{3})),\\
 & w(s_{3})=w(G)-(w(u_{1})+w(u_{2})),\\
 & w(t_{1})=2(w(u_{2})-w(p_{3})),\,w(t_{2})=3(w(u_{2})-w(p_{3})),\\
 & w(t_{123})=w(G)-(2w(p_{3})+2w(u_{1})-w(u_{2})),\,w(t_{124})=w(G)-(2w(p_{3})+w(u_{1})),\\
 & w(t_{125})=w(G)-(w(p_{3})+2w(u_{1})),\,w(t_{126})=w(G)-(w(p_{3})+w(u_{1})+w(u_{2})),\\
 & w(t_{135})=w(G)-3w(u_{1}),\,w(t_{136})=w(G)-(2w(u_{1})+w(u_{2})),\\
 & w(t_{245})=w(G)-(w(u_{1})+2w(u_{2})),\,w(t_{246})=w(G)-3w(u_{2}),\\
 & {\rm \mathbf{Weights\,for\,equations}}\\
 & d_{1}:=w(F_{1})=w(G)-(3w(u_{2})-2w(p_{3})),\,d_{2}:=w(F_{2})=w(G)-(w(u_{1})+2w(u_{2})-2w(p_{3})),\\
 & d_{3}:=w(F_{3})=w(G)-(2w(u_{2})-w(p_{3})),\,d_{4}:=w(F_{4})=w(G)-(w(u_{1})+w(u_{2})-w(p_{3})),\\
 & d_{5}:=w(F_{5})=w(G)-w(u_{2}),\,d_{6}:=w(F_{6})=w(G)-w(u_{1}),\\
 & d_{7}:=w(F_{7})=2(w(G)-w(u_{1})-2w(u_{2})+w(p_{3})),\\
 & d_{8}:=w(F_{8})=2w(G)-(2w(u_{1})+3w(u_{2})-w(p_{3})),\\
 & d_{9}:=w(F_{9})=2\big(w(G)-(w(u_{1})+w(u_{2}))\big).
\end{align*}

We need this result in the proofs of Lemma \ref{lem:primePi} and
Proposition \ref{prop:UFDPi}.

\subsection{${\rm GL_{2}}$-action on $\Pi_{\mathbb{A}}^{15}$\label{sec:--and--actions}}
\begin{prop}
\label{prop:GL}For any element $g=\left(\begin{array}{cc}
a & b\\
c & d
\end{array}\right)$ of ${\rm GL}_{2}$, we define $\widehat{g}$ as in Proposition \ref{prop:For-any-elementGL2GL2Up},
{\small{}and denote by }$g^{\dagger}$ the adjoint matrix of $g$.
The group ${\rm GL_{2}}$ acts on the affine scheme $\Pi_{\mathbb{A}}^{15}$
by the following rule: for any $g\in{\rm GL_{2}}$, we set

{\small{}
\begin{equation}
\begin{cases}
\left(\begin{array}{cc}
p_{1} & p_{3}\\
p_{2} & p_{4}\\
u_{1} & u_{2}
\end{array}\right)\mapsto\left(\begin{array}{cc}
p_{1} & p_{3}\\
p_{2} & p_{4}\\
u_{1} & u_{2}
\end{array}\right)g,\left(\begin{array}{c}
s_{1}\\
s_{2}\\
s_{3}
\end{array}\right)\mapsto\det g\left(\begin{array}{c}
s_{1}\\
s_{2}\\
s_{3}
\end{array}\right),\\
\left(\begin{array}{cc}
t_{123} & t_{125}\\
t_{124} & t_{126}
\end{array}\right)\mapsto g^{\dagger}\left(\begin{array}{cc}
t_{123} & t_{125}\\
t_{124} & t_{126}
\end{array}\right),\\
\left(\begin{array}{ccc}
t_{136} & t_{245} & t_{246}\\
-t_{135} & -t_{136} & -t_{245}
\end{array}\right)\mapsto\empty^{t}\!\!g\left(\begin{array}{ccc}
t_{136} & t_{245} & t_{246}\\
-t_{135} & -t_{136} & -t_{245}
\end{array}\right)\hat{g}, & t_{1}\mapsto t_{1},\,t_{2}\mapsto t_{2}.
\end{cases}\label{eq:GL2act}
\end{equation}
}{\small\par}
\end{prop}

\begin{proof}
By direct computations, we may check that {\small{}
\begin{align*}
 & \left(\begin{array}{ccc}
F_{1} & F_{3} & F_{5}\\
F_{2} & F_{4} & F_{6}
\end{array}\right)\mapsto(\det g)\,\empty^{t}\!\!g\left(\begin{array}{ccc}
F_{1} & F_{3} & F_{5}\\
F_{2} & F_{4} & F_{6}
\end{array}\right),\\
 & \left(\begin{array}{c}
F_{7}\\
F_{8}\\
F_{9}
\end{array}\right)\mapsto(\det g)^{2}\left(\begin{array}{c}
F_{7}\\
F_{8}\\
F_{9}
\end{array}\right)
\end{align*}
}for $g\in{\rm GL_{2}}${\small{}.} Therefore the above rule defines
a ${\rm GL_{2}}$-action on $\Pi_{\mathbb{A}}^{15}$.
\end{proof}
\begin{rem}
The $GL_{2}$-action on $\Pi_{\mathbb{A}}^{15}$ is more visible by
the following presentations of $F_{1}$--$F_{6}$. 

We set{\small{}
\begin{equation}
\begin{cases}
\mathsf{L} & :=\left(\begin{array}{ccc}
t_{136} & t_{245} & t_{246}\\
-t_{135} & -t_{136} & -t_{245}
\end{array}\right),\\
\mathsf{U} & :=\left(\begin{array}{c}
U_{1}\\
U_{2}\\
U_{3}
\end{array}\right)=\left(\begin{array}{c}
p_{1}^{2}+t_{1}p_{2}^{2}+p_{2}u_{1}\\
2(p_{1}p_{3}+t_{1}p_{2}p_{4})+(p_{4}u_{1}+p_{2}u_{2})\\
p_{3}^{2}+t_{1}p_{4}^{2}+p_{4}u_{2}
\end{array}\right),\\
\mathsf{I} & :=\left(\begin{array}{c}
I_{1}\\
I_{2}\\
I_{3}
\end{array}\right)=\left(\begin{array}{c}
t_{2}p_{2}^{2}+p_{1}u_{1}\\
2t_{2}p_{2}p_{4}+(p_{3}u_{1}+p_{1}u_{2})\\
t_{2}p_{4}^{2}+p_{3}u_{2}
\end{array}\right),\\
\mathsf{A} & :=\left(\begin{array}{c}
A_{1}\\
A_{2}\\
A_{3}
\end{array}\right)=\left(\begin{array}{c}
t_{1}p_{2}u_{1}-t_{2}p_{1}p_{2}+u_{1}^{2}\\
t_{1}(p_{4}u_{1}+p_{2}u_{2})-t_{2}(p_{2}p_{3}+p_{1}p_{4})\\
t_{1}p_{4}u_{1}-t_{2}p_{3}p_{4}+u_{2}^{2}
\end{array}+2u_{1}u_{2}\right),
\end{cases}\label{eq:LUIA}
\end{equation}
}{\small\par}

{\small{}
\begin{equation}
\begin{cases}
\bm{{v}}_{1} & :=(-p_{2}p_{3}+p_{1}p_{4})\left(\begin{array}{c}
t_{126}\\
-t_{125}
\end{array}\right)+\mathsf{LU},\\
\bm{{v}}_{2} & :=(-p_{2}p_{3}+p_{1}p_{4})\left(\begin{array}{c}
t_{124}\\
-t_{123}
\end{array}\right)+\mathsf{LI},\\
\bm{{v}}_{3} & :=-\left(\begin{array}{cc}
U_{2} & -U_{1}\\
U_{3} & -U_{2}
\end{array}\right)\left(\begin{array}{c}
t_{124}\\
-t_{123}
\end{array}\right)+\left(\begin{array}{cc}
I_{2} & -I_{1}\\
I_{3} & -I_{2}
\end{array}\right)\left(\begin{array}{c}
t_{126}\\
-t_{125}
\end{array}\right)+\mathsf{LA}\\
 & +1/2(-p_{3}u_{1}+p_{1}u_{2})\left(\begin{array}{c}
t_{126}\\
-t_{125}
\end{array}\right)+1/2(p_{4}u_{1}-p_{2}u_{2})\left(\begin{array}{c}
t_{124}\\
-t_{123}
\end{array}\right).
\end{cases}\label{eq:v1v2v3}
\end{equation}
}{\small\par}

Then we have the following simple presentations of $F_{1}$--$F_{6}$:

{\small{}
\begin{align*}
 & \left(\begin{array}{ccc}
F_{1} & F_{3} & F_{5}\\
F_{2} & F_{4} & F_{6}
\end{array}\right)=\\
 & \quad\left(\begin{array}{ccc}
u_{1} & -p_{1} & -p_{2}\\
u_{2} & -p_{3} & -p_{4}
\end{array}\right)\left(\begin{array}{ccc}
s_{1} & s_{2} & s_{3}\\
s_{2} & s_{3} & t_{2}s_{1}-t_{1}s_{2}\\
s_{3} & t_{2}s_{1}-t_{1}s_{2} & t_{2}s_{2}-t_{1}s_{3}
\end{array}\right)-(\bm{{v}}_{1}\bm{\,{v}}_{2}\,\bm{{v}}_{3}).
\end{align*}
} It is desirable that $F_{7}$--$F_{9}$ are also interpreted in
this context. 
\end{rem}

\subsection{Singular locus of $\Pi_{\mathbb{A}}^{15}$\label{sec:Singular-locus-Pi}}


Let $\mathbb{A}_{T}$ be the affine $10$-space with the coordinates
\[
t_{1},\,t_{2},\,t_{123},\,t_{124},\,t_{125},\,t_{126},\,t_{135},\,t_{136},\,t_{245},\,t_{246}.
\]
 In this section, we consider $\mathbb{A}_{T}$ is contained in the
ambient affine space $\mA_{\Pi}$ of $\Pi_{\mathbb{A}}^{15}$. Then
we see that $\mathbb{A}_{T}\subset\Pi_{\mathbb{A}}^{15}$. 
\begin{prop}
\label{prop:SingPi} The open subset $\Pi_{\mathbb{A}}^{15}\setminus\mA_{T}$
of $\Pi_{\mathbb{A}}^{15}$ is $15$-dimensional and irreducible.
Its singular locus is an $8$-dimensional locally closed subset $\Delta$
such that

{\small{}
\begin{align*}
 & \Delta\cap\{p_{2}p_{4}\not=0\}=\\
 & \{\,p_{2}p_{4}\not=0,\\
 & p_{1}=(p_{2}p_{3})/p_{4},\\
 & u_{1}=(p_{2}u_{2})/p_{4},\,u_{2}=-(p_{3}^{2}+t_{1}p_{4}^{2})/p_{4},\,t_{2}=-(p_{3}u_{2})/p_{4}^{2},\\
 & s_{1}=-(t_{136}p_{2}^{2}+2t_{245}p_{2}p_{4}+t_{246}p_{4}^{2})/p_{2},\,s_{2}=(p_{3}/p_{4})s_{1},\,s_{3}=-((2p_{3}^{2}+t_{1}p_{4}^{2})/p_{4}^{2})s_{1},\\
 & t_{123}=-(t_{126}p_{2}p_{3}p_{4}^{2}+t_{136}(3p_{3}^{2}+2t_{1}p_{4}^{2})p_{2}^{2}+3t_{245}(p_{3}^{2}+t_{1}p_{4}^{2})p_{2}p_{4}+t_{246}t_{1}p_{4}^{4})/(p_{2}^{2}p_{4}^{2}),\\
 & t_{124}=(t_{126}p_{2}p_{3}p_{4}-t_{245}(3p_{3}^{2}+t_{1}p_{4}^{2})p_{2}-t_{246}(3p_{3}^{2}+t_{1}p_{4}^{2})p_{4})/(p_{2}p_{4}^{2}),\\
 & t_{125}=(-t_{126}p_{2}p_{4}^{2}+3t_{136}p_{2}^{2}p_{3}+6t_{245}p_{2}p_{3}p_{4}+3t_{246}p_{3}p_{4}^{2})/(p_{2}^{2}p_{4}),\\
 & t_{135}=-(p_{4}(3t_{136}p_{2}^{2}+3t_{245}p_{2}p_{4}+t_{246}p_{4}^{2}))/p_{2}^{3}\},
\end{align*}
}and $\Pi_{\mathbb{A}}^{15}\setminus\mA_{T}$ has $c({\rm G}(2,5))$-singularities
along $\Delta$, where we call a singularity isomorphic to the vertex
of the cone over ${\rm G}(2,5)$ a $c({\rm G}(2,5))$-singularity. 
\end{prop}

\begin{proof}
Assume that $p_{2}p_{4}\not=0$. Then, we can solve the equation $F_{2}=0$
with respect to $s_{3}$, and hence we may erase the coordinate $s_{3}$.
Now we consider the following coordinate change: the 8 coordinates
\begin{equation}
p_{2},\,p_{3},\,p_{4},\,t_{1},\,t_{126},\,t_{136},\,t_{245},\,t_{246}\label{eq:8coord}
\end{equation}
 are unchanged while the 10 coordinates 
\begin{align}
p_{1},u_{1},u_{2},s_{1},s_{2},t_{2},t_{123},t_{124},t_{125},t_{135}\label{eq:10coord}
\end{align}
 are changed to

{\small{}
\begin{align*}
 & m_{12}=\nicefrac{1}{p_{2}p_{4}^{3}}\,(t_{123}p_{2}p_{4}^{2}+t_{126}p_{3}p_{4}^{2}+p_{3}p_{4}s_{2}+t_{135}p_{2}p_{4}u_{1}\\
 & \qquad+t_{136}(p_{2}p_{3}^{2}+p_{1}p_{3}p_{4}+p_{2}p_{4}u_{2})+2t_{245}p_{3}^{2}p_{4}),\\
 & m_{13}=-\nicefrac{1}{p_{2}^{2}p_{4}^{2}}\,(t_{135}(p_{1}^{2}+t_{1}p_{2}^{2})p_{2}^{2}-t_{125}(p_{2}p_{3}-p_{1}p_{4})p_{2}^{2}-t_{126}p_{2}^{2}p_{3}p_{4}\\
 & \qquad+t_{246}(p_{2}p_{3}+p_{1}p_{4})p_{3}p_{4}+t_{124}p_{2}^{2}p_{4}^{2}+p_{2}(-p_{2}p_{3}+p_{1}p_{4})s_{2}+t_{135}p_{2}^{2}u_{1}\\
 & \qquad+p_{2}^{2}s_{1}u_{2}+t_{136}(2p_{1}p_{3}+2t_{1}p_{2}p_{4}+2p_{4}u_{1}+p_{2}u_{2})p_{2}^{2}\\
 & \qquad+t_{245}(p_{2}p_{3}^{2}+2p_{1}p_{3}p_{4}+t_{1}p_{2}p_{4}^{2}+2p_{2}p_{4}u_{2})p_{2}),\\
 & m_{14}=\nicefrac{1}{p_{2}^{2}p_{4}}\,(p_{2}p_{4}s_{2}+t_{136}p_{2}^{2}p_{3}+2t_{245}p_{2}p_{3}p_{4}+t_{246}p_{3}p_{4}^{2}),\\
 & m_{15}=-(p_{3}u_{1}+p_{2}p_{4}t_{2})/p_{4},\\
 & m_{23}=\nicefrac{1}{p_{2}p_{4}^{2}}\,(-(p_{2}p_{3}+p_{1}p_{4})s_{1}-p_{2}p_{4}s_{2}-t_{125}p_{2}^{2}p_{4}-t_{126}p_{2}p_{4}^{2}\\
 & \qquad+t_{136}(-p_{2}p_{3}+p_{1}p_{4})p_{2}-t_{245}(p_{2}p_{3}+p_{1}p_{4})p_{4}),\\
 & m_{24}=\nicefrac{1}{p_{2}p_{4}}\,(-t_{135}p_{2}^{2}+p_{4}s_{1}-2t_{136}p_{2}p_{4}-t_{245}p_{4}^{2}),\\
 & m_{25}=-\nicefrac{p_{2}}{p_{4}^{2}}\,U_{3},\\
 & m_{34}=(p_{2}s_{1}+t_{136}p_{2}^{2}+2t_{245}p_{2}p_{4}+t_{246}p_{4}^{2})/p_{2},\\
 & m_{35}=-U_{1}\\
 & m_{45}=-p_{2}p_{3}+p_{1}p_{4}.
\end{align*}
}Indeed, these define a coordinate change since these equations can
be solve regularly with respect to the 10 coordinates (\ref{eq:10coord}).
Then, by direct computations for this coordinate change, we see that
$\Pi_{\mathbb{A}}^{15}\cap\{p_{2}p_{4}\not=0\}$ is isomorphic to
the variety defined by the five $4\times4$ Pfaffians of the skew-symmetric
matrix {\small{}
\[
\left(\begin{array}{ccccc}
0 & m_{12} & m_{13} & m_{14} & m_{15}\\
 & 0 & m_{23} & m_{24} & m_{25}\\
 &  & 0 & m_{34} & m_{35}\\
 &  &  & 0 & m_{45}\\
 &  &  &  & 0
\end{array}\right)
\]
}in the open subset $\{p_{2}p_{4}\not=0\}$ of the affine space $\mathbb{A}^{18}$
with the 10 coordinates $m_{ij}\,(1\leq i<j\leq5),$ and 8 coordinates
(\ref{eq:8coord}). Therefore $\Pi_{\mathbb{A}}^{15}\cap\{p_{2}p_{4}\not=0\}$
is an irreducible $15$-dimensional variety and is singular along
the locus corresponding to $\{m_{ij}=0\,(1\leq i<j\leq5)\}$ and has
$c({\rm G}(2,5))$-singularities there. We may verify $\{m_{ij}=0\,(1\leq i<j\leq5)\}\cap\{p_{2}p_{4}\not=0\}=S\cap\{p_{2}p_{4}\not=0\}$
as in the statement of the proposition.

Since points with $p_{2}\not=0$ or $p_{4}\not=0$ can be transformed
by the ${\rm GL}_{2}$-action on $\Pi_{\mathbb{A}}^{15}$ to points
with $p_{2}p_{4}\not=0$ ((\ref{eq:GL2act}) in Proposition \ref{prop:GL}),
the open subset $\Pi_{\mathbb{A}}^{15}\cap\{p_{2}\not=0,\ \text{or}\ p_{4}\not=0\}$
is an irreducible $15$-dimensional variety and $\Pi_{\mA}^{15}$
has $c({\rm G}(2,5))$-singularities along an $8$-dimensional locus
in there. 

We show that $\Pi_{\mA}^{15}\cap\{U_{3}\not=0\}$ is an irreducible
$15$-dimensional variety and is smooth. If $U_{3}\not=0,$ then we
may eliminate $t_{124},t_{245},t{}_{246}$ from the equations of $\Pi_{\mathbb{A}}^{15}$.
Therefore $\Pi_{\mathbb{A}}^{15}$ is isomorphic to a hypersurface
on the open subset $\{U_{3}\not=0\}$. In particular, $\Pi_{\mA}^{15}\cap\{U_{3}\not=0\}$
is purely $15$-dimensional. By the Jacobian criterion, if the hypersurface
is singular, then $p_{2}=0$ or $p_{4}=0$ implies that $U_{3}=0$,
a contradiction. We see that $\Pi_{\mathbb{A}}^{15}\cap\{U_{3}\not=0\}\cap\{p_{2}p_{4}\not=0\}$
is smooth since this is equal to $\Pi_{\mathbb{A}}^{15}\cap\{p_{2}p_{4}\not=0\}\cap\{m_{25}\not=0\}$
in the above consideration. Therefore $\Pi_{\mathbb{A}}^{15}$ is
smooth on the open subset $\{U_{3}\not=0\}.$ Note that there exists
some positive weights of coordinates for $\Pi_{\mA}^{15}$ such that
all the equations of $\Pi_{\mA}^{15}$ are weighted homogeneous (cf.~Subsection
\ref{subsec:Weights-for-variablesPi}). Then the hypersurface is also
weighted homogeneous since so is $U_{3}$. Therefore $\Pi_{\mA}^{15}\cap\{U_{3}\not=0\}$
is connected. Therefore, we have shown that $\Pi_{\mA}^{15}\cap\{U_{3}\not=0\}$
is an irreducible $15$-dimensional variety and is smooth.

By the ${\rm GL}_{2}$-symmetry of $\Pi_{\mathbb{A}}^{15}$, the result
in the previous paragraph implies that $\Pi_{\mA}^{15}\cap\{U_{1}\not=0\}$
is an irreducible $15$-dimensional variety and is smooth.

Since $\Pi_{\mA}^{15}\cap\{U_{1}\not=0\}$, $\Pi_{\mA}^{15}\cap\{U_{3}\not=0\}$
and $\Pi_{\mA}^{15}\cap\{p_{2}\not=0\ \text{or}\ p_{4}\not=0\}$ mutually
intersect nontrivially, and are irreducible and $15$-dimensional,
$\Pi_{\mA}^{15}\cap\{p_{2}\not=0\ \text{or}\ p_{4}\not=0\ \text{or}\ U_{1}\not=0\ \text{or}\ U_{3}\not=0\}$
is an irreducible $15$-dimensional variety. 

Finally, we investigate $\Pi_{\mathbb{A}}^{15}$ along the closed
subset 

\[
\Gamma:=\{U_{1}=U_{3}=p_{2}=p_{4}=0\}\cap\Pi_{\mathbb{A}}^{15}.
\]
It holds that 
\[
\Gamma=(\Pi_{\mathbb{A}}^{15}\cap\{p_{1}=p_{2}=p_{3}=p_{4}=0\})\cap(\{u_{1}=u_{2}=0\}\cup\{s_{1}=s_{2}=0\}).
\]
We can calculate the Zariski tangent space at any point of these loci
outside $\mathbb{A}_{\Pi}$, and then can show that $\Pi_{\mathbb{A}}^{15}$
is $15$-dimensional and is smooth there. For example, let $\mathsf{p}$
be a point of $(\Pi_{\mathbb{A}}^{15}\cap\{p_{1}=p_{2}=p_{3}=p_{4}=u_{1}=u_{2}=0\})\setminus\mathbb{A}_{\Pi}^{10}$.
Then we have one of $s_{1},s_{2},s_{3}\not=0$. If $s_{1}=0$, then
the equality $F_{8}=0$ implies that $s_{2}=0$, and, then the equality
$F_{9}=0$ implies that $s_{3}=0,$ a contradiction. Therefore, we
have $s_{1}\not=0$. By the equality $F_{7}=F_{8}=F_{9}=0$, we obtain
$s_{3}=\nicefrac{s_{2}^{2}}{s_{1}},t_{2}=\nicefrac{s_{2}(s_{2}^{2}+s_{1}^{2}t_{1})}{s_{1}^{3}}.$
Hence we have
\begin{align*}
(\Pi_{\mathbb{A}}^{15}\cap\{p_{1} & =p_{2}=p_{3}=p_{4}=u_{1}=u_{2}=0\})\setminus\mathbb{A}_{\Pi}\\
 & =\{p_{1}=p_{2}=p_{3}=p_{4}=u_{1}=u_{2}=0,s_{1}\not=0,s_{3}=\nicefrac{s_{2}^{2}}{s_{1}},t_{2}=\nicefrac{s_{2}(s_{2}^{2}+s_{1}^{2}t_{1})}{s_{1}^{3}}\}.
\end{align*}
Using this descriptions, we can easily calculate the Zariski tangent
space of $\Pi_{\mathbb{A}}^{15}$ at $\mathsf{p}$ using the polynomials
$F_{1}$--$F_{9}$.

Since $\dim\Gamma<15$, we conclude that $\Pi_{\mA}^{15}$ is the
union of the closed subset $\mA_{T}$ and an irreducible $15$-dimensional
variety with $c({\rm G}(2,5))$-singularities. 
\end{proof}

\subsection{Gorensteinness and graded $9\times16$ minimal free resolution of
the ideal of $\Pi_{\mA}^{15}$}

we set 
\[
\delta:=3\left(w(p_{3})-w(u_{1})+2w(u_{2})\right)+4w(t_{246}).
\]
Then we observe that the summation of all the weights of the equations
of $\Pi_{\mathbb{A}}^{15}$ is equal to

\[
\sum_{i=1}^{9}d_{i}=3\delta.
\]

\begin{prop}
\label{prop:916}Let $S_{\Pi}$ be the polynomial ring over $\mC$
whose variables are (\ref{eq:coordpi}). Let $I_{\Pi}$ be the ideal
of the polynomial ring $S_{\Pi}$ generated by the $9$ equations
$F_{1},\dots,F_{9}$ of $\Pi_{\mathbb{A}}^{15}$.

\begin{enumerate}[$(1)$]

\item We give nonnegative weights for coordinates of $S_{\Pi}$ such
that all the equations of $\Pi_{\mA}^{15}$ are weighted homogeneous,
and we denote by $w(*)$ the weight of the monomial $*$. 

\begin{enumerate}[$({1}\text{-}1)$]

\item We denote by $\mP$ the corresponding weighted projective space,
where we allow some coordinates being nonzero constants (thus $\dim\mP$
could be less than $18$). It holds that
\[
\omega_{\mP}=\sO_{\mP}(12w(u_{1})-23w(u_{2})+2w(p_{3})-11w(t_{246})).
\]

\item We set
\[
\delta:=3\left(w(p_{3})-w(u_{1})+2w(u_{2})\right)+4w(t_{246}).
\]
The ideal of $I_{\Pi}$ of $\Pi_{\mathbb{A}}^{15}$ has the following
graded minimal $S_{\Pi}$-free resolution: 
\begin{equation}
0\leftarrow R_{\Pi}\leftarrow P_{0}\leftarrow P_{1}\leftarrow P_{2}\leftarrow P_{3}\leftarrow P_{4}\leftarrow0,\label{eq:minresPi}
\end{equation}
where {\small{}
\begin{align*}
 & P_{0}=S_{\Pi},P_{1}=\oplus_{i=1}^{9}S_{\Pi}(-d_{i}),\\
 & P_{2}=S_{\Pi}(-(d_{6}+w(p_{2})))\oplus S_{\Pi}(-(d_{6}+w(p_{1})))\oplus S_{\Pi}(-(d_{8}+w(p_{2})))\oplus S_{\Pi}(-(d_{8}+w(p_{4})))\\
 & \oplus S_{\Pi}(-(d_{8}+w(p_{1})))^{\oplus2}\oplus S_{\Pi}(-(d_{8}+w(p_{3})))^{\oplus2}\\
 & \oplus S_{\Pi}(-\delta+(d_{8}+w(p_{1})))^{\oplus2}\oplus S_{\Pi}(-\delta+(d_{8}+w(p_{3})))^{\oplus2}\\
 & \oplus S_{\Pi}(-\delta+(d_{6}+w(p_{2})))\oplus S_{\Pi}(-\delta+(d_{6}+w(p_{1})))\oplus S_{\Pi}(-\delta+(d_{8}+w(p_{2})))\oplus S_{\Pi}(-\delta+(d_{8}+w(p_{4}))),\\
 & P_{3}=\oplus_{i=1}^{9}S_{\Pi}(-\delta+d_{i}),P_{4}=S_{\Pi}(-\delta).
\end{align*}
}{\small\par}

\end{enumerate}

\item $I_{\Pi}$ is a Gorenstein ideal of codimension $4$. 

\item $\Pi_{\mA}^{15}$ is irreducible and reduced, thus $I_{\Pi}$
is a prime ideal.

\item $\Pi_{\mA}^{15}$ is normal.

\end{enumerate}
\end{prop}

\begin{proof}
We may compute the $S_{\Pi}$-free resolution (\ref{eq:minresPi})
of $I_{\Pi}$ by \textsc{Singular} \cite{DGPS}. For the remaining
assertions, the proof of \cite[Prop.4.8]{Tak5} works verbatim.
\end{proof}

\subsection{Factoriality of $\Pi_{\mathbb{A}}^{15}$\label{sec:Factoriality-of-the}}

In this section, we show that the affine coordinate ring $R_{\Pi}$
of $\Pi_{\mathbb{A}}^{15}$ is a UFD. We use the following polynomial
$\mathsf{b}$ frequently: 
\begin{equation}
\mathsf{b}:=p_{3}^{2}+t_{1}p_{4}^{2}+p_{4}u_{2},\label{eq:boundary}
\end{equation}
which nothing but $U_{3}$ as in (\ref{eq:LUIA}). We denote by $\bar{\mathsf{b}}$
the image of $\mathsf{b}$ in $R_{\Pi}$. 
\begin{lem}
\label{lem:primePi}The element $\bar{\mathsf{b}}\in R_{\Pi}$ is
a prime element, equivalently, $\Pi_{\mathbb{A}}^{15}\cap\{\mathsf{b}=0\}$
is irreducible and reduced. 
\end{lem}

\begin{proof}
By the proof of Lemma \ref{lem:prime-Up}, it suffices to show that
$\Pi_{\mathbb{A}}^{15}\cap\{\mathsf{b}=0\}$ is normal. Since it is
Gorenstein by Proposition \ref{prop:916} (2) and $\bar{\mathsf{b}}$
is not a zero divisor by ibid.~(3), we have only to show that $\Pi_{\mathbb{A}}^{15}\cap\{\mathsf{b}=0\}$
is regular in codimension 1. Assume that $p_{1}p_{4}\not=0$. Then,
by the proof of Proposition \ref{prop:SingPi}, $\Pi_{\mathbb{A}}^{15}\cap\{\mathsf{b}=0\}$
is isomorphic to the divisor $\{m_{25}=0\}$. Therefore it is easy
to verify that it is regular in codimension 1. Then we have only to
check that ${\rm codim_{\Pi_{\mA}^{15}}}\left(\{p_{i}=0\}\cap{\rm Sing}(\Pi_{\mathbb{A}}^{15}\cap\{\mathsf{b}=0\})\right)\geq7$
for $i=2,4$. In a similar way to the proof of Proposition \ref{prop:SingPi},
we can show that $\Pi_{\mathbb{A}}^{15}$ is isomorphic to a hypersurface
on the open set $\Pi_{\mathbb{A}}^{15}\cap\{U_{1}\not=0\}$. Using
this description, we may verify that $(\Pi_{\mathbb{A}}^{15}\cap\{\mathsf{b}=0\})$
is regular in codimension 1 on the open set $\{U_{1}\not=0\}$ by
the Jacobian criterion. It remains to show that ${\rm codim_{\Pi_{\mA}^{15}}}\left(\{p_{i}=0,U_{1}=0\}\cap{\rm Sing}(\Pi_{\mathbb{A}}^{15}\cap\{\mathsf{b}=0\})\right)\geq7$
for $i=2,4$. Actually, we can verify ${\rm codim_{\Pi_{\mA}^{15}}}\left(\{p_{i}=0,U_{1}=0,\mathsf{b}=0\}\cap\Pi_{\mathbb{A}}^{15}\right)\geq7$
for $i=2,4$ by a straightforward calculations. 
\end{proof}
In the proof of Lemma \ref{lem:primePi}, we have proved the following: 
\begin{cor}
\label{cor:boundarynormal.}$\Pi_{\mathbb{A}}^{15}\cap\{\mathsf{b}=0\}$
is normal. 
\end{cor}

\begin{prop}
\label{prop:UFDPi}The affine coordinate ring $R_{\Pi}$ of $\Pi_{\mathbb{A}}^{15}$
is a UFD. 
\end{prop}

\begin{proof}
As we saw in Subsection \ref{subsec:Weights-for-variablesPi}, $R_{\Pi}$
can be positively graded. In this situation, by \cite[Prop.~7.4]{UFD},
$R_{\Pi}$ is UFD if and only if so is the localization $R_{\Pi,o}$
of $R_{\Pi}$ with respect to the maximal irrelevant ideal. Thus we
have only to prove that $R_{\Pi,o}$ is UFD. Abusing notation, we
denote $R_{\Pi,o}$ by $R_{\Pi}$. We already know that $R_{\Pi}$
is a domain by Proposition \ref{prop:916} (3). Then, by Nagata's
theorem \cite[Thm.~20.2]{Matsumura}, it suffices to show that $\bar{\mathsf{b}}$
is a prime element, and the localization $(R_{\Pi})_{\bar{b}}$ is
a UFD. We already show the first assertion in Lemma \ref{lem:primePi}.
As we see in the proof of Proposition \ref{prop:SingPi}, the open
subset $\Pi_{\mathbb{A}}^{15}\cap\{\mathsf{b}\not=0\}$ is smooth.
Therefore $(R_{\Pi})_{\bar{\mathsf{b}}}$ is a complete intersection
(a hypersurface) and its localizations with respect to all the prime
ideals are regular. This implies that $(R_{\Pi})_{\bar{\mathsf{b}}}$
is a UFD by \cite[Exp. XI, Cor.~3.10 and Thm.~3.13 (ii)]{sga2}. 
\end{proof}
In \cite{Tak3}, we construct prime $\mQ$-Fano 3-folds with the following
subvariety of $\Pi_{\mA}^{15}$ as the key variety:

\[
\Pi_{\mA}^{14}:=\Pi_{\mA}^{15}\cap\{t_{246}=1\}.
\]

To obtain information on $\Pi_{\mathbb{A}}^{14}$ from $\Pi_{\mA}^{15}$,
we use the following: 
\begin{prop}
\label{prop:Pi14}There exists an isomorphism from $\Pi_{\mathbb{A}}^{15}\cap\{t_{246}\not=0\}$
to $\Pi_{\mathbb{A}}^{14}\times(\mathbb{A}^{1})^{*}$. Moreover, this
induces an isomorphism from $\Pi_{\mathbb{A}}^{15}\cap\{t_{246}\not=0,\mathsf{b}\not=0\}$
to $(\Pi_{\mathbb{A}}^{14}\cap\{\mathsf{b}\not=0\})\times(\mathbb{A}^{1})^{*}$. 
\end{prop}

\begin{proof}
We can define a morphism $\iota_{1}\colon\Pi_{\mathbb{A}}^{15}\cap\{t_{246}\not=0\}\to\Pi_{\mathbb{A}}^{14}$
by the following coordinate change: 
\[
\alpha\mapsto\begin{cases}
t_{246}^{-2}\alpha & :\alpha=s_{1},s_{2},s_{3}\\
\alpha & :\alpha=t_{1},t_{2}\\
t_{246}^{-1}\alpha & :\text{{otherwise}}.
\end{cases}
\]
Let $\iota_{2}\colon\Pi_{\mathbb{A}}^{15}\cap\{t_{246}\not=0\}\to(\mathbb{A}^{1})^{*}$
be the projection to $(\mathbb{A}^{1})^{*}$ with $t_{246}$ as coordinate.
We can show easily that the morphism $\iota_{1}\times\iota_{2}$ is
an isomorphism and the latter assertion holds. 
\end{proof}
\begin{prop}
\label{prop:Pic1}Let $\Pi_{\mathbb{P}}^{14}$ be the weighted projectivization
of $\Pi_{\mathbb{\mathbb{A}}}^{15}$ with some positive weights of
coordinates, and $\Pi_{\mathbb{P}}^{14}$ the weighted cone over $\Pi_{\mathbb{P}}^{13}$.
It holds that

\vspace{3pt}

\noindent $(1)$ any prime Weil divisors on $\Pi_{\mathbb{P}}^{13}$
and $\Pi_{\mathbb{P}}^{14}$ are the intersections between $\Pi_{\mathbb{P}}^{13}$
and $\Pi_{\mathbb{P}}^{14}$ respectively, and weighted hypersurfaces.
In particular, $\Pi_{\mathbb{P}}^{13}$ and $\Pi_{\mathbb{P}}^{14}$
are $\mathbb{Q}$-factorial and have Picard number one, and

\noindent $(2)$ let $X$ be a quasi-smooth threefold such that $X$
is a codimension $10$ or $11$ weighted complete intersection in
$\Pi_{\mathbb{P}}^{13}$ or $\Pi_{\mathbb{P}}^{14}$ respectively.
Assume moreover that $X\cap\{\mathsf{b}=0\}$ is a prime divisor.
Then any prime Weil divisor on $X$ is the intersection between $X$
and a weighted hypersurface. In particular, $X$ is $\mathbb{Q}$-factorial
and has Picard number one. 
\end{prop}

\begin{proof}
The proofs of \cite[Cor.4.10]{Tak5} work verbatim using Propositions
\ref{prop:UFDPi} and \ref{prop:Pi14}. 
\end{proof}

\subsection{Relation between $\Pi_{\mathbb{A}}^{15}$ and $\mathscr{H}_{\mathbb{A}}^{13}$. }

We establish the following relation between $\Pi_{\mathbb{A}}^{15}$
and $\mathscr{H}_{\mathbb{A}}^{13}$. We define some notation. Let
$\mathbb{A}_{T}$ be the affine $10$-space as in the section \ref{sec:Singular-locus-Pi}
(but we do not consider here that $\mathbb{A}_{T}\subset\Pi_{\mathbb{A}}^{15}$).
Let $\mathbb{A}'_{T}$ be the affine space with the coordinates $t_{123}$,
$t_{124}$, $t_{125}$, $t_{126}$, $t_{135}$, $t_{136}$, $t_{245}$,
$t_{246}$ and $\mathbb{A}_{A\!B}$ the affine space with the coordinates
$A$, $B$. We set $\widetilde{\mathbb{A}}_{T}=\mathbb{A}_{A\!B}\times\mathbb{A}'_{T}$.
Let $b\colon\widetilde{\mathbb{A}}_{T}\to\mathbb{A}_{T}$ be the morphism
defined with 
\[
t_{1}=-A^{2}-AB-B^{2},\,t_{2}=-AB(A+B)
\]
and the remaining coordinates being unchanged. The morphism $b$ is
a finite morphism of degree six since $b$ is identified with the
morphism $\{A+B+C=0\}\times\mathbb{\mA}'_{T}\to\mathbb{A}_{T}$ defined
with $t_{1}=AB+BC+CA$, $t_{2}=ABC$ and the remaining coordinates
being unchanged, where $\{A+B+C=0\}$ is a closed subset of the affine
$3$-space with the coordinates $A,B,C$. Let $\widetilde{\Pi}_{\mathbb{A}}^{15}:=\Pi_{\mathbb{A}}^{15}\times_{\mathbb{A}_{T}}\widetilde{\mathbb{A}}_{T}$,
and $b_{\Pi}\colon\widetilde{\Pi}_{\mathbb{A}}^{15}\to\Pi_{\mathbb{A}}^{15}$
and $p_{\Pi}\colon\widetilde{\Pi}_{\mathbb{A}}^{15}\to\widetilde{\mathbb{A}}_{T}$
the naturally induced morphisms. Let $p_{\mathscr{H}}\colon\mathscr{H}_{\mathbb{A}}^{13}\to\mathbb{A_{\mathsf{P}}}$
be the natural projection. Let $U_{A\!B}$ be the open subset $\{(A-B)(2A+B)(A+2B)\not=0\}$
of $\mathbb{A}_{A\!B}$. Note that the morphism $b$ is unramified
on $U_{A\!B}\times\mathbb{A}_{T}\subset\widetilde{\mathbb{A}}_{T}$. 
\begin{prop}
\label{prop:isomHPi}There exists an isomorphism from $U_{A\!B}\times\mathbb{A}_{T}$
to $U_{A\!B}\times\mathbb{A}_{\mathsf{P}}$, and an isomorphism from
$\widetilde{\Pi}_{\mathbb{A}}^{15}\cap\{(A-B)(2A+B)(A+2B)\not=0\}$
to $U_{A\!B}\times\mathscr{H}_{\mathbb{A}}^{13}$ which fit into the
following commutative diagram:{\small{}
\[
\xymatrix{\widetilde{\Pi}_{\mathbb{A}}^{15}\cap\{(A-B)(2A+B)(A+2B)\not=0\}\ar[r]\ar[d]_{p_{\Pi}} & U_{A\!B}\times\mathscr{H}_{\mathbb{A}}^{13}\ar[d]^{p_{\mathscr{H}}\times{\rm id}}\\
U_{A\!B}\times\mathbb{A}'_{T}\ar[r] & U_{A\!B}\times\mathbb{A}_{\mathsf{P}}
}
\]
}{\small\par}
\end{prop}

\begin{proof}
We can directly check that a morphism from $U_{A\!B}\times\mathbb{A}'_{T}$
to $U_{A\!B}\times\mathbb{A}_{\mathsf{P}}$, and a morphism from $\widetilde{\Pi}_{\mathbb{A}}^{15}\cap\{(A-B)(2A+B)(A+2B)\not=0\}$
to $U_{A\!B}\times\mathscr{H}_{\mathbb{A}}^{13}$ can be defined by
the following equalities and they are actually isomorphisms with the
desired properties, where we set $C:=-A-B$ to write them symmetrically.

{\small{}
\begin{align*}
p_{111} & =-3t_{246},\\
p_{112} & =3t_{245}+\nicefrac{3(t_{124}-Ct_{126})}{(A-C)(B-C)},\\
p_{121} & =3t_{245}+\nicefrac{3(-t_{124}+Bt_{126})}{(A-B)(B-C)},\\
p_{122} & =-3t_{136}+\nicefrac{3(t_{123}-At_{125})}{(A-B)(A-C)},\\
p_{211} & =3t_{245}+\nicefrac{3(t_{124}-At_{126})}{(A-B)(A-C)},\\
p_{212} & =-3t_{136}+\nicefrac{3(-t_{123}+Bt_{125})}{(A-B)(B-C)},\\
p_{221} & =-3t_{136}+\nicefrac{3(t_{123}-Ct_{125})}{(A-C)(B-C)},\\
p_{222} & =3t_{135},
\end{align*}
\begin{align*}
x_{11} & =-(Ap_{1}+A^{2}p_{2}-u_{1}),\,x_{21}=-(Ap_{3}+A^{2}p_{4}-u_{2}),\\
x_{12} & =Bp_{1}+B^{2}p_{2}-u_{1},\,x_{22}=Bp_{3}+B^{2}p_{4}-u_{2},\\
x_{13} & =Cp_{1}+C^{2}p_{2}-u_{1},\,x_{23}=Cp_{3}+C^{2}p_{4}-u_{2},
\end{align*}
}{\small\par}

{\small{}
\begin{align*}
u_{1} & =-3\left(BCs_{1}+As_{2}+s_{3}\right)+\nicefrac{3}{(A-B)(A-C)}[\left(2At_{123}+BCt_{125}\right)p_{1}\\
 & -\left((A^{2}+B^{2}+C^{2}+BC)t_{123}+2ABCt_{125}\right)p_{2}+\left(2At_{124}+BCt_{126}\right)p_{3}\\
 & -((A^{2}+B^{2}+C^{2}+BC)t_{124}+2ABCt_{126})p_{4}+(t_{123}-At_{125})u_{1}+(t_{124}-At_{126})u_{2}],
\end{align*}
\begin{align*}
u_{2} & =-3(ACs_{1}+Bs_{2}+s_{3})-\nicefrac{3}{(A-B)(B-C)}[\left(2BL_{123}+ACt_{125}\right)p_{1}\\
 & -\left((A^{2}+B^{2}+C^{2}+AC)t_{123}+2ABCt_{125}\right)p_{2}+\left(2Bt_{124}+ACt_{126}\right)p_{3}\\
 & -\left((A^{2}+B^{2}+C^{2}+AC)t_{124}+2ABCt_{126}\right)p_{4}+(t_{123}-Bt_{125})u_{1}+(t_{124}-Bt_{126})u_{2}],\\
u_{3} & =3(ABs_{1}+Cs_{2}+s_{3})-\nicefrac{3}{(A-C)(B-C)}[\left(2Ct_{123}+ABt_{125}\right)p_{1}\\
 & -\left((A^{2}+B^{2}+C^{2}+AB)t_{123}+2ABCt_{125}\right)p_{2}+\left(2Ct_{124}+ABt_{126}\right)p_{3}\\
 & -\left((A^{2}+B^{2}+C^{2}+AB)t_{124}+2ABCt_{126}\right)p_{4}+\left(t_{123}-Ct_{125}\right)u_{1}+\left(t_{124}-Ct_{126}\right)u_{2})].
\end{align*}
}{\small\par}
\end{proof}

\subsection{$\mathbb{P}^{2}\times\mathbb{P}^{2}$-fibration associated to $\Pi_{\mathbb{A}}^{15}$\label{subsec:-fibration-associated-toPi}}

Note that any equation of $\text{\ensuremath{\Pi_{\mathbb{A}}^{15}}}$
is of degree two if we regard the coordinates of $\mathbb{A}_{T}$
as constants. Therefore, considering the coordinates $p_{1},p_{2},p_{3},p_{4},u_{1},u_{2},s_{1},s_{2},s_{3}$
as projective coordinates, we obtain a quasi-projective variety with
the same equations as $\text{\ensuremath{\Pi_{\mathbb{A}}^{15}}}$,
which we denote by $\widehat{\Pi}$. We also denote by $\rho_{\Pi}\colon\widehat{\Pi}\to\mathbb{A}_{T}$
the natural projection.

We note that the image of $U_{A\!B}\times\mathbb{A}_{T}'$ by the
morphism $b$ is the open subset $\{4t_{1}^{3}+27t_{2}^{2}\not=0\}\times\mathbb{A}_{T}'$
of $\mathbb{A}_{T}$. Let $\{D'_{\mathscr{H}}=0\}\subset\mA_{T}$
be the image by $b$ of the pull-back on $U_{AB}\times\mathbb{A}_{T}'$
of $U_{A\!B}\times\{D_{\mathscr{H}}=0\}$. 
\begin{prop}
\label{prop:P2P2Pi}Let $\rho_{\Pi}\colon\widehat{\Pi}\to\mathbb{A}_{T}$
be the natural projection. The $\rho_{\Pi}$-fibers over points of
$\{D'_{\mathscr{H}}\not=0\}\cap(\{4t_{1}^{3}+27t_{2}^{2}\not=0\}\times\mathbb{A}'_{T})$
are isomorphic to $\mathbb{P}^{2}\times\mathbb{P}^{2}$, and the $\rho_{\Pi}$-fibers
over points of a nonempty open subset of $\{D'_{\mathscr{H}}=0\}$
are isomorphic to $\mathbb{P}^{2,2}$. 
\end{prop}

\begin{proof}
The proof of Proposition \ref{prop:P2P2Up} works verbatim. 
\end{proof}
It is possible to describe all the $\rho_{\Pi}$-fibers by using the
${\rm GL}_{2}$-action on $\Pi_{\mA}^{15}$ but it needs more works
and pages; Proposition \ref{prop:P2P2Pi} is sufficient for an application
in this paper (cf.~Proposition \ref{prop:MoreSing}). 

\section{\textbf{Singularities of $\text{\ensuremath{\Upsilon_{\mA}^{14}} and }\Pi_{\mathbb{A}}^{15}$}\label{sec:Singularities-ofUpPi}}
\begin{prop}
\label{prop:MoreSing} Let $\Sigma$ be one of the varieties $\Upsilon_{\mA}^{14}$
and $\Pi_{\mathbb{A}}^{15}$, which is defined over $\mC$. The variety
$\Sigma$ has only terminal singularities with the following descriptions:

\begin{enumerate}[$(1)$]

\item $\Sigma$ has $c({\rm G}(2,5))$-singularities along a locally
closed subset $\Delta$ of codimension $7$ in $\Sigma$.

\item The singular locus of $\Sigma$ coincides with $\Delta\cup\mA_{{\rm B}}$,
where $\mA_{{\rm B}}=\mA_{S\!T}$ for $\Sigma=\Upsilon_{\mA}^{14}$
or $\mA_{{\rm B}}=\mA_{S}$ for $\Sigma=\Pi_{\mA}^{15}$. 

\item There exists a primitive $K$-negative divisorial extraction
$f\colon\widetilde{\Sigma}\to\Sigma$ such that

\begin{enumerate}[$(a)$] 

\item singularities of $\widetilde{\Sigma}$ are only $c({\rm G}(2,5))$-singularities
along the strict transform of the closure of $\Delta$, and

\item for the $f$-exceptional divisor $E_{\Sigma}$, the morphism
$f|_{E_{\Sigma}}$ can be identified with $\rho_{\Upsilon}$ if $\Sigma=\Upsilon_{\mA}^{14}$
or $\rho_{\Pi}$ if $\Sigma=\Pi_{\mA}^{15}$. 

\end{enumerate}

\end{enumerate}
\end{prop}

\begin{proof}
The proof of \cite[Prop.4.11]{Tak5} works verbatim using the following
assertions:

\noindent $\Sigma=\Upsilon_{\mA}^{14}$: Theorem \ref{thm:PropUp},
Proposition \ref{prop:The-singular-locusUp}, Corollary \ref{cor:The-variety-normal UP},
and Propositions \ref{prop:UFD-Up} and \ref{prop:P2P2Up}.

\noindent $\Sigma=\Pi_{\mA}^{15}$: Propositions \ref{prop:SingPi},
\ref{prop:916}, \ref{prop:UFDPi} and \ref{prop:P2P2Pi}.
\end{proof}

\end{document}